\documentclass[11pt]{amsart}

\usepackage{latexsym,amsfonts,amssymb,amscd,graphics,axodraw,verbatim}
\usepackage[all]{xy}
 
\theoremstyle{plain}
\newtheorem{theorem}{Theorem}[section]
\newtheorem{lemma}[theorem]{Lemma}
\newtheorem{prop}[theorem]{Proposition}
\newtheorem{conj}[theorem]{Conjecture}
 
\theoremstyle{definition}
\newtheorem{definition}[theorem]{Definition}
\newtheorem{rem}[theorem]{Remark}
\newtheorem{ex}[theorem]{Example}

\numberwithin{equation}{section}

\newcommand{\gggg}{\mathfrak{g}}
\newcommand{\wt}{\mathrm{wt}}
\newcommand{\inner}[2]{\langle\,#1\,,\,#2\,\rangle}
\newcommand{\BB}{\mathbb{B}}
\newcommand{\VV}{\mathbb{V}}
\newcommand{\cl}{\mathrm{cl}}
\newcommand{\af}{\mathrm{af}}
\newcommand{\Pfin}{P_{\mathrm{fin}}}
\newcommand{\Pcl}{P_{\cl}}
\newcommand{\Pclell}{(\Pcl^+)_\ell}
\newcommand{\Paf}{P_{\af}}
\newcommand{\alcl}{\alpha^{\cl}}
\newcommand{\alb}{\overline{\alpha}}
\newcommand{\Lacl}{\La^{\cl}}
\newcommand{\wtcl}{\wt_{\cl}}
\newcommand{\HH}{\mathcal{H}}
\newcommand{\Bmin}{\Bc_{\min}}
\newcommand{\nut}{\widetilde{\nu}}
\newcommand{\mut}{\widetilde{\mu}}
\newcommand{\sigg}{\sigma}
\newcommand{\bb}{\overline{b}}
\newcommand{\PP}{\mathbb{P}}
\newcommand{\Eb}{\overline{E}}
	
\newcommand{\sln}{sl_n}
\newcommand{\slnhat}{\widehat{sl}_n}
\newcommand{\pslnhat}{\widehat{sl}'_n}
\newcommand{\RSK}{\mathrm{RSK}}

\newcommand{\std}{\mathrm{std}}
\newcommand{\shape}{\mathrm{shape}}

\newcommand{\la}{\lambda}
\newcommand{\La}{\Lambda}
\newcommand{\lab}{\overline{\la}}
\newcommand{\Lab}{\overline{\La}}
\newcommand{\Z}{\mathbb{Z}}
\newcommand{\N}{\mathbb{N}}
\newcommand{\ST}{\mathrm{ST}}
\newcommand{\CLR}{\mathrm{CLR}}
\newcommand{\RLR}{\mathrm{RLR}}
\newcommand{\LRT}{\mathrm{LR}}
\newcommand{\CST}{\mathrm{CST}}
\newcommand{\SCST}{\mathrm{SCST}}
\newcommand{\SLRT}{\mathrm{SLR}}
\newcommand{\Conf}{\mathrm{C}}
\newcommand{\RC}{\mathrm{RC}}
\newcommand{\Bc}{\mathcal{B}}
\newcommand{\word}{\mathrm{word}}
\newcommand{\cc}{\mathrm{cc}}
\newcommand{\rows}{\mathrm{r}}
\newcommand{\charge}{c}
\newcommand{\content}{\mathrm{content}}

\newcommand{\PU}{\mathcal{P}}

\newcommand{\nub}{\overline{\nu}}
\newcommand{\Jb}{\overline{\mathrm{J}}}

\newcommand{\tr}{\mathrm{tr}}
\newcommand{\LRtr}{\mathrm{tr}_{\mathrm{LR}}}
\newcommand{\RCtr}{\mathrm{tr}_{\mathrm{RC}}}

\newcommand{\qbin}[2]{\genfrac{[}{]}{0pt}{}{#1}{#2}}
\newcommand{\qbins}[2]{{\textstyle\genfrac{[}{]}{0pt}{}{#1}{#2}}}

\newcommand{\Rl}{R^<}
\newcommand{\il}{i^<}
\newcommand{\jl}{j^<}
\newcommand{\ick}{i^\vee}
\newcommand{\Ick}{I^\vee}

\newcommand{\lt}{\widetilde{\ell}}
\newcommand{\ltb}{\overline{\lt}}
\newcommand{\psib}{\overline{\psi}}
\newcommand{\psit}{\widetilde{\psi}}
\newcommand{\phib}{\overline{\phi}}
\newcommand{\phit}{\widetilde{\phi}}

\newcommand{\Lb}{\overline{L}}

\newcommand{\bs}{\boldsymbol}
\newcommand{\vm}{\bs{m}}
\newcommand{\vn}{\bs{n}}
\newcommand{\vL}{\bs{L}}
\newcommand{\vu}{\bs{u}}
\newcommand{\vf}{\bs{f}}
\newcommand{\ve}{\bs{e}}
\newcommand{\vv}{\bs{v}}
\newcommand{\vP}{\bs{P}}
\newcommand{\vLb}{\overline{\bs{L}}}

\newcommand{\x}{x}
\newcommand{\xb}{\overline{x}}
\newcommand{\tb}{\overline{t}}
\newcommand{\s}{s}

\newcommand{\db}{\overline{\delta}}

\newcommand{\Wb}{\overline{W}}

\newcommand{\hLap}{\widehat{\La'}}
\newcommand{\hLa}{\widehat{\La}}
\newcommand{\hatb}{\widehat{b}}
\newcommand{\hatc}{\widehat{c}}
\newcommand{\haty}{\widehat{y}}
\newcommand{\hatx}{\widehat{x}}
\newcommand{\hatw}{\widehat{w}}

\begin{document}

\title{Fermionic formulas for level-restricted generalized
Kostka polynomials and coset branching functions}
 
\author{Anne Schilling}
\address{Department of Mathematics\\
Massachusetts Institute of Technology\\
77 Massachusetts Avenue\\
Cambridge, MA 02139\\
U.S.A.}
\email{anne@math.mit.edu}

\author{Mark Shimozono}
\address{Department of Mathematics\\
Virginia Tech\\
Blacksburg, VA 24061-0123\\
U.S.A.}
\email{mshimo@math.vt.edu}
\thanks{The second author is partially supported by NSF grant DMS-9800941.}


\begin{abstract}
Level-restricted paths play an important r\^{o}le in crystal theory.
They correspond to certain highest weight vectors of modules of 
quantum affine algebras. We show that the recently established 
bijection between Littlewood--Richardson tableaux and rigged 
configurations is well-behaved with respect to level-restriction 
and give an explicit characterization of level-restricted rigged
configurations. As a consequence a new general fermionic 
formula for the level-restricted generalized Kostka polynomial
is obtained. Some coset branching functions of type $A$ are computed
by taking limits of these fermionic formulas.
\end{abstract}

\maketitle

\pagestyle{plain}

\section{Introduction}
Generalized Kostka polynomials \cite{KS,SW,S,S2,S3,ShW}
are $q$-analogues of the tensor product multiplicity
\begin{equation}\label{multi}
  c^\la_R = \dim \mathrm{Hom}_{\sln}
   (V^\la,V^{R_1} \otimes \dots \otimes V^{R_L}),
\end{equation}
where $\la$ is a partition, $R=(R_1,\ldots,R_L)$ is a sequence 
of rectangles and $V^\la$ is the irreducible integrable highest weight
module of highest weight $\la$ over the quantized enveloping
algebra $U_q(\sln)$. The generalized Kostka polynomials can be expressed 
as generating functions of classically restricted paths \cite{NY,SW,S3}.
In terms of the theory of $U_q(\sln)$-crystals~\cite{Kas,KN} these paths 
correspond to the highest weight vectors of tensor products of perfect 
crystals. The statistic is given by the energy function on paths. 

The $U_q(\sln)$-crystal structure can be extended to a
$U_q(\widehat{sl}'_n)$-crystal structure~\cite{KKMMNN}.
For particular weights, the highest weight vectors of the 
$U_q(\widehat{sl}'_n)$-modules correspond to level-restricted paths.
Hence it is natural to consider the generating functions of
level-restricted paths, giving rise to level-restricted
generalized Kostka polynomials which will take a lead r\^{o}le in this paper.
The notion of level-restriction is also very important in the
context of restricted-solid-on-solid (RSOS) models in statistical
mechanics~\cite{Ba} and fusion models in conformal field 
theory~\cite{V}. The one-dimensional configuration sums of
RSOS models are generating functions of level-restricted paths
(see for example~\cite{ABF,DJKMO,JMO}). The structure constants of the 
fusion algebras of Wess--Zumino--Witten conformal field theories
are exactly the level-restricted analogues of the Littlewood--Richardson
coefficients in \eqref{multi} as shown by Kac~\cite[Exercise 13.35]{Kac} 
and Walton~\cite{Wal,Wal1}. $q$-Analogues of these level-restricted
Littlewood--Richardson coefficients in terms of ribbon tableaux were
proposed in ref.~\cite{FLOT}.

The generalized Kostka polynomial admits a fermionic (or quasi-particle) 
formula~\cite{KSS}. Fermionic formulas originate from the Bethe 
Ansatz~\cite{Bethe} which is a technique to construct eigenvectors and 
eigenvalues of row-to-row transfer matrices of statistical mechanical models. 
Under certain assumptions (the string hypothesis) it is possible
to count the solutions of the Bethe equations resulting in
fermionic expressions which look like sums of products of binomial 
coefficients. The Kostka numbers arise in the study of the
$XXX$ model in this way~\cite{KKR,K,KR}. Fermionic formulas are of interest
in physics since they reflect the particle structure of the underlying
model~\cite{KKMM1,KKMM2} and also reveal information about the
exclusion statistics of the particles~\cite{BM,BCR,BS}.

The fermionic formula of the Kostka polynomial
can be combinatorialized by taking a weighted sum over sets of 
rigged configurations~\cite{KKR,K,KR}. In ref.~\cite{KSS}
the fermionic formula for the generalized Kostka polynomial was proven
by establishing a statistic-preserving bijection between
Littlewood--Richardson tableaux and rigged configurations.
In this paper we show that this bijection is well-behaved with respect
to level-restriction and we give an explicit characterization
of level-restricted rigged configurations (see Definition~\ref{def levrc}
and Theorem~\ref{rect res}). This enables us to obtain a
combinatorial formula for the level-restricted generalized Kostka 
polynomials as the generating function of level-restricted
rigged configurations (see Theorem~\ref{main}). As an immediate
consequence this proves a new general fermionic 
formula for the level-restricted generalized Kostka polynomial
(see Theorem~\ref{thm qp} and Eq.~\eqref{Kmn}). Special cases of this 
formula were conjectured in refs.~\cite{Das97,HKKOTY,HKOTY,K98,SW,WP}.
As opposed to some definitions of ``fermionic formulas'' the expression 
of Theorem~\ref{thm qp} involves in general explicit negative signs. 
However, we would like to point out that because of the equivalent 
combinatorial formulation in terms of rigged configurations as given 
in Theorem~\ref{main} the fermionic sum is manifestly positive 
(i.e., a polynomial with positive coefficients). 

The branching functions of type $A$ can be described in terms of 
crystal graphs of irreducible integrable highest weight 
$U_q(\widehat{sl}_n)$-modules.
For certain triples of weights they can be expressed as limits 
of level-restricted generalized Kostka polynomials. The structure
of the rigged configurations allows one to take this limit, thereby
yielding a fermionic formula for the corresponding branching functions
(see Eq.~\eqref{bf ferm}). The derivation of this formula requires
the knowledge of the ground state energy, which is obtained
from the explicit construction of certain local isomorphisms of 
perfect crystals (see Theorem~\ref{iso}).
A more complete set of branching functions can be obtained by
considering ``skew'' level-restricted generalized Kostka
polynomials. We conjecture that rigged configurations
are also well-behaved with respect to skew shapes
(see Conjecture~\ref{conj skew}).

The paper is structured as follows. Section~\ref{sec note} sets out
notation used in the paper. In Section~\ref{sec path} we review
some crystal theory, in particular the definition of level-restricted
paths, which are used to define the level-restricted generalized
Kostka polynomials. Littlewood--Richardson tableaux and their
level-restricted counterparts are defined in Section~\ref{sec tab}.
The formulation of the generalized Kostka polynomials in terms
of Littlewood--Richardson tableaux with charge statistic is 
necessary for the proof of the fermionic formula which makes
use of the bijection between Littlewood--Richardson tableaux
and rigged configurations. The latter are subject of 
Section~\ref{sec rc} which also contains the new definition
of level-restricted rigged configurations and our main Theorem~\ref{main}.
The proof of this theorem is reserved for Section~\ref{sec proof}.
The fermionic formulas for the level-restricted Kostka polynomial
and the type $A$ branching functions are given in 
Sections~\ref{sec qp Kostka} and \ref{sec bf}, respectively.

\subsection*{Acknowledgements}
We are deeply indebted to Anatol Kirillov for stimulating discussions.
We would also like to thank Peter Bouwknegt, Srinandan Dasmahapatra,
Atsuo Kuniba, and Masato Okado for helpful discussions.

\section{Notation}
\label{sec note}

All partitions are assumed to have $n$ parts, some of which may
be zero.  Let $R=(R_1,R_2,\dots,R_L)$ be a sequence of partitions
whose Ferrers diagrams are rectangles.  Let $R_j$ have
$\mu_j$ columns and $\eta_j$ rows for $1\le j\le L$.
We adopt the English notation for partitions and tableaux.
Unless otherwise specified, all tableaux are assumed to be
column-strict (that is, the entries in each row weakly increase
from left to right and in each column strictly increase from top to bottom).

\section{Paths}
\label{sec path}

The main goal of this section is to define the level-restricted
generalized Kostka polynomials.  These polynomials are
defined in terms of certain finite $U_q(\widehat{sl}_n)$-crystal graphs
whose elements are called paths.  The theory of crystal graphs was invented 
by Kashiwara \cite{Kas}, who showed that
the quantized universal enveloping algebras of Kac-Moody algebras and
their integrable highest weight modules admit special bases
whose structure at $q=0$ is specified by a colored graph known as the
crystal graph.  The crystal graphs for the finite-dimensional
irreducible modules for the classical Lie algebras were computed explicitly
by Kashiwara and Nakashima \cite{KN}.  The theory of perfect crystals
gave a realization of the crystal graphs of the irreducible integrable
highest weight modules for affine Kac-Moody algebras, as certain
eventually periodic sequences of elements taken from
finite crystal graphs \cite{KKMMNN2}.  This realization is
used for the main application, some new explicit formulas for
coset branching functions of type $A$.

\subsection{Crystal graphs}

Let $U_q(\gggg)$ be the quantized universal enveloping algebra
for the Kac-Moody algebra $\gggg$.  Let $I$ be an indexing
set for the Dynkin diagram of $\gggg$, $P$ the weight lattice
of $\gggg$, $P^*$ the dual lattice,
$\{\alpha_i\mid i\in I\}$ the (not necessarily linearly
independent) simple roots, $\{h_i\mid i\in I\}$ the
simple coroots, and $\{\La_i\mid i\in I\}$ the fundamental weights.
Let $\inner{\cdot}{\cdot}$ denote the natural pairing of $P^*$ and $P$.

Suppose $V$ is a $U_q(\gggg)$-module with
crystal graph $\Bc$.  Then $\Bc$ is a directed graph whose
vertex set (also denoted $\Bc$) indexes a basis of weight vectors of $V$,
and has directed edges colored by the elements of the set $I$.
The edges may be viewed as a combinatorial version of
the action of Chevalley generators.
This graph has the property that for every $b\in \Bc$ and $i\in I$,
there is at most one edge colored $i$ entering (resp. leaving) $b$.
If there is an edge $b\rightarrow b'$ colored $i$, denote this by
$f_i(b)=b'$ and $e_i(b')=b$.  If there is no edge colored $i$
leaving $b$ (resp. entering $b'$) then say that
$f_i(b)$ (resp. $e_i(b')$) is undefined.  The $f_i$ and $e_i$ are called
Kashiwara lowering and raising operators.
Define $\phi_i(b)$ (resp. $\epsilon_i(b)$) to be the
maximum $m\in\N$ such that $f_i^m(b)$ (resp. $e_i^m(b)$) is defined.
There is a weight function $\wt:\Bc\rightarrow P$ that satisfies the
following properties:
\begin{equation} \label{crystal string}
\begin{split}
  \wt(f_i(b)) &= \wt(b)-\alpha_i, \\
  \wt(e_i(b)) &= \wt(b)+\alpha_i, \\
  \inner{h_i}{\wt(b)} &= \phi_i(b) - \epsilon_i(b).
\end{split}
\end{equation}
$\Bc$ is called a $P$-weighted $I$-crystal.

Let $P^+=\{\La\in P \mid \inner{h_i}{\La}\ge0,\,\, \forall i\in I\}$
be the set of dominant integral weights.
For $\La\in P^+$ denote by $\VV(\La)$ the irreducible integrable
highest weight $U_q(\gggg)$-module of highest weight $\La$.
Let $\BB(\La)$ be its crystal graph.

Say that an element $b\in \Bc$ of the $P$-weighted $I$-crystal $\Bc$
is a highest weight vector if $\epsilon_i(b)=0$ for all $i\in I$.

Let $u_\La$ be the highest weight vector in $\BB(\La)$.
By \eqref{crystal string}, for all $i\in I$,
\begin{equation} \label{highest}
\begin{split}
  \epsilon_i(u_\La) &= 0,  \\
  \phi_i(u_\La) &= \inner{h_i}{\La}.
\end{split}
\end{equation}

Let $\Bc'$ be the crystal graph of a $U_q(\gggg)$-module $V'$.
A morphism of $P$-weighted $I$-crystals
is a map $\tau:\Bc\rightarrow \Bc'$ such that
$\wt(\tau(b))=\wt(b)$ and 
$\tau(f_i(b))=f_i(\tau(b))$ for all $b\in \Bc$ and $i\in I$.
In particular $f_i(b)$ is defined if and only if $f_i(\tau(b))$ is.

Suppose $V$ and $V'$ are $U_q(\gggg)$-modules with crystal graphs
$\Bc$ and $\Bc'$ respectively.  Then $V \otimes V'$ admits a crystal
graph denoted $\Bc \otimes \Bc'$ which is equal to the direct product
$\Bc \times \Bc'$ as a set.  We use the opposite of the convention
used in the literature.  Define
\begin{equation} \label{tensor lower}
  f_i(b \otimes b') =
  \begin{cases}
    b \otimes f_i(b') & \text{if $\phi_i(b')>\epsilon_i(b)$,} \\
    f_i(b) \otimes b' & \text{if $\phi_i(b')\le \epsilon_i(b)$ and
    			$\phi_i(b)>0$,} \\
    \text{undefined} & \text{otherwise.}    
  \end{cases}
\end{equation}
Equivalently,
\begin{equation} \label{tensor raise}
  e_i(b \otimes b') =
  \begin{cases}
    e_i(b) \otimes b' & \text{if $\phi_i(b')<\epsilon_i(b)$,} \\
    b \otimes e_i(b') & \text{if $\phi_i(b')\ge\epsilon_i(b)$ and
    			$\epsilon_i(b')>0$,} \\
    \text{undefined} & \text{otherwise.}
  \end{cases}
\end{equation}
One has
\begin{equation} \label{tensor crystal}
\begin{split}
 \phi_i(b \otimes b') &=
	\phi_i(b) + \max\{0,\phi_i(b')-\epsilon_i(b)\}, \\
 \epsilon_i(b \otimes b') &=
	\max\{0,\epsilon_i(b)-\phi_i(b')\}+\epsilon_i(b').
\end{split}
\end{equation}
Finally $\wt:\Bc\otimes \Bc'\rightarrow P$ is defined
by $\wt(b \otimes b') = \wt_\Bc(b)+\wt_{\Bc'}(b')$ where
$\wt_\Bc:\Bc\rightarrow P$ and
$\wt_{\Bc'}:\Bc'\rightarrow P$ are the weight functions for
$\Bc$ and $\Bc'$.

This construction is "associative", that is, the $P$-weighted
$I$-crystals form a tensor category.

\begin{rem} \label{act at position}
It follows from \eqref{tensor raise} that if $b=b_L\otimes\dots\otimes b_1$
and $e_i(b)$ is defined, then
$e_i(b)=b_L \otimes \dots \otimes b_{j+1} \otimes e_i(b_j) \otimes
b_{j-1} \otimes \dots \otimes b_1$ for some $1\le j\le L$.
\end{rem}

\subsection{$U_q(\sln)$-crystal graphs on tableaux}

Let $J = \{1,2,\dots,n-1\}$ be the indexing set for the
Dynkin diagram of type $A_{n-1}$, with weight lattice $\Pfin$,
simple roots $\{\alb_i\mid i\in J\}$, fundamental weights
$\{\Lab_i\mid i\in J\}$, and simple coroots $\{h_i\mid i\in J\}$.

Let $\la=(\la_1\ge\la_2\ge\dots\ge\la_n)\in\N^n$ be a partition.
There is a natural projection $\Z^n\rightarrow\Pfin$ denoted
$\la\mapsto\lab=\sum_{i=1}^{n-1} (\la_i-\la_{i+1}) \Lab_i$.
Let $V(\lab)$ be the irreducible integrable highest weight module
of highest weight $\lab$ over the quantized universal enveloping algebra 
$U_q(\sln)$~\cite{KN}.  By abuse of notation we shall write
$V^\la=V(\lab)$ and denote the crystal graph of $V^\la$ by $\Bc_\la$.

As a set $\Bc_{\la}$ may be realized as the set of tableaux of shape $\la$ 
over the alphabet $\{1,2,\ldots,n\}$.
Define the content of $b\in \Bc_\la$ by $\content(b)=(c_1,\dots,c_n)\in\N^n$
where $c_j$ is the number of times the letter $j$ appears in $b$.
The weight function $\wt:\Bc_\la\rightarrow\Pfin$ is given by
sending $b$ to the image of $\content(b)$ under the projection
$\Z^n\rightarrow\Pfin$.  The row-reading word of $b$ is defined by
$\word(b) = \cdots w_2 w_1$ where $w_r$ is the word obtained by
reading the $r$-th row of $b$ from left to right.
This definition is useful even in the context that $b$ is a skew tableau.

The edges of $\Bc_\la$ are given as follows.
First let $v$ be a word in the alphabet $\{1,2,\dots,n\}$.
View each letter $i$ (resp. $i+1$)
of $v$ as a closing (resp. opening) parenthesis, ignoring other letters.
Now iterate the following step: declare each adjacent pair
of matched parentheses to be invisible.
Repeat this until there are no matching pairs of visible parentheses.  
At the end the result must be a sequence of closing parentheses
(say $p$ of them) followed by a sequence of opening parentheses
(say $q$ of them).  The unmatched (visible) subword is of the form
$i^p (i+1)^q$.  If $p>0$ (resp. $q>0$) then $f_i(v)$
(resp. $e_i(v)$) is obtained from $v$ by replacing the unmatched
subword $i^p (i+1)^q$ by $i^{p-1} (i+1)^{q+1}$
(resp. $i^{p+1} (i+1)^{q-1}$).  Then
$\phi_i(v)=p$, $\epsilon_i(v)=q$, and $f_i(v)$ (resp. $e_i(v)$) is
defined if and only if $p>0$ (resp. $q>0$).

For the tableau $b\in\Bc_\la$, let $f_i(b)$ be undefined
if $f_i(\word(b))$ is; otherwise define $f_i(b)$ to be the unique
(not necessarily column-strict) tableau of shape $\la$ such that
$\word(f_i(b))=f_i(\word(b))$.  It is easy to verify that
when defined, $f_i(b)$ is a column-strict tableau.
Consequently $\phi_i(b)=\phi_i(\word(b))$.
The operator $e_i$ and the quantity $\epsilon_i(b)$ are
defined similarly.

\subsection{$U_q(\pslnhat)$-crystal structure on rectangular tableaux}

There is an inclusion of algebras $U_q(\sln)\subset U_q(\pslnhat)$
where $U_q(\pslnhat)$ is the quantized universal enveloping algebra
corresponding to the derived subalgebra $\pslnhat$ of the affine Kac-Moody
algebra $\slnhat$ \cite{Kac}.  Let $I=\{0,1,2,\dots,n-1\}$ be the
index set for the Dynkin diagram of $A_{n-1}^{(1)}$.
Let $\Pcl$ be the weight lattice of $\pslnhat$, with
(linearly dependent) simple roots $\{\alcl_i\mid i\in I\}$,
simple coroots $\{h_i\mid i\in I\}$, and fundamental weights
$\{\Lacl_i\mid i\in I\}$.  The simple roots satisfy the relation
$\alcl_0 = -\sum_{i\in J} \alcl_i$.  There is a natural projection
$\Pcl\rightarrow \Pfin$ with kernel $\Z\La_0$ such that
$\Lacl_i\mapsto \Lab_i$ for $i\in J$ and $\Lacl_0\mapsto 0$.
Let $\cl:\Pfin\rightarrow\Pcl$ be the section of the above
projection defined by $\cl(\Lab_i) = \Lacl_i - \Lacl_0$ for $i\in J$.
Let $c\in\pslnhat$ be the canonical central element.
The level of a weight $\La\in\Pcl$ is defined by $\inner{c}{\La}$.
Let $\Pclell=\{\La\in\Pcl^+\mid \inner{c}{\La}=\ell\}$.

Suppose $V$ is a finite-dimensional
$U_q(\pslnhat)$-module that has a crystal graph $\Bc$ (not all do);
$\Bc$ is a $\Pcl$-weighted $I$-crystal.
A weight function $\wtcl:\Bc\rightarrow\Pcl$ may be given by
$\wtcl(b) = \cl(\wt(b))$ where $\wt:\Bc\rightarrow\Pfin$
is the weight function on the set $\Bc$ viewed as a
$U_q(\sln)$-crystal graph.
In addition to being a $U_q(\sln)$-crystal graph,
$\Bc$ also has some edges colored $0$.  The action of $U_q(\sln)$
on $V^\la$ extends to an action of $U_q(\pslnhat)$ which admits a
crystal structure, if and only if the partition $\la$ is a rectangle
\cite{KKMMNN,NY}.  If $\la$ is the rectangle with $k$ rows
and $m$ columns, then write $V^{k,m}$ for the $U_q(\pslnhat)$-module with
$U_q(\sln)$-structure $V^\la$ and denote its crystal graph
by $\Bc^{k,m}$.  If one of $m$ or $k$ is $1$, then
it is easy to give $e_0$ and $f_0$ explicitly on $\Bc^{k,m}$,
for in this case the weight spaces of $V^{k,m}$
are one-dimensional, and the zero edges can be
deduced from \eqref{crystal string} \cite{KKMMNN}.
The general case is given as follows \cite{S3}.

We shall first define a content-rotating
bijection $\psi^{-1}: \Bc^{k,m}\rightarrow \Bc^{k,m}$.
Let $b\in \Bc^{k,m}$ be a tableau, say of content $(c_1,c_2,\dots,c_n)$.
$\psi^{-1}(b)$ will have content $(c_2,c_3,\dots,c_n,c_1)$.
Remove all the letters $1$ from $b$, leaving a vacant horizontal strip
of size $c_1$ in the northwest corner of $b$.  Compute Schensted's $P$
tableau \cite{Schen} of the row-reading word of this skew subtableau.
It can be shown that this yields a tableau of the shape obtained by
removing $c_1$ cells from the last row of the rectangle $(m^k)$.
Subtract one from the value of each entry of this tableau, and
then fill in the $c_1$ vacant cells in the last row of the rectangle
$(m^k)$ with the letter $n$.  It can be shown that
$\psi^{-1}$ is a well-defined bijection, whose inverse
$\psi$ can be given by a similar algorithm.  Then
\begin{equation} \label{psi and raising}
\begin{split}
  f_i &= \psi^{-1} \circ f_{i+1} \circ \psi, \\
  e_i &= \psi^{-1} \circ e_{i+1} \circ \psi
\end{split}
\end{equation}
for all $i$ where indices are taken modulo $n$; in particular
for $i=0$ this defines explicitly the operators $e_0$ and $f_0$.

\subsection{Sequences of rectangular tableaux}

For a sequence of rectangles $R$, consider the tensor product
$V^{R_L} \otimes \dots\otimes V^{R_1}$.  Its
$U_q(\pslnhat)$-crystal graph has underlying set
$\PU_R = \Bc_{R_L} \otimes \dots\otimes \Bc_{R_1}$,
where the tensor symbols denote the Cartesian product of sets.
A typical element of $\PU_R$ is called a path and is written
$b = b_L \otimes \dots\otimes b_2 \otimes b_1$
where $b_j\in \Bc_{R_j}$ is a tableau of shape $R_j$.

The edges of the crystal graph $\PU_R$ are given explicitly as follows.
Define the word of a path $b$ by
\begin{equation*}
  \word(b) = \word(b_L) \dotsm \word(b_2)\word(b_1).
\end{equation*}
Then for $i=1,2,\ldots,n-1$ (as in the definition of $f_i$ for
$b\in \Bc_\la$), if $f_i(\word(b))$ is undefined, let $f_i(b)$ be undefined;
otherwise it not hard to see that there is a unique
path $f_i(b)\in \PU_R$ such that $\word(f_i(b))=f_i(\word(b))$.
To define $f_0$, let $\psi(b) = \psi(b_L) \otimes \dots \otimes \psi(b_1)$
and $f_0 = \psi^{-1} \circ f_1 \circ \psi$.
This definition is equivalent to that given by taking the
above definition of $f_i$ on the crystals $\Bc_{R_j}$ and then applying the
rule for lowering operators on tensor products \eqref{tensor lower}.
The action of $e_i$ for $i\in I$ is defined analogously.

\subsection{Integrable affine crystals}

Consider the affine Kac-Moody algebra $\slnhat$,
with weight lattice $\Paf$, independent simple roots
$\{\alpha_i\mid i\in I\}$, simple coroots $\{h_i\mid i\in I\}$, and
fundamental weights $\{\La_i\mid i\in I\}$.  Let $\delta\in \Paf$ be
the null root.  There is a natural projection
which we shall by abuse of notation also call $\cl:\Paf\rightarrow\Pcl$
such that $\cl(\delta)=0$ and $\cl(\La_i)=\Lacl_i$ for $i\in I$.
Write $\af:\Pcl\rightarrow\Paf$ for the section of $\cl$
given by $\af(\Lacl_i) = \La_i$ for $i\in I$.

Let $\La\in \Pcl^+$ be a dominant integral weight
and $\BB(\La)$ the crystal graph of the irreducible integrable highest
weight $U_q(\pslnhat)$-module of highest weight $\La$.
If $\La\not=0$ then $\BB(\La)$ is infinite.  The set of weights in $\Paf$ that
project by $\cl$ to $\La$ are given by
$\cl^{-1}(\La) = \{\af(\La)+j\delta\mid j\in\Z\}$.
Now fix $j$.  The irreducible integrable highest weight
$U_q(\slnhat)$-crystal graph $\BB(\af(\La)+j \delta)$ may be identified
with $\BB(\La)$ as sets and as $I$-crystals (independent of $j$).
The weight functions for $\BB(\af(\La)+j\delta)$ and
$\BB(\af(\La))$ differ by the global constant $j\delta$.
The weight function $\BB(\La)\rightarrow \Z$
is obtained by composing the weight function
for $\BB(\af(\La)+j \delta)$, with the projection $\cl:\Paf\rightarrow\Pcl$.

The set $\BB(\La)$ is then endowed with an induced $\Z$-grading
$E:\BB(\La)\rightarrow\N$ defined by $E(b) = \inner{d}{\wt(b)}$
where $\BB(\La)$ is identified with $\BB(\af(\La))$,
$\wt:\BB(\af(\La))\rightarrow \Paf$ is the weight function and
$d\in \Paf^*$ is the degree generator.
The map $\inner{d}{\cdot}$ takes the coefficient of the element $\delta$
of an element in $\Paf$ when written in the basis
$\{\La_i\mid i\in I\} \cup \{\delta\}$.

\subsection{Energy function on finite paths}

The set of paths $\PU_R$ has a natural statistic called the energy
function.  The definitions here follow \cite{NY}.

Consider first the case that $R=(R_1,R_2)$ is a sequence of two rectangles.
Let $\Bc_j=\Bc_{R_j}$ for
$1\le j\le 2$.  Since $\Bc_2\otimes \Bc_1$ is a connected crystal graph,
there is a unique $U_q(\pslnhat)$-crystal graph isomorphism
\begin{equation} \label{local iso}
 \sigma: \Bc_2 \otimes \Bc_1 \cong \Bc_1 \otimes \Bc_2.
\end{equation}
This is called the local isomorphism (see Section \ref{sec auto} for an
explicit construction).  Write $\sigma(b_2\otimes b_1)=b_1'\otimes b_2'$.
Then there is a unique (up to a global additive constant) map
$H:\Bc_2 \otimes \Bc_1 \rightarrow \Z$ such that
\begin{equation} \label{loc energy}
  H(e_i(b_2 \otimes b_1)) =  H(b_2\otimes b_1) +
  \begin{cases}
    -1 & \text{if $i=0$, $e_0(b_2\otimes b_1)=e_0 b_2 \otimes b_1$} \\
       & \text{and $e_0(b_1'\otimes b_2')=e_0 b_1' \otimes b_2'$,} \\
    1 & \text{if $i=0$, $e_0(b_2\otimes b_1)=b_2 \otimes e_0 b_1$} \\
    & \text{and $e_0(b_1'\otimes b_2')=b_1' \otimes e_0 b_2'$,} \\
    0 & \text{otherwise.}
  \end{cases}
\end{equation}
This map is called the local energy function.  By definition it
is invariant under the local isomorphism and under $f_i$ and $e_i$
for $i\in J$.  Let us normalize
it by the condition that $H(u_2 \otimes u_1)=|R_1\cap R_2|$ where
$u_j$ is the $U_q(\sln)$ highest weight vector of $\Bc_j$ for $1\le j\le 2$,
$R_1\cap R_2$ is the intersection of the Ferrers diagrams of $R_1$
and $R_2$, and $|R_1 \cap R_2|$ is the number of cells in this
intersection.  Explicitly $|R_1\cap R_2|=\min\{\eta_1,\eta_2\}
\min\{\mu_1,\mu_2\}$.  If $\eta_1+\eta_2\le n$ then
the local energy function attains precisely the values from $0$ to
$|R_1\cap R_2|$.

Now let $R=(R_1,\dots,R_L)$ be a sequence of rectangles and 
$b=b_L \otimes \dots \otimes b_1\in \PU_R$.  For $1\le p\le L-1$
let $\sigma_p$ denote the local isomorphism that exchanges
the tensor factors in the $p$-th and $(p+1)$-th positions.
For $1\le i<j\le L$, let $b_j^{(i+1)}$ be the $(i+1)$-th tensor factor in
$\sigma_{i+1} \sigma_{i+2} \dots \sigma_{j-1} (b)$.
Then define the energy function
\begin{equation}\label{energy}
  E(b) = \sum_{1\le i<j\le L} H(b_j^{(i+1)}\otimes b_i).
\end{equation}
The value of the energy function is unchanged under local
isomorphisms and under $e_i$ and $f_i$ for $i\in J$,
since the local energy function has this property.

The next lemma follows from the definition of the local energy function.

\begin{lemma} \label{energy down} Suppose $b=b_L\otimes \dots\otimes b_1
\in \PU_R$ is such that $e_0(b)$ is defined and
for any image $b'=b'_L \otimes \dots \otimes b'_1$ of $b$
under a composition of local isomorphisms,
$e_0(b')=b'_L \otimes \dots \otimes b'_{j+1} \otimes
e_0(b'_j) \otimes b'_{j-1} \otimes \dots \otimes b'_1$
where $j\not=1$.  Then $E(e_0(b))=E(b)-1$.
\end{lemma}

If all rectangles $R_j$ are the same then each of the local isomorphisms
is the identity and 
\begin{equation}\label{hom energy}
  E(b) = \sum_{1\le i \le L-1} (L-i) H(b_{i+1}\otimes b_i).
\end{equation}

Say that $b\in \PU_R$ is \textit{classically restricted} if
it is an $\sln$-highest weight vector, that is,
$\epsilon_i(b)=0$ for all $i\in J$.  Equivalently,
$\word(b)$ is a (reverse) lattice permutation (every final subword
has partition content).  Let $\PU_{\La R}$ be
the set of classically restricted paths in $\PU_R$ of weight $\La\in\Pcl$.

It was shown in~\cite{S3} that the generalized Kostka polynomial
(which was originally defined in terms of Littlewood--Richardson
tableaux; see \eqref{gK}) can be expressed as
\begin{equation}\label{K path}
  K_{\la R}(q) = \sum_{b\in \PU_{\cl(\lab) R}} q^{E(b)}.
\end{equation}
This extends the path formulation of the Kostka
polynomial by Nakayashiki and Yamada~\cite{NY}.

\subsection{Level-restricted paths}
\label{sec restr. paths}

Let $\Bc$ be any $\Pcl$-weighted $I$-crystal and $\La\in\Pcl^+$.
Say that $b\in \Bc$ is $\La$-restricted
if $b \otimes u_\La$ is a highest weight vector in the
$\Pcl$-weighted $I$-crystal $\Bc \otimes \BB(\La)$, that is,
$\epsilon_i(b\otimes u_\La)=0$ for all $i\in I$.
Equivalently $\epsilon_i(b) \le \inner{h_i}{\La}$ for all $i \in I$
by \eqref{tensor crystal} and \eqref{highest}.  Denote by
$\HH(\La,\Bc)$ the set of elements $b\in \Bc$ that are $\La$-restricted.
If $\La'\in\Pcl^+$ has the same level as $\La$,
define $\HH(\La,\Bc,\La')$ to be the set of $b\in \HH(\La,\Bc)$ such that
$\wt(b)=\La'-\La\in\Pcl$, that is, the set of $b\in\Bc$ such that
$b \otimes u_\La$ is a highest weight vector of weight $\La'$.
Say that the element $b$ is restricted of level $\ell$ if
it is $(\ell\La_0)$-restricted.  Such paths are also classically
restricted since $\inner{h_i}{\ell\La_0}=0$ for $i\in J$.
Let $\PU^\ell_{\La R}$ denote the
set of paths in $\PU_{\La R}$ that are restricted of level $\ell$.
Letting $\Bc=\PU_R$, this is the same as saying
$\PU^\ell_{\La R} = \HH(\ell\La_0,\Bc,\La+\ell\La_0)$.

Define the level-restricted generalized Kostka polynomial by
\begin{equation}\label{K path lev}
  K^\ell_{\la R}(q) = \sum_{b\in \PU^\ell_{\cl(\lab) R}} q^{E(b)}.
\end{equation}

\subsection{Perfect crystals}

This section is needed to compute the coset branching functions
in Section \ref{sec bf}.
We follow \cite{KKMMNN2}, stating the definitions in the
case of $\pslnhat$.  For any $U_q(\pslnhat)$-crystal $\Bc$,
define $\epsilon,\phi:\Bc\rightarrow\Pcl$ by
$\epsilon(b)=\sum_{i\in I} \epsilon_i(b) \La_i$
and $\phi(b)=\sum_{i\in I} \phi_i(b) \La_i$.

Now let $\ell$ be a positive integer and 
$\Bc$ the crystal graph of a finite dimensional
irreducible $U_q(\pslnhat)$-module $V$.  Say that $\Bc$ is perfect
of level $\ell$ if
\begin{enumerate}
\item $\Bc\otimes \Bc$ is connected.
\item There is a weight $\La'\in\Pcl$ such that $\Bc$ has a unique
vector of weight $\La'$ and all other vectors in $\Bc$ have lower
weight in the Chevalley order, that is,
$\wt(\Bc)\subset\La'-\sum_{i\in J} \N \alpha_i$.
\item $\ell = \min_{b\in \Bc} \inner{c}{\epsilon(b)}$.  
\item The maps $\epsilon$ and $\phi$ restrict to bijections
$\Bmin  \rightarrow \Pclell$ where $\Bmin\subset \Bc$
is the set of $b\in \Bc$ achieving the minimum in 3.
\end{enumerate}

For $\pslnhat$ the perfect crystals of level $\ell$ are precisely those
of the form $\Bc^{k,\ell}$ for $1\le k\le n-1$ \cite{KKMMNN,NY}.
Let $\Bc=\Bc^{k,\ell}$.  The weight $\La'$ can be taken to be 
$\ell (\Lacl_k-\Lacl_0)$.

\begin{ex} \label{min bijection}
We describe the bijections $\epsilon,\phi:\Bmin\rightarrow
\Pclell$ in this example.  Let $\Bc=\Bc^{k,\ell}$.
For this example let $n=6$, $k=3$, $\ell=5$,
and consider the weight $\La=2\La_0+\La_1+\La_2+\La_4$.
As usual subscripts are identified modulo $n$.  The unique
tableau $b\in \Bc^{k,\ell}$ such that $\phi(b)=\La$ is constructed
as follows.  First let $T$ be the following tableau of shape $(\ell^k)$.
Its bottom row contains $\inner{h_i}{\La}$ copies of the letter $i$
for $1\le i\le n$ (here it is $12466$ since the sequence of
$\inner{h_i}{\La}$ for $1\le i\le 6$ is $(1,1,0,1,0,2)$).
Let every letter in $T$ have value one smaller than
the letter directly below it.  Here we have
\begin{equation*}
  T = \begin{matrix}
    -1&0&2&4&4\\
    0&1&3&5&5 \\
    1&2&4&6&6
  \end{matrix}.
\end{equation*}
Let $T_-$ be the subtableau of $T$ consisting of the entries that
are nonpositive and $T_+$ the rest.  Say $T_-$ has shape $\nu$
(here $\nu=(2,1)$).  Let $\nut=(\ell^k)-(\nu_k,\nu_{k-1},\dots,\nu_1)$
(here $\nut=(5,4,3)$).  The desired tableau $b$ is defined as follows.
The restriction of $b$ to the shape $\nut$ is
$P(T_+)$, or equivalently, the tableau obtained
by taking the skew tableau $T_+$ and first pushing all letters
straight upwards to the top of the bounding rectangle $(\ell^k)$,
and then pushing all letters straight to the left inside $(\ell^k)$.
The restriction of $b$ to $(\ell^k)/\nut$ is the tableau of that
skew shape in the alphabet $\{1,2,\dots,n\}$ with maximal entries,
that is, its bottom row is filled with the letter $n$,
the next-to-bottom row is filled with the letter $n-1$, etc. In the example,
\begin{equation*}
  b = \begin{matrix}
    1&1&2&4&4\\
    2&3&5&5&5\\
    4&6&6&6&6
  \end{matrix}.
\end{equation*}

To construct the unique element
$b'\in \Bc^{k,\ell}$ such that $\epsilon(b')=\La$,
let $U$ be the tableau whose first row
has $\inner{h_i}{\La}$ copies of the letter $i+1$ for $1\le i\le n$,
again identifying subscripts modulo $n$; here
$U$ has first row $11235$.  Now let the rest of $U$ be
defined by letting each entry have value one greater than
the entry above it.  So
\begin{equation*}
  U = \begin{matrix}
    1&1&2&3&5\\
    2&2&3&4&6\\
    3&3&4&5&7
  \end{matrix}.
\end{equation*}
Let $U_-$ be the subtableau of $U$ consisting of the
values that are at most $n$.
Let $\mu$ be the shape of $U_-$ and
$\mut = (\ell^k)-(\mu_k,\mu_{k-1},\dots,\mu_1)$.
Here $\mu=(5,5,4)$ and $\mut=(1,0,0)$.  The element $b'$ is defined
as follows.  Its restriction to the skew shape $(\ell^k)/\mut$
is the unique skew tableau $V$ of that
shape such that $P(V)=U_-$, or equivalently, this restriction
is obtained by taking the tableau $U_-$, pushing all letters
directly down within the rectangle $(\ell^k)$ and then
pushing all letters to the right within $(\ell^k)$.
The restriction of $b'$ to the shape $\mut$ is
filled with the smallest letters possible, so that
the first row of this subtableau consists of ones,
the second row consists of twos, etc.  Here
\begin{equation*}
  b' = \begin{matrix}
    1&1&1&2&3\\
    2&2&3&4&5\\
    3&3&4&5&6
  \end{matrix}.
\end{equation*}
\end{ex}

The main theorem for perfect crystals is:

\begin{theorem} \label{perfect} \cite{KKMMNN2} Let $\Bc$ be a perfect crystal
of level $\ell'$ and $\La\in\Pclell$ with $\ell\ge\ell'$.
Then there is an isomorphism of $U_q(\pslnhat)$-crystals
\begin{equation} \label{res iso}
  \Bc \otimes \BB(\La) \cong \bigoplus_{b\in\HH(\La,\Bc)} \BB(\La+\wt(b)).
\end{equation}
\end{theorem}

Suppose now that $\Bc$ is perfect of level $\ell$ and $\La\in\Pclell$.
Write $b(\La)$ for the unique element of $\Bc$ such that
$\phi(b(\La))=\La$.  Theorem \ref{perfect} 
(with $\La$ therein replaced by $\La'=\epsilon(b(\La))$) says that
$\Bc \otimes \BB(\epsilon(b(\La))) \cong \BB(\La)$ with corresponding
highest weight vectors $b(\La) \otimes u_{\epsilon(b(\La))} 
\mapsto u_\La$.  This isomorphism can be iterated.
Let $\sigg:\Bmin\rightarrow\Bmin$ be the unique bijection defined by
$\phi \circ \sigg = \epsilon$.  Then there are isomorphisms
$\Bc^{\otimes N} \otimes \BB(\phi(\sigg^N(b(\La)))) \cong \BB(\La)$
such that the highest weight vector of the left-hand side is
given by $b(\La) \otimes \sigg(b(\La)) \otimes \sigg^2(b(\La)) \otimes
\dots \otimes \sigg^{N-1}(b(\La)) \otimes u_{\phi(\sigg^N(b(\La)))}$.
For the $U_q(\pslnhat)$ perfect crystals $\Bc^{k,\ell}$,
it can be shown that the map $\sigg$ is none other than the power
$\psi^{-k}$ of the content rotating map $\psi$.  Moreover
if $\sigg$ is extended to a bijection
$\sigg:\Bc^{k,\ell}\rightarrow \Bc^{k,\ell}$ by defining $\sigg=\psi^{-k}$,
then the extended function
also satisfies $\phi(\sigg(b))=\epsilon(b)$ for all $b\in \Bc^{k,\ell}$
not just for $b\in\Bmin$.  Since the bijection $\psi$ on $\Bc^{k,\ell}$
has order $n$, the bijection $\sigg$ has order $n/\gcd(n,k)$.
The ground state path for the pair $(\La,\Bc)$ is by definition
the infinite periodic sequence $\bb=\bb_1\otimes \bb_2 \otimes \dots$
where $\bb_i = \sigg^{i-1}(b(\La))$.

Let $\PP(\La,\Bc)$ be the set of all semi-infinite sequences
$b=b_1\otimes b_2 \otimes \dots$ of elements in $\Bc$ such that
$b$ eventually agrees with the ground state path $\bb$ for $(\La,\Bc)$.
Then the set $\PP(\La,\Bc)$ has the structure of the crystal
$\BB(\La)$ with highest weight vector $u_\La=\bb$ and weight function
$\wt(b) = \sum_{i\ge1} (\wt(b_i)-\wt(\bb_i))$.
To recover the weight function of the
$U_q(\slnhat)$-crystal $\BB(\af(\La))$,
define the energy function on $\PP(\La,\Bc)$ by
\begin{equation} \label{energy seq}
  E(b) = \sum_{i\ge1} i (H(b_i \otimes b_{i+1})-H(\bb_i \otimes \bb_{i+1}))
\end{equation}
and define the map $\BB(\af(\la))\rightarrow \Paf$ by
$b\mapsto \wt(b) - E(b) \delta$ where $\wt:\BB(\La)\rightarrow\Pcl$.

$\PP(\La,\Bc)$ can be regarded as a direct limit of the finite
crystals $\Bc^{\otimes N}$.  Define the embedding
$i_N:\Bc^{\otimes N} \rightarrow\PP(\La,\Bc)$ by
\begin{equation*}
 b_1 \otimes \dots \otimes b_N \mapsto 
 b_1 \otimes b_2 \otimes b_N \otimes \bb_{N+1} \otimes \bb_{N+2}
\otimes \dots
\end{equation*}
Define $E_N:\Bc^{\otimes N} \rightarrow\Z$ by
$E_N(b_1 \otimes \dots \otimes b_N) =
E(b_1 \otimes \dots \otimes b_N \otimes \bb_{N+1})$
where the $E$ on the right hand side is the energy function
for the finite path space $\Bc^{\otimes N+1}$.
By definition for all $p=b_1\otimes \dots \otimes b_N\in \Bc^{\otimes N}$,
$E(i_N(p))=E_N(p)-E_N(\bb_1 \otimes \dots \otimes \bb_N)$.
Note that the last fixed step $\bb_{N+1}$
is necessary to make the energy function on the finite paths
stable under the embeddings into $\PP(\La,\Bc)$.

\subsection{Standardization embeddings}
\label{sec path embed}

We require certain embeddings of finite path spaces.
Given a sequence of rectangles $R$, let $\rows(R)$
denote the sequence of rectangles given by splitting the
rectangles of $R$ into their constituent rows.  For example,
if $R=((1),(2,2))$ then $\rows(R)=((1),(2),(2))$.
There is a unique embedding
\begin{equation} \label{path embed}
 i_R:\PU_R\hookrightarrow\PU_{\rows(R)}
\end{equation}
defined as follows.  Its explicit computation is based
on transforming $R$ into $\rows(R)$ using two kinds of steps.
\begin{enumerate}
\item Suppose $R_1$ has more than one row ($\eta_1>1$).
Then use the transformation $R\rightarrow
\Rl=((\mu_1),(\mu_1^{\eta_1-1}),R_2,R_3,\dots,R_L)$.
Informally, $\Rl$ is obtained from $R$ by splitting off the first row
of $R_1$.  There is an associated embedding of $U_q(\pslnhat)$-crystal 
graphs $\il_R:\PU_R \rightarrow \PU_{\Rl}$ defined by the
property that $\word(\il(b))=\word(b)$ for all $b\in\PU_R$.
Here it is crucial that the rectangle being split horizontally, is
the first one, for otherwise the embedding does not preserve the
edges labeled by $0$.
\item If $\eta_1=1$, then use a transformation of the form
$R\rightarrow s_p R$ for some $p$.  Here $s_p R$ denotes
the sequence of rectangles obtained by exchanging the
$p$-th and $(p+1)$-th rectangles in $R$.  The associated
isomorphism of $U_q(\pslnhat)$-crystal graphs is the local isomorphism
$\sigma_p:\PU_R\rightarrow \PU_{s_p R}$ defined before.
\end{enumerate}
It is clear that one can transform $R$ into $\rows(R)$ using these
two kinds of steps.  Now fix one such sequence of steps leading from
$R$ to $\rows(R)$, say $R=R^{(0)}\rightarrow R^{(1)}\rightarrow
\dots\rightarrow R^{(N)}=\rows(R)$ where each $R^{(m)}$ is a sequence
of rectangles and each step $R^{(m-1)}\rightarrow R^{(m)}$
is one of the two types defined above.  Define the map
$i^{(m)}:\PU_{R^{(m-1)}}\hookrightarrow \PU_{R^{(m)}}$
by $i^{(m)} = \il_{R^{(m-1)}}$ if the step is of the first kind,
and by $i^{(m)}=\sigma_p$ if it is of the second kind.
Let $i_R:\PU_R\rightarrow \PU_{\rows(R)}$ be the composition
$i_R = i^{(N)} \circ \dots \circ i^{(1)}$.
It can be shown that the map $i_R$ does not depend on
the sequence of the $R^{(m)}$; this is proven in the
equivalent language of Littlewood--Richardson tableaux in \cite{S2}.

\section{Littlewood--Richardson tableaux}
\label{sec tab}

We now review some formulations of type $A$ tensor product
multiplicities that use tableaux.  These tableaux,
which we call Littlewood--Richardson (LR) tableaux, are the
intermediate combinatorial objects between paths
and rigged configurations, which give rise to fermionic expressions.
For the most part, the material in this section is taken from
\cite{SW,S,S2,S3}.

\subsection{Three formulations}
Let $I_1,I_2,\dots,I_L$ be intervals of integers such that
if $i<j$, $x\in I_i$ and $y\in I_j$, then $x<y$.
Set $I = \bigcup_{j=1}^L I_j$.  For each $1\le j \le L$,
fix a tableau $Z_j$ of shape $R_j$ in the alphabet $I_j$.
Define the set $\SLRT(\la;Z)$ to be the set of tableaux $Q$
of shape $\la$ in the alphabet $I$ such that $P(Q|_{I_j})=Z_j$ for all $j$,
where $Q|_{I_j}$ denotes the skew subtableau of $Q$ obtained by restricting
to the alphabet $I_j$, and $P(S)$ denotes the Schensted $P$-tableau
\cite{Schen} of the row-reading word of the skew tableau $S$.
It is well-known that $|\SLRT(\la;Z)|=c_R^\la$, where $c_R^\la$ was
defined in \eqref{multi}.

We shall define three kinds of LR tableaux given by $\SLRT(\la;Z)$
for various choices of intervals $I_j$ and tableaux $Z_j$.
\begin{enumerate}
\item $\LRT(\la;R)$:  Define the set of intervals of integers
$I_j=A_j=[\eta_1+\dots+\eta_{j-1}+1,\eta_1+\dots+\eta_{j-1}+\eta_j]$.
Let $Z_j=Y_j$ be the tableau of shape $R_j$
whose $r$-th row is filled with copies of the $r$-th largest
letter of $A_j$, namely, $\eta_1+\dots+\eta_{j-1}+r$.
Define $\LRT(\la;R) := \SLRT(\la;Y)$.
When $R$ consists of single rows (that is, $\eta_j=1$ for all $j$),
then $\LRT(\la;R)=\CST(\la;\mu)$, the (column-strict) tableaux of shape
$\la$ and content $\mu$.
\item $\CLR(\la;R)$ (Columnwise LR): Let $ZC_1$ be the standard
tableau of shape $R_1$ obtained by placing the numbers $1$ through
$\eta_1$ down the first column, the next $\eta_1$ numbers down the second
column, etc. Continue this process to obtain $ZC_2$, starting
with the next available number, namely, $\eta_1\mu_1 + 1$.
Explicitly, for $1\le j\le L$,
the $(r,c)$-th entry in the $j$-th tableau $ZC_j$ is equal to
$\eta_1\mu_1+\dots+\eta_{j-1}\mu_{j-1}+(c-1)\eta_j+r$.
Let $B_j$ be the interval consisting of the entries of the tableau
$ZC_j$.  Define $\CLR(\la;R) := \SLRT(\la;ZC)$.
\item $\RLR(\la;R)$ (Rowwise LR):
Define this similarly to $\CLR(\la;R)$ but label
by rows, so that the $(r,c)$-th entry of $ZR_j$ is
$\eta_1\mu_1+\dots+\eta_{j-1}\mu_{j-1}+(r-1)\mu_j+c$.
Then let $\RLR(\la;R):=\SLRT(\la;ZR)$.
\end{enumerate}

\begin{ex} Let $R=((1),(2,2))$ and $\la=(3,2)$.  Here
$A_1=\{1\}$, $A_2=\{2,3\}$, and
\begin{equation*}
Y_1=\begin{matrix} 1 \end{matrix}\quad\text{and}
\quad Y_2=\begin{matrix} 2&2\\3&3 \end{matrix}.
\end{equation*}
We have $B_1=\{1\}$, $B_2=\{2,3,4,5\}$, 
\begin{equation*}
ZC_1=\begin{matrix} 1 \end{matrix}\quad\text{and}
\quad ZC_2=\begin{matrix} 2&4\\3&5 \end{matrix}
\end{equation*}
and
\begin{equation*}
ZR_1=\begin{matrix} 1 \end{matrix}\quad\text{and}
\quad ZR_2=\begin{matrix} 2&3\\4&5 \end{matrix}.
\end{equation*}
Observe that
\begin{equation*}
T=\begin{matrix} 1&2&4\\3&5&\end{matrix}
\end{equation*}
is in $\CLR(\la;R)$ since $P(T|_{B_1})=1=ZC_1$ and
$P(T|_{B_2})=ZC_2$.  On the other hand
$T=\begin{matrix}1&3&5\\2&4\end{matrix}$
is not in $\CLR(\la;R)$ since
$P(T|_{B_2})=\begin{matrix}2&3&5\\4&&\end{matrix}\not=ZC_2$.
\end{ex}

\subsection{Obvious bijections among the various LR tableaux}

There are trivial relabeling bijections between the various
kinds of LR tableaux defined above.  We give them explicitly here
for later use.

\begin{enumerate}
\item The bijection $\gamma_R:\CLR(\la;R)\rightarrow \RLR(\la;R)$ is
given by the following relabeling.  Consider an entry $x$
in a standard tableau $S\in\CLR(\la;R)$.  Then $x$ appears in
one of the $ZC$ tableaux, say, it is the $(r,c)$-th entry of $ZC_j$.
Let $y$ be the $(r,c)$-th entry of the rowwise tableau $ZR_j$.
Then replace $x$ by $y$ in $S$.  Performing all such replacements
simultaneously yields $\gamma_R(S)\in\RLR(\la;R)$.
\item The bijection $\std:\LRT(\la;R)\rightarrow \RLR(\la;R)$ is
given by Schensted's standardization map \cite{Schen}.  Let 
$Q\in \LRT(\la;R)$ and $i$ be some entry in $Q$.
Suppose $i$ is the $r$-th largest value in the subinterval $A_j$.
Replace the occurrences of the letter $i$ in $Q$ from left to right
by the consecutive integers given by the $r$-th row of $ZR_j$.
The result of these substitutions is $\std(Q)\in\RLR(\la;R)$.
\item Define a bijection $\beta_R:\LRT(\la;R)\rightarrow\CLR(\la;R)$
by $\gamma_R^{-1} \circ \std$.
\item Observe that ordinary transposition of standard tableaux
restricts to a bijection $\tr:\RLR(\la;R)\leftrightarrow\CLR(\la^t;R^t)$
where $\la^t$ denotes the transpose partition of $\la$
and $R^t=(R_1^t,R_2^t,\dots,R_L^t)$.
\item There is a bijection $\LRtr:\CLR(\la;R)\rightarrow\CLR(\la^t;R^t)$
defined by $\LRtr = \tr \circ \gamma_R$.
\end{enumerate}

\subsection{Paths to tableau pairs}
\label{sec RSK}

The Robinson--Schensted--Knuth correspondence allows one to pass from
paths to pairs of tableaux.  This bijection gives a combinatorial
decomposition of the crystal graph of $\PU_R$ into $U_q(\sln)$ irreducible
components and encodes the energy function in the recording
tableau.

The column insertion version of the Robinson--Schensted--Knuth
correspondence, restricts to a bijection
\begin{equation} \label{path to tabs}
  \RSK:\PU_R \rightarrow \bigcup_\la \CST(\la;\cdot) \times \LRT(\la;R)
\end{equation}
as follows.  Let $b=b_L\otimes\dots\otimes b_2\otimes b_1\in \PU_R$.
Define $P(b) := P(\word(b))$.  
This can be computed by the column insertion of 
$\word(b)$ starting from the right end.
Recall that $b_j$ and $Y_j$ are column-strict tableaux of shape $R_j$.
Let $Q(b)$ be the tableau obtained by recording the insertion of
a letter in $b_j$ by the letter in the corresponding position in $Y_j$.
It can be shown that $Q(b)\in \LRT(\la;R)$, and that
the map \eqref{path to tabs} given by $b\mapsto (P(b),Q(b))$ is a bijection.

\begin{rem} \label{RSK morphism}\mbox{}
\begin{enumerate}
\item This bijection is a morphism of $U_q(\sln)$-crystal graphs
in the sense that $P(e_i(b))=e_i(P(b))$ for $i\in J$.
In particular, $b\in\PU_R$ is classically restricted if and only if
$P(b)$ is a Yamanouchi tableau, that is, its $r$-th row is filled with
copies of the letter $r$ for all $1\le r\le n$.
\item The energy function on paths can be transferred easily to
a statistic on $\LRT(\la;R)$ called the generalized charge (written
$\charge_R$)
such that $\charge_R(Q(b)) = E(b)$.  The generalized charge
is defined explicitly in \eqref{charge} below.
\end{enumerate}
\end{rem}

\begin{ex} Let $R=((1),(2,2))$ and $b\in \PU_R$ given by
\begin{equation*}
  b = \begin{matrix} 1&1\\2&2\end{matrix} \otimes
  \begin{matrix} 1 \end{matrix}.
\end{equation*}
Then $\word(b) = 2211\,1$ and
\begin{equation*}
  P(b) = \begin{matrix} 1&1&1\\2&2&\end{matrix}
  \qquad
  Q(b) = \begin{matrix} 1&2&2\\3&3&\end{matrix}.
\end{equation*}
\end{ex}

\subsection{Generalized Automorphisms of Conjugation}
\label{sec auto}

For the moment let $R=(R_1,R_2)$ and
$\Bc_j = \Bc_{R_j}$ for $1\le j\le 2.$
Recall that the local isomorphism \eqref{local iso}
is the unique isomorphism of $U_q(\pslnhat)$-crystal graphs
$\Bc_2 \otimes \Bc_1 \rightarrow \Bc_1 \otimes \Bc_2$ or
equivalently $\PU_{(R_1,R_2)}\rightarrow \PU_{(R_2,R_1)}$.
Let us make this more explicit.  By Remark \ref{RSK morphism}
we have a commutative diagram of bijections
\begin{equation*}
\begin{CD}
\PU_{(R_1,R_2)} @>{\RSK}>>
	\bigcup_\la \CST(\la) \times \LRT(\la;(R_1,R_2)) \\
  @V{\sigma}VV @VV{\bigcup 1\times s}V \\
 \PU_{(R_2,R_1)} @>>{\RSK}>
	\bigcup_\la \CST(\la) \times \LRT(\la;(R_2,R_1))
\end{CD}
\end{equation*}
such that $P(\sigma(b))=P(b)$.  This induces a bijection
$s:\LRT(\la;(R_1,R_2))\rightarrow \LRT(\la;(R_2,R_1))$ for each
$\la$.  The tensor product $V^{R_2}\otimes V^{R_1}$
is multiplicity-free.  Therefore
the domain and codomain of $s$ are both empty or both
singletons.  Hence the bijection $s$ is unique and can be computed
from the definition of the set $\LRT$.
Then $\sigma(b)$ can be computed by applying RSK to obtain
$(P(b),Q(b))$, then applying $s$ to get $(P(b),s(Q(b))$, and finally,
the inverse of RSK to obtain $\sigma(b)$.

The local energy function is recovered using only the shape
of the tableau pair.
For a tableau $Q\in\LRT(\la;(R_1,R_2))$ let $d(Q)$ be the
number of cells in $Q$ that lie strictly to the right of the
$\max\{\mu_1,\mu_2\}$-th column, or equivalently, strictly to the
right of the shape $R_1\cup R_2$.  Then $H(b) = d(Q(b))$.

Then the $U_q(\pslnhat)$-crystal graph
isomorphism $\sigma_p: \PU_R \rightarrow \PU_{s_p R}$
induces involutions $s_p:\LRT(\la;R)\rightarrow \LRT(\la;s_p R)$
such that the diagram commutes:
\begin{equation*}
\begin{CD}
\PU_R @>{\RSK}>>
	\bigcup_\la \CST(\la) \times \LRT(\la;R) \\
  @V{\sigma_p}VV @VV{\bigcup 1\times s_p}V \\
 \PU_{s_p R} @>>{\RSK}>
	\bigcup_\la \CST(\la) \times \LRT(\la;s_p R).
\end{CD}
\end{equation*}
The map $s_p$ is computed explicitly as follows \cite{S3}.
Let $Q\in \LRT(\la;R)$ and $A_j$ be the alphabets
as in the definition of $\LRT(\la;R)$.
Remove the skew subtableau $U=Q|_{A_p \cup A_{p+1}}$.
Use the usual column insertion of its row reading word,
obtaining a pair of tableaux $(P',Q')$ where
$P'\in \LRT(\rho;(R_p,R_{p+1}))$ for some partition $\rho$ and $Q'$ is
the standard column insertion tableau.  Next replace
$P'$ by $s(P')$ where $s$ is the unique bijection
$\LRT(\rho;(R_p,R_{p+1}))\rightarrow \LRT(\rho;(R_{p+1},R_p))$.
Finally, pull back the pair of tableaux $(s(P'),Q')$
under column insertion to obtain a word which turns out to be 
the row reading word of a skew column-strict tableau $V$
of the same shape as $U$.  Then $s_p(Q)$ is obtained
by replacing $U$ by $V$.

The bijections $s_p$ specialize to the automorphisms of
conjugation of Lascoux and Sch\"utzenberger \cite{LS1}
in the case that $R$ consists of single rows.

It is shown in \cite{S3} that the bijections
$\sigma_p$ and $s_p$ define an action of the symmetric group $S_L$
on paths and LR tableaux respectively.  Specifically, for $w\in S_L$
let $w=s_{i_1} s_{i_2} \dots s_{i_N}$ be
any factorization of $w$ into adjacent transpositions
$s_i=(i,i+1)$.  For $b\in \PU_R$, define
$w b = \sigma_{i_1} \sigma_{i_2} \dots \sigma_{i_N} b \in \PU_{w R}$.
For $Q\in\LRT(\la;R)$ define $w Q =
s_{i_1} s_{i_2} \dots s_{i_N} Q \in \LRT(\la;w R)$.

\subsection{Generalized charge}

The generalized charge on $Q\in\LRT(\la;R)$ is defined by~\cite{S,SW}
\begin{equation}\label{charge} 
  \charge_R(Q)=\frac{1}{L!}\sum_{w \in S_L}\sum_{i=1}^{L-1}
    (L-i) d_{i,w R}(w Q).
\end{equation}
where $d_{i,R}(Q)=d(P(\word(Q|_{A_i \cup A_{i+1}})))$ where
$d$ is understood to be the function
$d:\LRT(\rho;(R_i,R_{i+1}))\rightarrow\N$.

It was shown in \cite[Section 6]{SW} and \cite{S} that 
$\LRT(R)=\cup_{\la} \LRT(\la;R)$ has the structure of a graded poset
with covering relation given by the $R$-cocyclage and grading function
given by the generalized charge.  The generalized Kostka polynomial is
by definition the generating function of LR tableaux
with the charge statistic~\cite{SW,S}
\begin{equation}\label{gK}
K_{\la R}(q)=\sum_{T\in\LRT(\la;R)} q^{\charge_R(T)}.
\end{equation}
This extends the charge representation of the Kostka polynomial
$K_{\la\mu}(q)$ of Lascoux and Sch\"utzenberger~\cite{LS,LS1}.

For a path $b\in \PU_R$ one has $E(b)=\charge_R(Q(b))$ \cite{S3},
so the formulas \eqref{K path} and \eqref{gK} are equivalent.

\subsection{Embeddings of LR tableaux}
\label{sec LR embed}

The embeddings \eqref{path embed} of sets of paths,
induce embeddings
\begin{equation} \label{LR embed}
  i_R:\LRT(\la;R)\hookrightarrow\LRT(\la;\rows(R))
\end{equation}
via RSK.  These maps are defined in \cite{SW, S2}.  In
the notation of \cite[Section 8.4]{KSS} they are denoted
$\theta^{\rows(R)}_R$.
They are given by compositions of the generalized
automorphisms of conjugation $s_p$ and by the embeddings
of the form $\il_R:\LRT(\la;R)\rightarrow\LRT(\la;\Rl)$
(which is just the inclusion map).  These embeddings preserve
the $R$-cocyclage poset structure and the generalized charge,
since they are induced by maps that preserve the $U_q(\pslnhat)$-crystal 
graph structure.

\subsection{Level-restricted LR tableaux}

Say that a tableau $Q\in\LRT(\la;R)$ is restricted of level $\ell$
if there is a level-restricted path $b\in\PU^\ell_{\la R}$ such that
$Q=Q(b)$.  Denote the set of such tableaux by $\LRT^\ell(\la;R)$.

\begin{ex} \label{lev CST} Suppose
each rectangle is a single row so that
$\LRT(\la;R)=\CST(\la;\mu)$.  In this case let us write
$\CST^\ell(\la;\mu)=\LRT^\ell(\la;R)$.  The following
explicit rule appears in \cite{GW}.  Let $Q\in \CST(\la;\mu)$.
The tableau $Q$ may be viewed as a sequence of shapes
$\emptyset=\la^{(0)} \subset \la^{(1)} \subset \dots \subset\la^{(L)}=\la$
where $\la^{(j)}$ is the shape of $Q|_{[1,j]}$.
Then $Q$ is restricted of level $\ell$ if
\begin{equation} \label{CST lev}
  \la^{(j)}_1 - \la^{(j-1)}_n \le \ell \qquad\text{ for all $1\le j\le L$.}
\end{equation}
In the further special case that $R_j=(1)$ for all $j$,
write $\ST(\la)=\LRT(\la;R)$ for the set of standard tableaux of
shape $\la$ and write $\ST^\ell(\la)=\LRT^\ell(\la;R)$ for the 
level-restricted subset.  For $S\in\ST(\la)$, associate the chain of shapes
$\la^{(j)}$ as above.  Since passing from $\la^{(j-1)}$ to
$\la^{(j)}$ adds only one additional cell, the condition \eqref{CST lev}
simplifies to
\begin{equation} \label{ST lev}
  \la^{(j)}_1 - \la^{(j)}_n \le \ell \qquad\text{ for all $1\le j\le L$.}
\end{equation}
\end{ex}

For general $R$ it is possible to transfer the condition of 
level-restriction on paths to an explicit condition on LR tableaux.  However
for our purposes it is more convenient to use the following description
of $\LRT^\ell(\la;R)$.  Since the embedding \eqref{LR embed}
is induced by the embedding \eqref{path embed} that
preserves the $U_q(\pslnhat)$-crystal graph structure, it follows that
\begin{equation} \label{lev LR}
  \LRT^\ell(\la;R) = \{Q\in\LRT(\la;R) \mid i_R(Q)\in 
    \CST^\ell(\la;\rows(R))\}.
\end{equation}
Hence an expression for the level-restricted generalized Kostka
polynomials equivalent to \eqref{K path lev} is
\begin{equation*}
K_{\la R}^\ell(q)=\sum_{T\in\LRT^\ell(\la;R)}q^{c_R(T)}.
\end{equation*}

\section{Rigged configurations}
\label{sec rc}

This section follows~\cite[Section 2.2]{KSS}, with the notational
difference that here $R_j$ is a rectangle with $\mu_j$
columns and $\eta_j$ rows.  The reason for this is that
here we work with $\RC(\la;R)$ rather than
$\RC(\la^t;R^t)$ as in \cite{KSS}.

\subsection{Review of definitions}

A $(\la;R)$-configuration is a sequence of partitions
$\nu = (\nu^{(1)},\nu^{(2)},\dots)$ with the size constraints
\begin{equation} \label{R config def}
  |\nu^{(k)}| = \sum_{j>k} \la_j -
  	\sum_{a=1}^L \mu_a \max\{\eta_a-k,0\}
\end{equation}
for $k\ge 0$ where by convention $\nu^{(0)}$ is the empty partition.
If $\la$ has at most $n$ parts all partitions $\nu^{(k)}$ for $k\ge n$
are empty.
For a partition $\rho$, define $m_i(\rho)$ to be the number
of parts equal to $i$ and
\begin{equation*}
Q_i(\rho)=\rho^t_1+\rho^t_2+\dots+\rho^t_i=\sum_{j\ge1} \min\{i,\rho_j\},
\end{equation*}
the size of the first $i$ columns of $\rho$.
Let $\xi^{(k)}(R)$ be the partition whose parts are the widths
of the rectangles in $R$ of height $k$.
The vacancy numbers for the $(\la;R)$-configuration $\nu$
are the numbers (indexed by $k \ge 1$ and $i\ge 0$) defined by
\begin{equation} \label{R vacancy def}
  P^{(k)}_i(\nu) = Q_i(\nu^{(k-1)}) - 2 Q_i(\nu^{(k)}) +
   Q_i(\nu^{(k+1)}) + Q_i(\xi^{(k)}(R)).
\end{equation}
In particular $P^{(k)}_0(\nu) = 0$ for all $k\ge 1$.
The $(\la;R)$-configuration $\nu$ is said to be
admissible if $P^{(k)}_i(\nu) \ge 0$ for all $k,i \ge 1$, and
the set of admissible $(\la;R)$-configurations is denoted
by $\Conf(\la;R)$.  Following \cite[(3.2)]{KS}, set
\begin{equation*}
  \cc(\nu)=\sum_{k,i\ge 1} \alpha_i^{(k)}(\alpha_i^{(k)}-\alpha_i^{(k+1)})
\end{equation*}
where $\alpha_i^{(k)}$ is the size of the $i$-th column in
$\nu^{(k)}$.  Define the charge $\charge(\nu)$ of a configuration
$\nu\in\Conf(\la;R)$ by
\begin{equation*}
  \charge(\nu) = ||R|| - \cc(\nu) - |P| 
\end{equation*}
\begin{equation*}
\text{with}\quad
  ||R|| = \sum_{1\le i<j\le L} |R_i\cap R_j|
\quad\text{and}\quad
  |P| = \sum_{k,i\ge1} m_i(\nu) P_i^{(k)}(\nu).
\end{equation*}
Observe that $\charge(\nu)$ depends on both $\nu$ and $R$ but
$\cc(\nu)$ depends only on $\nu$.

\begin{ex}\label{ex config}
Let $\la=(3,2,2,1)$ and $R=((2),(2,2),(1,1))$. Then 
$\nu=((2),(2,1),(1))$ is a $(\la;R)$-configuration with
$\xi^{(1)}(R)=(2)$ and $\xi^{(2)}(R)=(2,1)$. The configuration
$\nu$ may be represented as
\begin{equation*}
\begin{array}{r|c|c|} \cline{2-3} 1&&\\ \cline{2-3} \end{array}
\qquad
\begin{array}{r|c|c|} \cline{2-3} 0&&\\ \cline{2-3}
      0&&\multicolumn{1}{l}{}\\ \cline{2-2} \end{array}
\qquad
\begin{array}{r|c|} \cline{2-2} 0&\\ \cline{2-2} \end{array}
\end{equation*}
where the vacancy numbers are indicated to the left of each part.
In addition $\cc(\nu)=3$, $\|R\|=5$, $|P|=1$ and $\charge(\nu)=1$.
\end{ex}

Define the $q$-binomial by
\begin{equation*}
\qbin{m+p}{m}=\frac{(q)_{m+p}}{(q)_m(q)_p}
\end{equation*}
for $m,p\in\N$ and zero otherwise where $(q)_m=(1-q)(1-q^2)\cdots
(1-q^m)$.  The following fermionic or quasi-particle expression of the 
generalized Kostka polynomials, is a variant of \cite[Theorem 2.10]{KSS}.

\begin{theorem}\label{theo_qp}
For $\la$ a partition and $R$ a sequence of rectangles
\begin{equation}\label{qp}
  K_{\la R}(q)=\sum_{\nu\in\Conf(\la;R)} q^{\charge(\nu)}
   \prod_{k,i\ge 1} \qbin{P_i^{(k)}(\nu)+m_i(\nu^{(k)})}{m_i(\nu^{(k)})}.
\end{equation}
\end{theorem}

Expression \eqref{qp} can be reformulated as the generating function
over rigged configurations. To this end we need
to define certain labelings of the rows of the
partitions in a configuration.
For this purpose one should view a partition as
a multiset of positive integers.
A rigged partition is by definition a finite multiset of
pairs $(i,x)$ where $i$ is a positive integer and
$x$ is a nonnegative integer.  The pairs $(i,x)$ are referred to
as strings; $i$ is referred to as the
length of the string and $x$ as the label or
quantum number of the string.  A rigged partition is
said to be a rigging of the partition $\rho$ if
the multiset consisting of the lengths of the strings,
is the partition $\rho$.  So a rigging of $\rho$
is a labeling of the parts of $\rho$ by nonnegative integers,
where one identifies labelings that differ only by
permuting labels among equal-sized parts of $\rho$.

A rigging $J$ of the $(\la;R)$-configuration $\nu$ is
a sequence of riggings of the partitions $\nu^{(k)}$ such that
for every part of $\nu^{(k)}$ of length $i$ and label $x$,
\begin{equation} \label{rigging def}
  0 \le x \le P^{(k)}_i(\nu).
\end{equation}
The pair $(\nu,J)$ is called a rigged configuration.
The set of riggings of admissible $(\la;R)$-configurations
is denoted by $\RC(\la;R)$.
Let $(\nu,J)^{(k)}$ be the $k$-th rigged partition
of $(\nu,J)$.  A string $(i,x)\in (\nu,J)^{(k)}$
is said to be singular if $x=P^{(k)}_i(\nu)$, that is,
its label takes on the maximum value.

Observe that the definition of the set $\RC(\la;R)$ is completely
insensitive to the order of the rectangles in the sequence $R$.
However the notation involving the sequence $R$
is useful when discussing the bijection between LR tableaux
and rigged configurations,
since the ordering on $R$ is essential in the definition of LR tableaux.

Define the cocharge and charge of $(\nu,J)\in\RC(\la;R)$ by
\begin{equation*}
\begin{split}
  \cc(\nu,J)&=\cc(\nu)+|J| \\
  \charge(\nu,J)&=\charge(\nu)+|J| \\
  |J| &= \sum_{k,i\ge 1} |J_i^{(k)}|
\end{split}
\end{equation*}
where $J_i^{(k)}$ is the partition inside the rectangle
of height $m_i(\nu^{(k)})$ and width $P_i^{(k)}(\nu)$ given
by the labels of the parts of $\nu^{(k)}$ of size $i$.

Since the $q$-binomial $\qbins{m+p}{m}$ is the generating
function of partitions with at most $m$ parts each not
exceeding $p$ \cite[Theorem 3.1]{Andrews76}, Theorem \ref{theo_qp} is 
equivalent to the following theorem.

\begin{theorem}\label{theo_rc}
For $\la$ a partition and $R$ a sequence of rectangles
\begin{equation}\label{rc}
  K_{\la R}(q)=\sum_{(\nu,J)\in\RC(\la;R)} q^{\charge(\nu,J)}.
\end{equation}
\end{theorem}

\subsection{Switching between quantum and coquantum numbers}

Let $\theta_R:\RC(\la;R)\rightarrow\RC(\la;R)$
be the involution that complements quantum numbers.  More precisely,
for $(\nu,J)\in\RC(\la;R)$, replace every string
$(i,x)\in (\nu,J)^{(k)}$ by $(i,P_i^{(k)}(\nu)-x)$.
The notation here differs from that in \cite{KSS}, in which
$\theta_R$ is an involution on $\RC(\la^t;R^t)$.

\begin{lemma} \label{theta and c}
$\charge(\theta_R(\nu,J))=||R||-\cc(\nu,J)$ for all $(\nu,J)\in\RC(\la;R)$.
\end{lemma}
\begin{proof} Let $\theta_R(\nu,J)=(\nu',J')$.  It follows immediately
from the definitions that $\nu'=\nu$.  In particular $\nu$ and $\nu'$
have the same vacancy numbers and $|J'|=|P|-|J|$. Then
\begin{equation*}
\begin{split}
  \charge(\theta_R(\nu,J)) &= \charge(\nu',J')= ||R||-\cc(\nu')-|P|+|J'| \\
  &= ||R||-\cc(\nu)-|J|=||R||-\cc(\nu,J).
\end{split}
\end{equation*}
\end{proof}

There is a bijection $\RCtr:\RC(\la;R)\rightarrow\RC(\la^t;R^t)$
that has the property
\begin{equation} \label{RCtr and cc}
 \cc(\RCtr(\nu,J))=||R||-\cc(\nu,J)
\end{equation}
for all $(\nu,J)\in\RC(\la;R)$;
see the proof of \cite[Proposition 11]{KS}.

\subsection{RC's and level-restriction}
\label{sec lev rc}

Here we introduce the most important new definition in this paper,
namely, that of a level-restricted rigged configuration.

Say that a partition $\la$ is restricted of level $\ell$
if $\la_1-\la_n\le \ell$, recalling that it is assumed that all
partitions have at most $n$ parts, some of which may be zero.
Fix a shape $\la$ and a sequence of rectangles $R$
that are all restricted of level $\ell$.  Define
$\lt = \ell-(\la_1-\la_n)$, which is nonnegative by assumption.

Set $\la'=(\la_1-\la_n,\ldots,\la_{n-1}-\la_n)^t$ and 
denote the set of all column-strict tableaux of shape $\la'$ over the alphabet
$\{1,2,\ldots,\la_1-\la_n\}$ by $\CST(\la')$. Define a table of modified
vacancy numbers depending on $\nu\in\Conf(\la;R)$ and $t\in\CST(\la')$ by
\begin{equation} \label{t vacancy}
  P_i^{(k)}(\nu,t) =
  P_i^{(k)}(\nu) - \sum_{j=1}^{\la_k-\la_n} \chi(i\ge\lt+t_{j,k})
+ \sum_{j=1}^{\la_{k+1}-\la_n} \chi(i\ge\lt+t_{j,k+1})
\end{equation}
for all $i,k\ge1$, 
where $\chi(S)=1$ if the statement $S$ is true and
$\chi(S)=0$ otherwise, and $t_{j,k}$ is the $(j,k)$-th entry of $t$.
Finally let $\x_i^{(k)}$ be the largest part of the partition
$J_i^{(k)}$; if $J_i^{(k)}$ is empty set $x_i^{(k)}=0$.

\begin{definition}\label{def levrc}
Say that $(\nu,J)\in\RC(\la;R)$ is restricted of level $\ell$
provided that
\begin{enumerate}
\item $\nu_1^{(k)} \le \ell$ for all $k$.
\item There exists a tableau $t\in\CST(\la')$, such that for every 
$i,k\ge 1$,
\begin{equation*}
  \x_i^{(k)} \le P_i^{(k)}(\nu,t).
\end{equation*}
\end{enumerate}
Let $\Conf^\ell(\la;R)$ be the set of all $\nu\in \Conf(\la;R)$ such that
the first condition holds, and denote by $\RC^\ell(\la;R)$ the set of 
$(\nu,J)\in\RC(\la;R)$ that are restricted of level $\ell$. 
\end{definition}
Note in particular that the second condition requires that
$P_i^{(k)}(\nu,t)\ge 0$ for all $i,k\ge 1$.

\begin{ex}\label{ex kir}
Let us consider Definition \ref{def levrc} for two classes of shapes $\la$
more closely:
\begin{enumerate}
\item \label{e rect}
Vacuum case: Let $\la=(a^n)$ be rectangular with $n$ rows.
Then $\la'=\emptyset$ and $P_i^{(k)}(\nu,\emptyset)=P_i^{(k)}(\nu)$
for all $i,k\ge 1$ so that the modified vacancy numbers are equal
to the vacancy numbers.
\item Two-corner case: Let $\la=(a^\alpha,b^\beta)$ with $\alpha+\beta=n$
and $a>b$. Then $\la'=(\alpha^{a-b})$ and there is only one tableau $t$
in $\CST(\la')$, namely the Yamanouchi tableau of shape $\la'$.
Since $t_{j,k}=j$ for $1\le k\le \alpha$ we find that
\begin{equation*}
P_i^{(k)}(\nu,t)=P_i^{(k)}(\nu)-\delta_{k,\alpha}\max\{i-\lt,0\}
\end{equation*}
for $1\le i\le \ell$ and $1\le k<n$. We wish to thank Anatol Kirillov 
for communicating this formula to us~\cite{K98}.
\end{enumerate}
\end{ex}

Our main result is the following
formula for the level-restricted generalized Kostka polynomial:
\begin{theorem} \label{main}
Let $\ell$ be a positive integer. For $\la$ a partition and $R$
a sequence of rectangles both restricted of level $\ell$,
\begin{equation*}
  K^\ell_{\la R}(q) = \sum_{(\nu,J)\in\RC^\ell(\la;R)} q^{\charge(\nu,J)}.
\end{equation*}
\end{theorem}
The proof of this theorem is given in Section \ref{sec proof}.

\begin{ex} Consider $n=3$, $\ell=2$, $\la=(3,2,1)$ and $R=((2),(1)^4)$.
Then
\begin{equation}\label{conf ex}
\begin{array}{r|c|} \cline{2-2} 0&\\ \cline{2-2} 0&\\ \cline{2-2} 0&\\
 \cline{2-2} \end{array}
\quad
\begin{array}{r|c|} \cline{2-2} 1&\\ \cline{2-2} \end{array}
\qquad \text{and} \qquad
\begin{array}{r|c|c|} \cline{2-3} 1&&\\ \cline{2-3}
      2&&\multicolumn{1}{l}{}\\ \cline{2-2} \end{array}
\quad
\begin{array}{r|c|} \cline{2-2} 0&\\ \cline{2-2} \end{array}
\end{equation}
are in $\Conf^\ell(\la;R)$ where again the vacancy numbers are indicated 
to the left of each part.
The set $\CST(\la')$ consists of the two elements
\begin{equation*}
\begin{array}{|c|c|} \cline{1-2} 1&1\\ \cline{1-2}
      2&\multicolumn{1}{l}{}\\ \cline{1-1} \end{array}
\qquad\text{and}\qquad
\begin{array}{|c|c|} \cline{1-2} 1&2\\ \cline{1-2}
      2&\multicolumn{1}{l}{}\\ \cline{1-1} \end{array}\;.
\end{equation*}
Since $\lt=0$ the three rigged configurations
\begin{equation*}
\begin{array}{|c|l} \cline{1-1} &0\\ \cline{1-1} &0\\ \cline{1-1} &0\\
 \cline{1-1} \end{array}
\quad
\begin{array}{|c|l} \cline{1-1} &0\\ \cline{1-1} \end{array},
\qquad
\begin{array}{|c|c|l} \cline{1-2} &&0\\ \cline{1-2}
      &\multicolumn{2}{l}{0}\\ \cline{1-1} \end{array}
\quad
\begin{array}{|c|l} \cline{1-1} &0\\ \cline{1-1} \end{array}
\qquad\text{and}\qquad
\begin{array}{|c|c|l} \cline{1-2} &&0\\ \cline{1-2}
      &\multicolumn{2}{l}{1}\\ \cline{1-1} \end{array}
\quad
\begin{array}{|c|l} \cline{1-1} &0\\ \cline{1-1} \end{array}
\end{equation*}
are restricted of level 2 with charges $2,3,4$, respectively.
The riggings are given on the right of each part. 
Hence $K_{\la R}^2(q)=q^2+q^3+q^4$.

In contrast to this, the Kostka polynomial $K_{\la\mu}(q)$
is obtained by summing over both configurations in \eqref{conf ex}
with all possible riggings below the vacancy numbers.
This amounts to $K_{\la\mu}(q)=q^2+2q^3+2q^4+2q^5+q^6$.
\end{ex}

In Section~\ref{sec bf} we will use Theorem~\ref{main}
to obtain explicit expressions for type $A$ branching functions.
The results suggest that it is also useful to consider the
following sets of rigged configurations with imposed minima on 
the set of riggings. 

Let $\rho\subset \la$ be a
partition and $R_\rho=((1^{\rho_1^t}),(1^{\rho_2^t}),\ldots,(1^{\rho_n^t}))$,
the sequence of single columns of height $\rho_i^t$.
Set $\rho'=(\rho_1-\rho_n,\ldots,\rho_{n-1}-\rho_n)^t$ and
\begin{equation*}
M_i^{(k)}(t)=\sum_{j=1}^{\rho_k-\rho_n}\chi(i\le \rho_1-\rho_n-t_{j,k})
-\sum_{j=1}^{\rho_{k+1}-\rho_n}\chi(i\le \rho_1-\rho_n-t_{j,k+1})
\end{equation*}
for all $t\in\CST(\rho')$.
Then define $\RC^\ell(\la,\rho;R)$ to be the set of all
$(\nu,J)\in\RC^\ell(\la;R_\rho\cup R)$ such that there exists
a $t\in\CST(\rho')$ such that $M_i^{(k)}(t)\le x$
for $(i,x)\in (\nu,J)^{(k)}$ and $M_i^{(k)}(t)\le P_i^{(k)}(\nu)$
for all $i,k\ge 1$. Note that the second condition is obsolete
if $i$ occurs as a part in $\nu^{(k)}$ since by definition
$M_i^{(k)}(t)\le x\le P_i^{(k)}(\nu)$ for all $(i,x)\in (\nu,J)^{(k)}$.

Conjecture~\ref{conj skew} asserts that the set
$\RC^\ell(\la,\rho;R)$ corresponds to the set of all level-$\ell$
restricted Littlewood--Richardson tableaux with a fixed
subtableaux of shape $\rho$.

\section{Fermionic expression of level-restricted generalized
Kostka polynomials} \label{sec qp Kostka}

\subsection{Fermionic expression}
Similarly to the Kostka polynomial case, one can rewrite the expression
of the level-restricted generalized Kostka polynomials of Theorem \ref{main}
in fermionic form.

\begin{lemma}\label{lem Pl}
For all $\nu\in\Conf^\ell(\la,R)$, $t\in\CST(\la')$ and $1\le k<n$, 
we have $P_i^{(k)}(\nu,t)=0$ for $i\ge \ell$.
\end{lemma}
\begin{proof}
Since $\nu^{(k)}_1\le \ell$ it follows from \cite[(11.2)]{KS} that
$P_i^{(k)}(\nu)=\la_k-\la_{k+1}$ for $i\ge \ell$.
Since $t$ is over the alphabet $\{1,2,\ldots,\la_1-\la_n\}$ this implies
for $i\ge \ell$
\begin{equation*}
\begin{split}
P_i^{(k)}(\nu,t)=&P_i^{(k)}(\nu)
   -\sum_{j=1}^{\la_k-\la_n}\chi(i\ge \lt+t_{j,k})
   +\sum_{j=1}^{\la_{k+1}-\la_n}\chi(i\ge \lt+t_{j,k+1})\\
=& \la_k-\la_{k+1}-(\la_k-\la_n)+(\la_{k+1}-\la_n)=0.
\end{split}
\end{equation*}
\end{proof}

Let $\SCST(\la')$ be the set of all nonempty subsets of $\CST(\la')$.
Furthermore set $P_i^{(k)}(\nu,S)=\min\{P_i^{(k)}(\nu,t)|t\in S\}$ for
$S\in\SCST(\la')$. Then by inclusion-exclusion the set of allowed
rigging for a given configuration $\nu\in\Conf^\ell(\la;R)$ is given by
\begin{equation*}
\sum_{S\in\SCST(\la')}(-1)^{|S|+1} \{J|\x_i^{(k)}\le P_i^{(k)}(\nu,S)\}.
\end{equation*}
Since the $q$-binomial $\qbins{m+p}{m}$ is the generating function of 
partitions with at most $m$ parts each not exceeding $p$ and since
$P_\ell^{(k)}(\nu,S)=0$ by Lemma \ref{lem Pl}
the level-$\ell$ restricted generalized Kostka polynomials
has the following fermionic form.
\begin{theorem}\label{thm qp}
\begin{equation*}
K_{\la R}^\ell(q)=\sum_{S\in\SCST(\la')} (-1)^{|S|+1}
\sum_{\nu\in\Conf^\ell(\la;R)} q^{\charge(\nu)}
\prod_{i=1}^{\ell-1}\prod_{k=1}^{n-1} 
\qbin{m_i(\nu^{(k)})+P_i^{(k)}(\nu,S)}{m_i(\nu^{(k)})}.
\end{equation*}
\end{theorem}

In Section \ref{sec bf} we will derive new expressions for 
branching functions of type $A$ as limits of the level-restricted
generalized Kostka polynomials. To this end we need to reformulate
the fermionic formula of Theorem \ref{thm qp} in terms of a
so-called $(\vm,\vn)$-system. Set
\begin{align*}
m_i^{(a)}&=P_i^{(a)}(\nu,S)=P_i^{(a)}(\nu)+f_i^{(a)}(S),\\
n_i^{(a)}&=m_i(\nu^{(a)}),
\end{align*}
and $L_i^{(a)}=\sum_{j=1}^L \chi(i=\mu_j)\chi(a=\eta_j)$
for $1\le i\le \ell$ and $1\le a \le n$ which is the number of 
rectangles in $R$ of shape $(i^a)$. Then 
\begin{equation*}
\begin{split}
&-m_{i-1}^{(a)}+2m_i^{(a)}-m_{i+1}^{(a)}
-n_i^{(a-1)}+2n_i^{(a)}-n_i^{(a+1)}\\
=&(\alpha_i^{(a-1)}-2\alpha_i^{(a)}+\alpha_i^{(a+1)})
-(\alpha_{i+1}^{(a-1)}-2\alpha_{i+1}^{(a)}+\alpha_{i+1}^{(a+1)})\\
&+\sum_{k=1}^L \delta_{a,\eta_k}(-\min\{i-1,\mu_k\}+2\min\{i,\mu_k\}
-\min\{i+1,\mu_k\})\\
&-f_{i-1}^{(a)}(S)+2f_i^{(a)}(S)-f_{i+1}^{(a)}(S)\\
&-(\alpha_i^{(a-1)}-\alpha_{i+1}^{(a-1)})+2(\alpha_i^{(a)}-\alpha_{i+1}^{(a)})
-(\alpha_i^{(a+1)}-\alpha_{i+1}^{(a+1)})\\
=&L_i^{(a)}-f_{i-1}^{(a)}(S)+2f_i^{(a)}(S)-f_{i+1}^{(a)}(S).
\end{split}
\end{equation*}
At this stage it is convenient to introduce vector notation.
For a matrix $v_i^{(a)}$ with indices $1\le i\le \ell-1$ and 
$1\le a\le n-1$ define
\begin{equation*}
\vv=\sum_{i=1}^{\ell-1}\sum_{a=1}^{n-1}v_i^{(a)}\ve_i\otimes\ve_a,
\end{equation*}
where $\ve_i$ and $\ve_a$ are the canonical basis vectors of $\Z^{\ell-1}$
and $\Z^{n-1}$, respectively. Define
\begin{equation*}
u_i^{(a)}(S)=-f_{i-1}^{(a)}(S)+2f_i^{(a)}(S)-f_{i+1}^{(a)}(S)
\end{equation*}
which in vector notation reads
\begin{equation}\label{u}
\vu(S)=(C\otimes I)\vf(S)+\sum_{a=1}^{n-1}(\la_a-\la_{a+1})
 \ve_{\ell-1}\otimes\ve_a,
\end{equation}
where $C$ is the Cartan matrix of type $A$ and $I$ is the identity
matrix. Since $n_i^{(0)}=n_i^{(n)}=m_0^{(k)}=0$ and 
$m_\ell^{(k)}=0$ by Lemma \ref{lem Pl} it follows that
\begin{equation}\label{mn}
(C\otimes I) \vm + (I\otimes C) \vn = \vL + \vu(S).
\end{equation}
In terms of the new variables the condition \eqref{R config def} 
on $|\nu^{(a)}|$ becomes
\begin{equation}\label{nl}
n_\ell^{(a)}=-\ve_{\ell-1}\otimes\ve_a ( C^{-1}\otimes I ) \vn
 -\frac{1}{\ell}\sum_{j=1}^a\la_j
 +\frac{1}{\ell}\sum_{i=1}^\ell \sum_{b=1}^n i \min\{a,b\} L_i^{(b)},
\end{equation}
where we used $C_{ij}^{-1}=\min\{i,j\}-ij/\ell$ if $C$ is 
$(\ell-1)\times (\ell-1)$-dimensional and
$\sum_{b=1}^n \sum_{i=1}^\ell ibL_i^{(b)}=|\la|$.

\begin{lemma}\label{lem charge}
In terms of the above $(\vm,\vn)$-system
\begin{multline}
\charge(\nu)=\frac{1}{2}\vm (C\otimes C^{-1})\vm
-\vm(I\otimes C^{-1})\vu(S)\\
+\frac{1}{2}\vu(S)(C^{-1}\otimes C^{-1})\vu(S)+g(R,\la)
\end{multline}
where
\begin{equation*}
g(R,\la)=\|R\|
 -\frac{1}{2}\sum_{a,b=1}^{n-1}\sum_{j=1}^\ell C_{ab}^{-1}
   L_j^{(a)} \Lb_j^{(b)}
 +\frac{1}{2\ell}\sum_{j=1}^n (\la_j-\frac{1}{n}|\la|)^2
\end{equation*}
and $\Lb_i^{(a)}=\sum_{j=1}^\ell \min\{i,j\}L_j^{(a)}$.
\end{lemma}

\begin{proof}
By definition $\charge(\nu)=\|R\|-\cc(\nu)-|P|$. Note that
\begin{equation*}
\begin{split}
|P|&= \sum_{i=1}^{\ell}\sum_{k=1}^{n-1} m_i(\nu^{(k)})P_i^{(k)}(\nu)\\
&= \sum_{i=1}^\ell\sum_{k=1}^{n-1} (\alpha_i^{(k)}-\alpha_{i+1}^{(k)})
(\sum_{j=1}^i(\alpha_j^{(k-1)}-2\alpha_j^{(k)}+\alpha_j^{(k+1)})
+\Lb_i^{(k)})\\
&= -2\cc(\nu)+\sum_{i=1}^\ell\sum_{k=1}^{n-1} n_i^{(k)}\Lb_i^{(k)}.
\end{split}
\end{equation*}
Hence eliminating $\cc(\nu)$ in favor of $|P|$ yields
\begin{equation*}
\charge(\nu)=\|R\|-\frac{1}{2}|P|-\frac{1}{2}\sum_{i=1}^\ell
\sum_{k=1}^{n-1} n_i^{(k)} \Lb_i^{(k)}.
\end{equation*}
On the other hand, using $n_i^{(k)}=m_i(\nu^{(k)})$ and 
$P_\ell^{(k)}(\nu)=\la_k-\la_{k+1}$,
\begin{equation*}
|P|=\vn (I\otimes I) \vP(\nu)+\sum_{k=1}^{n-1} n_\ell^{(k)}(\la_k-\la_{k+1})
\end{equation*}
so that
\begin{equation}\label{charge inter}
\charge(\nu)=\|R\|-\frac{1}{2}\vn (I\otimes I) (\vP(\nu)+\vLb)
-\frac{1}{2}\sum_{k=1}^{n-1}n_\ell^{(k)}(\la_k-\la_{k+1}+\Lb_\ell^{(k)}).
\end{equation}
Eliminating $\vn$ in favor of $\vm$ using \eqref{mn} and substituting
$\vP(\nu)=\vm-\vf(S)$ yields
\begin{multline*}
-\frac{1}{2}\vn(I\otimes I)(\vP(\nu)+\vLb)=
 \frac{1}{2}\vm\{C\otimes C^{-1}(\vm+\vLb-\vf(S))
 -I\otimes C^{-1}(\vL+\vu(S))\}\\
 -\frac{1}{2}(\vL+\vu(S))(I\otimes C^{-1})(\vLb-\vf(S)).
\end{multline*}
Similarly, replacing $\vn$ by $\vm$ in \eqref{nl} we obtain
\begin{multline}\label{nl m}
n_\ell^{(a)}= \ve_{\ell-1}\otimes \ve_a (I\otimes C^{-1}\vm
 -C^{-1}\otimes C^{-1}\vu(S))\\
 -\frac{1}{\ell}\sum_{j=1}^a(\la_j-\frac{1}{n}|\la|)
 +\sum_{b=1}^{n-1} C_{ab}^{-1} L_\ell^{(b)}.
\end{multline}
Inserting these equations into \eqref{charge inter}, trading
$\vf(S)$ for $\vu(S)$ by \eqref{u} and using
\begin{equation*}
(C\otimes I)\vLb-\vL-\sum_{a=1}^{n-1}\ve_{\ell-1}\otimes \ve_a \Lb_\ell^{(a)}
=0
\end{equation*}
results in the claim of the lemma.
\end{proof}

As a corollary of Lemma \ref{lem charge} and Theorem \ref{thm qp}
we obtain the following expression for the level-restricted generalized 
Kostka polynomial
\begin{multline}\label{Kmn}
K_{\la R}^\ell(q)=q^{g(R,\la)}
 \sum_{S\in\SCST(\la')} (-1)^{|S|+1} 
 q^{\frac{1}{2}\vu(S)C^{-1}\otimes C^{-1}\vu(S)}\\
\times
\sum_{\vm} q^{\frac{1}{2}\vm C\otimes C^{-1}\vm-\vm I\otimes C^{-1} \vu(S)}
\qbin{\vm+\vn}{\vm}
\end{multline}
where $\vn$ is determined by \eqref{mn}, the sum over $\vm$ is such
that 
\begin{multline*}
\ve_{\ell-1}\otimes \ve_a (I\otimes C^{-1}\vm-C^{-1}\otimes C^{-1}\vu(S))\\
 -\frac{1}{\ell}\sum_{j=1}^a(\la_j-\frac{1}{n}|\la|)
 +\sum_{b=1}^{n-1} C_{ab}^{-1}L_\ell^{(b)}\in\Z,
\end{multline*}
for all $1\le a\le n-1$ 
and $\qbins{\vm+\vn}{\vm}=\prod_{i=1}^{\ell-1}\prod_{k=1}^{n-1}
\qbins{m_i^{(k)}+n_i^{(k)}}{m_i^{(k)}}$.

Now consider the second case of Example \ref{ex kir}, namely
$\la=(a^\alpha,b^\beta)$ with $a>b$ and $\alpha+\beta=n$.
Then $\SCST(\la')$ only contains the element $S=\{t\}$ where $t$ is 
the Yamanouchi tableau of shape $\la'$ and 
$\vu(S)=\ve_{\lt}\otimes\ve_\alpha$.
In the vacuum case, that is, when $\la=((\frac{|\la|}{n})^n)$, the set 
$\SCST(\la')$ only contains $S=\{\emptyset\}$ and $\vu(S)=\vf(S)=0$.
In this case \eqref{Kmn} simplifies to
\begin{equation*}
K_{\la R}^\ell(q)=q^{g(R,\la)}
 \sum_{\vm} q^{\frac{1}{2}\vm C\otimes C^{-1}\vm}\qbin{\vm+\vn}{\vm}.
\end{equation*}
When $R$ is a sequence of single boxes this proves 
\cite[Theorem 1]{Das97}\footnote{We believe that the proof given
in \cite{Das97} is incomplete.}. When $R$ is a sequence of single
rows or single columns this settles \cite[Conjecture 4.7]{HKKOTY}.

\subsection{Polynomial Rogers--Ramanujan-type identities}

Let $\Wb$ be the Weyl group of $\sln$,
$M=\{\beta\in\Z^n| \sum_{i=1}^n \beta_i=0\}$ be the root lattice,
$\rho$ the half-sum of the positive roots,
and $(\cdot|\cdot)$ the standard symmetric bilinear form.
Recall the energy function \eqref{energy}. It was shown in \cite{SS} that
\begin{equation}\label{Kbose}
K_{\la R}^\ell(q)=\sum_{\tau\in\Wb}\sum_{\beta\in M}
\sum_{\substack{b\in \PU_R \\
	\wt(b)=-\rho+\tau^{-1}(\lab-(\ell+n)\beta+\rho)}} (-1)^\tau 
q^{-\frac{1}{2}(\ell+n)(\beta|\beta)+(\lab+\rho|\beta)+E(b)}.
\end{equation}
Equating \eqref{Kmn} and \eqref{Kbose} gives rise to 
polynomial Rogers--Ramanujan-type identities.
For the vacuum case, that is, when the partition $\la$ is rectangular
with $n$ rows, this proves \cite[Eq. (9.2)]{SW}\footnote{The 
definition of level-restricted path as given in \cite[p. 394]{SW}
only works when $R$ (or $\mu$ therein) consists of single rows; 
otherwise the description of Section \ref{sec restr. paths}
should be used.}. 

\section{New expressions for type $A$ branching functions}
\label{sec bf}

The coset branching functions $b_{\La'\La''}^\La$
labeled by the three weights $\La,\La',\La''$ have a natural finitization
in terms of $(\La'+\La'')$-restricted crystals. For certain
triples of weights these can be reformulated in terms of level-restricted
paths, which in turn yield an expression of the type $A$ branching
functions as a limit of the level-restricted generalized Kostka
polynomials. Together with the results of the last section
this implies new fermionic expressions for type $A$
branching functions at certain triples of weights.

\subsection{Branching function in terms of paths}

Let $\La,\La',\La''\in\Pcl$ be dominant integral weights of levels
$\ell$, $\ell'$, and $\ell''$ respectively, where $\ell=\ell'+\ell''$.
The branching function $b^\La_{\La'\La''}(z)$ is the
formal power series defined by
\begin{equation*}
  b^\La_{\La'\La''}(z) = \sum_{m\ge0} z^m
  	c^{\af(\La)-m\delta}_{\af(\La'),\af(\La'')}
\end{equation*}
where $c^{\af(\La)-m\delta}_{\af(\La'),\af(\La'')}$
is the multiplicity of the irreducible integrable highest
weight $U_q(\slnhat)$-module $\VV(\af(\La)-m\delta)$
in the tensor product $\VV(\af(\La'))\otimes \VV(\af(\La''))$.

The desired multiplicity is equal to the number of
$\slnhat$-highest weight vectors of weight $\af(\La)-m\delta$
in the tensor product $\BB(\af(\La'))\otimes \BB(\af(\La''))$, that
is, the number of elements
$b'\otimes b''\in \BB(\af(\La'))\otimes \BB(\af(\La''))$ such that
$\wt(b'\otimes b'')=\af(\La)-m\delta$ and
$\epsilon_i(b'\otimes b'')=0$ for all $i\in I$.
By \eqref{tensor crystal}, $b''=u_{\La''}$,
$b'$ is $\La''$-restricted, and
$\wt(b') = \af(\La-\La'')-m\delta$.

Let $\Bc$ be a perfect crystal of level $\ell'$.
Using the isomorphism $\BB(\La')\cong \PP(\La',\Bc)$
let $b'=b'_1\otimes b'_2 \otimes \dotsm$
and $\bb\in\PP(\La',\Bc)$ be the ground state path.
Suppose $N$ is such that for all $j>N$, $b'_j = \bb_j$.
Write $b=b'_1\otimes \dots \otimes b'_N$.
In type $A_{n-1}^{(1)}$ the period of the ground state path
$\bb$ always divides $n$.
Choose $N$ to be a multiple of $n$, so that
$b' = b \otimes \bb$ and $\bb_{N+1}=\bb_1$.

Then the above desired highest weight vectors have the form
$b' \otimes b'' = (b \otimes u_{\La'}) \otimes u_{\La''}
\in \Bc^{\otimes N} \otimes \BB(\af(\La')) \otimes \BB(\af(\La''))$.
But there is an embedding
$\BB(\af(\La'+\La''))\hookrightarrow \BB(\af(\La'))\otimes \BB(\af(\La''))$
defined by $u_{\La'+\La''}\rightarrow u_{\La'} \otimes u_{\La''}$.
With this rephrasing of the conditions on $b$
and taking limits, we have
\begin{equation}\label{branch limit}
  b^\La_{\La'\La''}(z) =
  	\lim_{\substack{N\rightarrow\infty \\ N\in n\Z}}
  	z^{-E_N(\bb_1 \otimes \dots \otimes \bb_N)}
  	\sum_{b\in \HH(\La'+\La'',\Bc^{\otimes N},\La)}
		z^{E_N(b)}
\end{equation}
where $E_N:\Bc^{\otimes N}\rightarrow\Z$ is given by
$E_N(b) = E(b \otimes \bb_{N+1})=E(b\otimes \bb_1)$
and $E$ is the energy function on finite paths.

Our goal is to express \eqref{branch limit} in terms of
level-restricted generalized Kostka polynomials. We find
that this is possible for certain triples of weights.
Using the results of Section~\ref{sec qp Kostka} this provides
explicit formulas for the branching functions.

\subsection{Reduction to level-restricted paths}

The first step in the transformation of \eqref{branch limit}
is to replace the condition of $(\La'+\La'')$-restrictedness
by level $\ell$ restrictedness.  This is achieved at the cost
of appending a fixed inhomogeneous path.

Consider any tensor product $\Bc''$ of perfect crystals
each of which has level at most $\ell''$ (the level of $\La''$),
such that there is an element $y''\in\HH(\ell''\La_0,\Bc'',\La'')$.
We indicate how such a $\Bc''$ and $y''$ can be constructed explicitly.
Let $\la$ be the partition with strictly less than $n$ rows
with $\inner{h_i}{\La''}$ columns of length $i$ for $1\le i\le n-1$.
Let $Y_\la$ be the Yamanouchi tableau of shape $\la$.
Then any factorization (in the plactic monoid)
of $Y_\la$ into a sequence of rectangular tableaux, yields such a
$\Bc''$ and $y''$.

\begin{ex} Let $n=6$, $\ell''=5$, $\La''=\La_0+2\La_2+\La_3+\La_4$.
Then $\la=(4,4,2,1)$ (its transpose is $\la^t=(4,3,2,2))$ and
\begin{equation*}
  Y_\la = \begin{matrix}
    1&1&1&1\\
    2&2&2&2\\
    3&3& & \\
    4& & &
  \end{matrix}\,.
\end{equation*}
One way is to factorize into single columns:
$\Bc'' = \Bc^{2,1} \otimes \Bc^{2,1} \otimes \Bc^{3,1} \otimes \Bc^{4,1}$
and $y''=y_4 \otimes y_3 \otimes y_2 \otimes y_1$ where
each $y_j$ is an $\sln$ highest weight vector, namely, the $j$-th column
of $Y_\la$.  Another way is to factorize into the minimum number
of rectangles by slicing $Y_\la$ vertically.  This yields
$\Bc''=\Bc^{2,2} \otimes \Bc^{3,1} \otimes \Bc^{4,1}$;
again the factors of $y''=y_3\otimes y_2 \otimes y_1$
are the $\sln$ highest weight vectors, namely,
\begin{equation*}
  y_3 = \begin{matrix} 1&1 \\ 2&2 \end{matrix}\,, \qquad
  y_2 = \begin{matrix} 1 \\ 2 \\ 3 \end{matrix}\,, \qquad
  y_1 = \begin{matrix} 1 \\ 2 \\ 3 \\ 4 \end{matrix}\,.
\end{equation*}
\end{ex}

Consider also a tensor product $\Bc'$ of perfect crystals
such that there is an element $y'\in\HH(\ell'\La_0,\Bc',\La')$.
Then $y=y' \otimes y'' \in \HH(\ell\La_0,\Bc'\otimes \Bc'',\La'+\La'')$.
Instead of $b\in \HH(\La'+\La'',\Bc^{\otimes N},\La)$,
we work with $b \otimes y$ where
$b\otimes y$ is restricted of level $\ell$.

This trick doesn't help unless one can recover the
correct energy function directly from $b \otimes y$.
Let $p$ be the first $N$ steps of the ground state path
$\bb\in\PP(\La',\Bc)$.
Define the normalized energy function on $\Bc^{\otimes N}$ by
$\Eb(b) = E(b \otimes y') - E(p \otimes y')$.
A priori it depends on $\La'$, $\Bc$, and $y'$.
The energy function occurring in the branching function is
$E'(b)=E(b \otimes \bb_1)-E(p \otimes \bb_1)$.

\begin{lemma} $\Eb=E'$.
\end{lemma}
\begin{proof} It suffices to show that the function
$\Bc^{\otimes N} \rightarrow \Z$ given by
$b\mapsto E(b\otimes y') - E(b \otimes \bb_1)$ is constant.
Using the definition \eqref{energy} and the fact that 
$b$ is homogeneous of length $N$, we have
\begin{equation*}
  E(b \otimes y') = E(b) + N E(b_N \otimes y') - (N-1)E(y').
\end{equation*}
Similarly $E(b\otimes \bb_1) = E(b)+N E(b_N \otimes \bb_1)$.
Therefore $E(b\otimes y')-E(b\otimes \bb_1)=
N(E(b_N \otimes y') -E(b_N \otimes \bb_1)) - (N-1) E(y')$.
Thus it suffices to show that the function $\Bc\rightarrow\Z$
given by $b'\mapsto E(b'\otimes y')-E(b'\otimes\bb_1)$
is a constant function.

Suppose first that
$\epsilon_i(b') > \inner{h_i}{\La'}$ for some $1\le i\le n-1$.
By the construction of $y'$ and $\bb_1$,
$\phi_i(y')=\inner{h_i}{\La'}=\phi_i(\bb_1)$
for $1\le i\le n-1$, since $\phi(\bb_1)=\La'$.  Then
$e_i(b' \otimes y') = e_i(b') \otimes y'$
and $e_i(b' \otimes \bb_1)=e_i(b') \otimes \bb_1$
by \eqref{tensor raise}.  Passing from $b'$ to $e_i(b')$ repeatedly,
the values of the energy functions are constant, so 
it may be assumed that $b' \otimes y'$ is a $\sln$ highest
weight vector; in particular, $\epsilon_i(b')\le \inner{h_i}{\La'}$
for all $1\le i\le n-1$.

Next suppose that $\epsilon_0(b') > \inner{h_0}{\La'}$.
Now $\phi_0(y')=0$ and $\phi_0(\bb_1)=\inner{h_0}{\La'}$.
By \eqref{tensor raise} $e_0(b' \otimes \bb_1) = e_0(b') \otimes \bb_1$
and $e_0(b'\otimes y')=e_0(b') \otimes y'$.
By \eqref{loc energy} and the fact that the local isomorphism
on $\Bc \otimes \Bc$ is the identity, we have
$E(e_0(b'\otimes \bb_1))=E(b'\otimes \bb_1)-1$.

To show that $E(e_0(b'\otimes y'))=E(b'\otimes y')-1$
we check the conditions of Lemma \ref{energy down}.
By \eqref{crystal string}
$\epsilon_0(y')=\phi_0(y')-\inner{h_0}{\wt(y')}=
0 - \inner{h_0}{\La'-\ell'\La_0} =
\ell'-\inner{h_0}{\La'}$.  Also by \eqref{tensor crystal},
since $\phi_0(y')=0$, we have
$\epsilon_0(b'\otimes y') = \epsilon_0(b')+\epsilon_0(y')
> \inner{h_0}{\La'}+\ell'-\inner{h_0}{\La'} = \ell'$.
Let $z \otimes x$ be the image of $b' \otimes y'$ under
an arbitrary composition of local isomorphisms.
Since $b' \otimes y'$ is an $\sln$ highest weight vector,
so is $z \otimes x$ and $x$.  Now $x$ is the $\sln$-highest
weight vector in a perfect crystal of level at most $\ell'$,
so $\phi_0(x)=0$ and $\epsilon_0(x)\le \ell'$.
But $\ell'<\epsilon_0(b'\otimes y')=\epsilon_0(z \otimes x)
=\epsilon_0(z)+\epsilon_0(x)$ so that $\epsilon_0(z)>0$.
By \eqref{tensor raise} $e_0(z \otimes x)= e_0(z) \otimes x$.
So $E(e_0(b'\otimes y'))=E(b'\otimes y')-1$ by Lemma \ref{energy down}.

By induction we may now assume that
$\epsilon_0(b') \le \inner{h_0}{\La'}$.
But then $\sum_i \epsilon_i(b') \le \sum_i \inner{h_i}{\La'}$,
or $\inner{c}{\epsilon(b')} \le \inner{c}{\La'} = \ell'$.
Since $b'\in \Bc$ and $\Bc$ is a perfect crystal of level $\ell'$,
$b'$ must be the unique element of $\Bc$ such that
$\epsilon(b')=\La'$.  Thus the function $\Bc\rightarrow\Z$
given by $b'\mapsto E(b'\otimes y')-E(b'\otimes \bb_1)$
is constant on $\Bc$ if it is constant on the
singleton set $\{\epsilon^{-1}(\La')\}$, which it obviously is.
\end{proof}

\subsection{Explicit ground state energy}

To go further, an explicit formula for the value $E(p \otimes y')$
is required.  This is achieved in \eqref{ground energy}.
The derivation makes use of the following explicit construction
of the local isomorphism.

\begin{theorem} \label{iso} Let $\Bc=\Bc^{k,\ell}$
be a perfect crystal of level $\ell$, $\La,\La'\in\Pclell$,
$\Bc'$ a perfect crystal of level $\ell'\le \ell$,
and $b\in \HH(\La',\Bc',\La)$.
Let $x\in\Bc$ (resp. $y\in\Bc$) be the unique element such that
$\epsilon(x)=\La$ (resp. $\epsilon(y)=\La'$).
Then under the local isomorphism $\Bc \otimes \Bc' \cong
\Bc'\otimes \Bc$, we have $x \otimes b \cong \psi^k(b) \otimes y$.
\end{theorem}

The proof requires several technical lemmas and is given in
the next section.

\begin{ex} Let $n=5$, $\ell=4$, $k=2$, $\La'=\La_0+\La_1+\La_3+\La_4$,
$\La=\La_0+\La_1+\La_2+\La_4$, $\ell'=2$, $\Bc'=\Bc^{2,2}$.
Here the set $\HH(\La',\Bc',\La)$ consists of two elements, namely,
\begin{equation*}
  \begin{matrix}
    1&2 \\
    4&5
  \end{matrix} \qquad \text{and} \qquad
  \begin{matrix}
    1&4 \\
    2&5 
  \end{matrix}.
\end{equation*}
Let $b$ be the second tableau. The theorem says that
\begin{equation*}
  \begin{matrix}
    1&1&2&3\\
    2&3&4&5
  \end{matrix}\, \otimes \,
  \begin{matrix}
    1&4 \\
    2&5 
  \end{matrix} \, \cong \,
  \begin{matrix}
    1&3 \\
    2&4 
  \end{matrix} \, \otimes \,
  \begin{matrix}
     1&1&2&4\\
     2&3&5&5 
  \end{matrix}.
\end{equation*}
\end{ex}

\begin{prop} Let $\La\in\Pclell$, $\Bc=\Bc^{k,\ell}$
a perfect crystal of level $\ell$, $\bb\in\PP(\La,\Bc)$
the ground state path, $p$ a finite path (say of length $N$
where $N$ is a multiple of $n$)
such that $p \otimes \bb = \bb$, $\Bc'$ the tensor product of
perfect crystals each of level at most $\ell$,
and $y\in \HH(\ell\La_0,\Bc',\La)$.  Let $p'$ be the path of length
$N$ such that $p' \otimes \bb' = \bb'$ where
$\bb'\in\PP(\ell\La_0,\Bc)$ is the ground state path.
Then under the composition of local isomorphisms
$\Bc^{\otimes N} \otimes \Bc' \cong \Bc' \otimes \Bc^{\otimes N}$
we have $p \otimes y \cong y \otimes p'$.
\end{prop}
\begin{proof} Induct on the length of the path $y$.
Suppose $\Bc'=\Bc_1 \otimes \Bc_2$ and
$y=y_1 \otimes y_2$ where $y_j\in \Bc_j$
and $\Bc_j$ is a perfect crystal.  Let $\La'=\La-\wt(y_1)$.
By the definitions $y_2\in\HH(\ell\La_0,\Bc_2,\La')$.
By induction the first $N$ steps $p''$ of the ground state path
of $\PP(\La',\Bc)$ satisfy
$p'' \otimes y_2 \cong y_2 \otimes p'$ under the composition of local
isomorphisms $\Bc^{\otimes N} \otimes \Bc_2\cong\Bc_2 \otimes \Bc^{\otimes N}$.
Tensoring on the left with $y_1$, it remains to show that
$p \otimes y_1 \cong y_1 \otimes p''$ under the composition of local
isomorphisms 
$\Bc^{\otimes N} \otimes \Bc_1 \cong \Bc_1\otimes \Bc^{\otimes N}$.
Now $p_N\in\Bc$ and $p''_N\in \Bc$ are the unique elements such that
$\epsilon(p_N)=\La$ and $\epsilon(p''_N)=\La'$.
Applying Theorem \ref{iso} we obtain
$p_N \otimes y_1 \cong \psi^k(y_1) \otimes p''_N$.
Now $p_N \otimes y_1 \in \HH(\La',\Bc \otimes \Bc_1,\phi(p_N))$
so that $\psi^k(y_1) \otimes p''_N\in
\HH(\La',\Bc_1\otimes \Bc,\phi(p_N))$.
This implies that $\psi^k(y_1)\in\HH(\phi(p''_N),\Bc_1,\phi(p_N))$.
Now by definition $\epsilon(p''_{N-1})=\phi(p''_N)$
and $\epsilon(p_{N-1})=\phi(p_N)$.
Applying Theorem \ref{iso} we obtain
$p_{N-1} \otimes \psi^k(y_1) \cong \psi^{2k}(y_1)\otimes p''_{N-1}$.
Continuing in this manner it follows that
$p_{N-j} \otimes \psi^{j k}(y_1) \cong \psi^{(j+1)k}(y_1) \otimes p''_{N-j}$
for $0\le j\le N-1$.  Composing these local isomorphisms it follows that
$p \otimes y_1 \cong \psi^{N k}(y_1) \otimes p''$.
But $\psi^N$ is the identity since the order of $\psi$ divides $n$
which divides $N$.  Therefore
$p \otimes y_1 \cong y_1 \otimes p''$ under the composition of local
isomorphisms and we are done.
\end{proof}

In the notation in the previous section,
$E(p \otimes y')=E(y' \otimes p')$ where
$p'$ is the first $N$ steps of the ground state path
of $\PP(\ell'\La_0,\Bc)$.  Write $N=n M$ and $\Bc=\Bc^{k,\ell'}$.
Then using the generalized cocyclage
one may calculate explicitly the generalized charge
of the LR tableau corresponding to the level $\ell'$ restricted
(and hence classically restricted) path $y' \otimes p'$.
Let $|y'|$ denote the total number of cells in the tableaux
comprising $y'$.  Then
\begin{equation} \label{ground energy}
  E(y' \otimes p') = E(y') + |y'| k M + n \ell'\binom{kM}{2}.
\end{equation}

\begin{ex} Let $n=5$, $\ell'=3$, $\La'=\La_0+\La_3+\La_4$,
$k=2$ and $M=1$.  Then $p'$ is the path
\begin{equation*}
\begin{matrix} 4&4&4 \\ 5&5&5 \end{matrix} \otimes
\begin{matrix} 2&2&2 \\ 3&3&3 \end{matrix} \otimes
\begin{matrix} 1&1&1 \\ 5&5&5 \end{matrix} \otimes
\begin{matrix} 3&3&3 \\ 4&4&4 \end{matrix} \otimes
\begin{matrix} 1&1&1 \\ 2&2&2 \end{matrix}.
\end{equation*}
The element $y'$ can be taken to be the tensor product
\begin{equation*}
  \begin{matrix} 1\\2\\3
  \end{matrix}
  \otimes
  \begin{matrix} 1\\2\\3\\4
  \end{matrix}.
\end{equation*}
Let $\la=(8,8,8,7,6)$.
Then the tableau $Q\in\LRT(\la;R)$ (resp. $Y$)
that records the path $y' \otimes p'$ (resp. $y'$) is given by
\begin{equation*}
 Q = \begin{matrix}
   1&1&1&5&5&5&11&15\\
   2&2&2&7&7&7&12&16\\
   3&3&3&8&8&8&13&17\\
   4&4&4&9&9&9&14&  \\
   6&6&6&10&10&10& & 
 \end{matrix}\, , \qquad
 Y = \begin{matrix}
   1&5\\2&6\\3&7 \\4&\\ &
 \end{matrix}
\end{equation*}
with $R=((3,3),(3,3),(3,3),(3,3),(3,3),(1,1,1,1),(1,1,1))$
and subalphabets $\{1,2\}$, $\{3,4\}$, $\{5,6\}$,
$\{7,8\}$, $\{9,10\}$, $\{11,12,13,14\}$, $\{15,16,17\}$.
The generalized charge $\charge_R(Q)$ is equal to the
energy $E(y' \otimes p')$ \cite[Theorem 23]{S3}.
Here the widest rectangle in the path is of width $\ell'$.
For any tableau $T\in\LRT(\rho;R)$ for some partition $\rho$,
define $V(T) = P((w^R_0 T_e)(w^R_0 T_w))$
where $P$ is the Schensted $P$ tableau, $w_0^R$ is the
automorphism of conjugation that reverses each of the
subalphabets, and $T_w$ and $T_e$ are the west and east subtableaux
obtained by slicing $T$ between the $\ell'$-th and $(\ell'+1)$-th
columns.  It can be shown that there is a composition of
$|T_e|$ generalized $R$-cocyclages leading from $T$ to $V(T)$
where $|T_e|$ denotes the number of cells in $T_e$.
It follows from the ideas in \cite[Section 3]{S} and the
intrinsic characterization of $\charge_R$ in \cite[Theorem 21]{S} that
\begin{equation} \label{crank stat}
\charge_R(T)=\charge_R(V(T))+|T_e|.
\end{equation}
For the above tableau $Q$ we have
\begin{equation*}
  Q_w = \begin{matrix}
    1&1&1\\
    2&2&2\\
    3&3&3\\
    4&4&4\\
    6&6&6
  \end{matrix}\qquad
  w_0^R Q_w = \begin{matrix}
    1&1&1\\
    2&2&2\\
    3&3&3\\
    4&4&4\\
    5&5&5
  \end{matrix}
\end{equation*}
and
\begin{equation*}
  Q_e = \begin{matrix}
   5&5&5&11&15\\
   7&7&7&12&16\\
   8&8&8&13&17\\
   9&9&9&14&  \\
   10&10&10& & 
  \end{matrix} \qquad
  w_0^R Q_e = \begin{matrix}
   6&6&6&11&15\\
   7&7&7&12&16\\
   8&8&8&13&17\\
   9&9&9&14&  \\
   10&10&10& & 
  \end{matrix}\, .
\end{equation*}
Then
\begin{equation*}
 V(Q) =
 \begin{matrix}
   1&1&1&11&15\\
   2&2&2&12&16\\
   3&3&3&13&17\\
   4&4&4&14&  \\
   5&5&5&  &  \\
   6&6&6&  &  \\
   7&7&7&  &  \\
   8&8&8&  &  \\
   9&9&9&  &  \\
   10&10&10&  &
 \end{matrix}  \quad \text{and} \quad
V(V(Q)) =
 \begin{matrix}
   1&1&1 \\
   2&2&2 \\
   3&3&3 \\
   4&4&4 \\
   5&5&5 \\
   6&6&6 \\
   7&7&7 \\
   8&8&8 \\
   9&9&9 \\
   10&10&10 \\
   11&15& \\
   12&16& \\
   13&17& \\
   14&  &
 \end{matrix}\, .
\end{equation*}
We have $c_R(V(V(Q)))=c_R(Y)=E(y')$ by \cite[Theorem 21]{S} and
$c_R(Q)=c_R(V(Q))+|Q_e| = c_R(V(Q))+\ell'n+|Y|$,
and $c_R(V(Q))=c_R(V(V(Q)))+|Y|$ by \eqref{crank stat}.
This implies $c_R(Q)=\ell'n+ E(y')+2 |Y|$.
\end{ex}

\subsection{Proof of Theorem~\ref{iso}}

The proof of Theorem~\ref{iso} requires several lemmas.

Words of length $L$ in the alphabet $\{1,2,\dots,n\}$ are
identified with the elements of the crystal basis of the
$L$-fold tensor product $(\Bc^{1,1})^{\otimes L}$.

\begin{lemma} \label{hw} Let $u$ and $v$ be words such that
$uv$ is an $A_{n-1}$ highest weight vector.  Then
$v$ is an $A_{n-1}$ highest weight vector and
$\epsilon_j(u)\le \phi_j(v)$ for all $1\le j\le n-1$.
\end{lemma}
\begin{proof} Let $uv$ be an $A_{n-1}$ highest weight vector
and $1\le j \le n-1$.  By \eqref{tensor crystal}
\begin{equation*}
  0=\epsilon_j(uv)=\epsilon_j(v) + \max\{0,\epsilon_j(u)-\phi_j(v)\}.
\end{equation*}
Since both summands on the right hand side are nonnegative
and sum to zero they must both be zero.
\end{proof}

\begin{lemma} \label{letter} Let $w$ be a word in the alphabet
$\{1,2\}$ and $\hatw$ a word obtained by removing a letter
$i$ of $w$.  Then
\begin{enumerate}
\item \label{remove} $\epsilon_1(\hatw) \le \epsilon_1(w) + 1$
with equality only if $i=1$.
\item \label{insert} $\epsilon_1(w) \le \epsilon_1(\hatw)+1$ with equality
only if $i=2$.
\end{enumerate}
\end{lemma}
\begin{proof} Write $w=uiv$ and $\hatw=uv$.  By \eqref{tensor crystal}
\begin{equation} \label{ui}
\begin{split}
  \epsilon_1(u i) &= \epsilon_1(i) + \max\{0,\epsilon_1(u)-\phi_1(i)\} \\
  &= \begin{cases}
    \max\{0,\epsilon_1(u)-1\} & \text{if $i=1$} \\
    1+\epsilon_1(u) & \text{if $i=2$.}
  \end{cases}
\end{split}  
\end{equation}
In particular $\epsilon_1(ui)\ge \epsilon_1(u)-1$.
Applying \eqref{tensor crystal} to both $\epsilon_1(uv)$ and
$\epsilon_1(uiv)$ and subtracting, we obtain
\begin{equation*}
\begin{split}
  \epsilon_1(uv)-\epsilon_1(uiv) &=
  \max\{0,\epsilon_1(u)-\phi_1(v)\} -
  \max\{0,\epsilon_1(ui)-\phi_1(v)\} \\
  &\le \max\{0,\epsilon_1(u)-\phi_1(v)\} - 
  \max\{0,\epsilon_1(u)-1-\phi_1(v)\} \\
  &\le 1.
\end{split}
\end{equation*}
Moreover if $\epsilon_1(uv)-\epsilon_1(uiv)=1$ then
all of the inequalities are equalities.  In particular
it must be the case that $\epsilon_1(ui)=\epsilon_1(u)-1$,
which by \eqref{ui} implies that $i=1$, proving the first
assertion.

On the other hand, \eqref{ui} also implies
$\epsilon_1(ui) \le 1+\epsilon_1(u)$.
Subtracting $\epsilon_1(uv)$ from $\epsilon_1(uiv)$ and
computing as before, the second part follows.
\end{proof}

Say that $w$ is an almost highest weight vector with defect $i$
if there is an index $1\le i\le n-1$
such that $\epsilon_j(w)=\delta_{ij}$ for $1\le j\le n-1$,
and also $\epsilon_{i-1}(e_i(w))=0$ if $i>1$.

\begin{lemma} \label{almost hw} Let $w$ be an almost highest weight vector
with defect $i$ for $1\le i\le n-1$.
Then $e_i(w)$ is either an $A_{n-1}$ highest weight vector or an
almost highest weight vector of defect $i+1$.
\end{lemma}
\begin{proof} For $j\not\in\{i-1,i,i+1\}$, the restriction of the words
$w$ and $e_i(w)$ to the alphabet $\{j,j+1\}$ are identical, so that
$\epsilon_j(e_i(w))=\epsilon_j(w)=0$ by the definition of
an almost highest weight vector.
Also $\epsilon_i(w)=1$ implies that $\epsilon_i(e_i(w))=0$.
Again by the definition of an almost highest weight vector,
$\epsilon_{i-1}(e_i(w))=0$.

If $i=n-1$ we have shown that $e_i(w)$ is an $A_{n-1}$ highest weight vector.
So it may be assumed that $i<n-1$.
It is enough to show that one of the two following
possibilities occurs.
\begin{enumerate}
\item $\epsilon_{i+1}(e_i(w))=0$.
\item $\epsilon_{i+1}(e_i(w))=1$ and
$\epsilon_i(e_{i+1}e_i(w))=0$.
\end{enumerate}
Recall that $e_i(w)$ is obtained from $w$
by changing an $i+1$ into an $i$.  Write $w=u (i+1) v$
such that $e_i(w)=u i v$.  In this notation we have
$\phi_i(v)=0$ and $\epsilon_i(u)=0$.  
By Lemma \ref{letter} point \ref{remove} with $\{1,2\}$
replaced by $\{i+1,i+2\}$ and using that 
$w$ is an almost highest weight vector of defect $i$, we have
$\epsilon_{i+1}(e_i(w)) \le \epsilon_{i+1}(w)+1=1$.
It is now enough to assume that $\epsilon_{i+1}(e_i(w))=1$
and to show that $\epsilon_i(e_{i+1}e_i(w))=0$.
By \eqref{tensor crystal}
\begin{equation*}
\begin{split}
  0 &=\epsilon_{i+1}(w)=\epsilon_{i+1}(u (i+1) v) \\
  &= \epsilon_{i+1}(v) + \max\{0,\epsilon_{i+1}(u)-\phi_{i+1}((i+1)v)\}.
\end{split}
\end{equation*}
In particular $\epsilon_{i+1}(v)=0$.
Hence $e_{i+1}(e_i(w))=e_{i+1}(u i v)=e_{i+1}(u)iv$.
Similar computations starting with $\epsilon_i(w)=1$
and which use the fact that $\epsilon_i(u)=\phi_i(v)=0$,
yield $\epsilon_i(v)=0$.  We have
\begin{equation*}
\begin{split}
  \epsilon_i(e_{i+1}e_i(w)) &= \epsilon_i(e_{i+1}(u) i v ) \\
  &= \epsilon_i(iv)+\max\{0,\epsilon_i(e_{i+1}(u))-\phi_i(iv) \} \\
  &= 0+\max\{0,\epsilon_i(e_{i+1}(u))-1 \}.
\end{split}
\end{equation*}
But $\epsilon_i(u)=0$ and in passing from $u$ to $e_{i+1}(u)$
an $i+2$ is changed into an $i+1$.  By Lemma \ref{letter} point \ref{insert}
applied to the restriction of $u$ to the alphabet $\{i,i+1\}$,
we have $\epsilon_i(e_{i+1}(u)) \le \epsilon_i(u)+1=1$.
It follows that $\epsilon_i(e_{i+1}e_i(w))=0$,
and that $e_i(w)$ is an almost highest weight vector of defect $i+1$.
\end{proof}

\begin{lemma} \label{lose letter}
Suppose $w$ is an $A_{n-1}$ highest weight vector
and $\hatw$ is a word obtained by removing
a letter (say $i$) from $w$.  Then there is an index $r$
such that $i\le r\le n$ and $e_{r-1}e_{r-2} \dotsm e_i(\hatw)$
is an $A_{n-1}$ highest weight vector.
\end{lemma}
\begin{proof} By Lemma \ref{almost hw} it suffices to show that
$\hatw$ is either an $A_{n-1}$ highest weight vector or an almost
highest weight vector of defect $i$.

First it is shown that $\epsilon_j(\hatw)=0$ for $j\not=i$.
For $j\not\in\{i-1,i\}$, the restrictions of $w$ and $\hatw$
to the alphabet $\{j,j+1\}$ are the same, so that
$\epsilon_j(\hatw)=\epsilon_j(w)=0$.
For $j=i-1$, by Lemma \ref{letter} point \ref{remove} and the
assumption that $w$ is an $A_{n-1}$ highest weight vector, it follows that
$\epsilon_{i-1}(\hatw) \le \epsilon_{i-1}(w)+1=1$.
But equality cannot hold since the removed letter is $i$
as opposed to $i-1$.  Thus $\epsilon_{i-1}(\hatw)=0$.

Next we observe that $\epsilon_i(\hatw) \le \epsilon_i(w) + 1 = 1$
by Lemma \ref{letter} point \ref{remove} and the fact that
$w$ is an $A_{n-1}$ highest weight vector.

If $\epsilon_i(\hatw)=0$ then $\hatw$ is
an $A_{n-1}$ highest weight vector.
So it may be assumed that $\epsilon_i(\hatw)=1$.
It suffices to show that $\epsilon_{i-1}(e_i(\hatw))=0$.
Write $w=uiv$ and $\hatw=uv$.  Now $\epsilon_j(v)=0$ for all
$1\le j\le n-1$ by Lemma \ref{hw}
since $w$ is an $A_{n-1}$ highest weight vector.
In particular $\epsilon_i(v)=0$ so that
$e_i(\hatw)=e_i(uv)=e_i(u) v$.
We have
\begin{equation*}
\begin{split}
  \epsilon_{i-1}(e_i(\hatw)) &= \epsilon_{i-1}(e_i(u)v) \\
  &= \epsilon_{i-1}(v) + \max\{0,\epsilon_{i-1}(e_i(u))-\phi_{i-1}(v)\} \\
  &= \max\{0,\epsilon_{i-1}(e_i(u))-\phi_{i-1}(v)\}
\end{split}
\end{equation*}
since $\epsilon_{i-1}(v)=0$ by Lemma \ref{hw}.  It is enough to show that
$\epsilon_{i-1}(e_i(u)) \le \phi_{i-1}(v)$.  But
\begin{equation*}
  \epsilon_{i-1}(e_i(u)) \le \epsilon_{i-1}(u) + 1
  = \epsilon_{i-1}(ui) \le \phi_{i-1}(v).
\end{equation*}
The first inequality holds by an application of Lemma \ref{letter}
point \ref{insert} since the restrictions of $u$ and $e_i(u)$ to the alphabet
$\{i-1,i\}$ differ by inserting a letter $i$.  The last inequality
holds by Lemma \ref{hw} since $w=uiv$ is an $A_{n-1}$ highest weight vector.
\end{proof}

\begin{lemma} \label{hole} Let $\Bc=\Bc^{k,\ell'}$
be a perfect crystal of level $\ell'\le \ell$, $\La\in\Pclell$,
$\Bc'$ a finite (possibly empty) tensor product of perfect crystals of level
at most $\ell$, $x\in \Bc'$ and $b \in \Bc$ such that
$x\otimes b\in \HH(\La,\Bc'\otimes \Bc)$. Let $i\in J$ such that
$\inner{h_i}{\La}>0$ and set $\La'=\La-\La_i+\La_{i-1}$. Then there is an index
$0\le s\le k$ such that
\begin{equation} \label{raise right}
e_{i+s-1} \dotsm e_{i+1} e_i (x \otimes b)
=x \otimes e_{i+s-1} \dotsm e_{i+1} e_i (b)
\end{equation}
and $e_{i+s-1} \dotsm e_i(b) \in \HH(\La',\Bc)$
where the subscripts are taken modulo $n$.
Moreover if $\ell'=\ell$ then $s=k$.
\end{lemma}
\begin{proof} Since the Dynkin diagram $A_{n-1}^{(1)}$ has an automorphism
given by rotation, it may be assumed that $i=1$.
Let $\la$ be the partition of length less than $n$,
given by $\inner{h_j}{\La} = \la_j-\la_{j+1}$ for $1\le j\le n-1$
and $\la_n=0$.  Since $\inner{h_1}{\La}>0$ it follows that
$\la$ has a column of size $1$.  Let $m=\la_1$ and $y_i$ be the
$A_{n-1}$-highest weight vector in $\Bc^{\la^t_j,1}$ for $1\le j\le m$.
Write $y=y_m \otimes \dots \otimes y_1$ and
$\haty=y_{m-1} \otimes \dots \otimes y_1$.
Observe that $y \otimes u_{\ell\La_0}$ is 
an affine highest weight vector in
$\Bc^{\la^t_m,1} \otimes \dots \otimes \Bc^{\la^t_1,1} \otimes
\BB(\ell\La_0)$ and has weight $\La$ so its connected
component is isomorphic to $\BB(\La)$.
A similar statement holds for $\haty \otimes u_{\ell\La_0}$ and $\BB(\La')$.
In particular, $b\otimes y$ is an $A_{n-1}$ highest weight vector.
The map $x \otimes b \otimes y \mapsto \word(x)\word(b)\word(y)$
gives an embedding of $A_{n-1}$-crystals into a tensor product
of crystals $\Bc^{1,1}$.
By Lemma \ref{lose letter}, there exists an index $1\le r\le n$ such that
$e_{r-1} e_{r-2} \dotsm e_1 (\word(x) \word(b) \word(\haty))$ is an
$A_{n-1}$ highest
weight vector. Since $\haty$ is an $A_{n-1}$ highest weight vector
it follows that $e_{r-1} \dotsm e_1(\word(x) \word(b) \word(\haty)) =
e_{r-1} \dotsm e_1(\word(x)\word(b)) \word(\haty)$.

Let $p_j$ be the position of the letter in
$e_{j-1}\dots e_1(\word(x)\word(b))$
that changes from a $j+1$ to $j$ upon the application of $e_j$,
for $1\le j \le r-1$.
It follows from the proof of Lemma \ref{almost hw} that
\begin{equation} \label{decreasing position}
  p_{r-1} < p_{r-2} < \dotsm < p_2 < p_1.
\end{equation}
Let $s$ be the maximal index such that $p_s$ is located in
$\word(b)$.  Write $b'=e_s \dotsm e_1(b)$.
It follows that $e_s e_{s-1} \dotsm e_1(x \otimes b)=
x \otimes b'$ and that $b' \otimes \haty$ is an
$A_{n-1}$ highest weight vector.

It remains to show that
\begin{equation} \label{no zero}
  \epsilon_0(b' \otimes \haty \otimes u_{\ell\La_0}) = 0
\end{equation}
and that $s \le k$ with equality if $\ell'=\ell$.

Consider the corresponding positions in the tableau $b$.
Since $b\mapsto\word(b)$ is an $A_{n-1}$-crystal morphism, 
$e_s \dotsm e_1(\word(b))=\word(e_s \dotsm e_1(b))$.
Let $(i_1,j_1)$ be the position in the tableau $b$
corresponding to the position $p_1$ in $\word(b)$,
and analogously define $(i_2,j_2)$, $(i_3,j_3)$, and so on.
Since the rows of all tableaux (and in particular
$b$, $e_1(b)$, $e_2 e_1(b)$, etc.) are weakly increasing
and \eqref{decreasing position} holds, it follows that
$i_1<i_2<i_3<\dots<i_s$.  But $b$ has $k$ rows, so $s\le k$.

The next goal is to prove \eqref{no zero}.
Suppose first that $s<n-1$.  In this case the letters $1$ and $n$
are undisturbed in passing from
$e_1(b)$ to $e_s \dotsm e_1(b)$.  Using this and the
Dynkin diagram rotation it follows that
\begin{equation} \label{zero root}
\begin{split}
  \epsilon_0(e_s \dotsm e_2 e_1(b) \otimes \haty \otimes
  u_{\ell\La_0}) &=
  \epsilon_0(e_1(b) \otimes u_{\La'}) \\
  &=  \max\{0,\epsilon_0(e_1(b))-\phi_0(u_{\La'})\} \\
  &=  \max\{0,\epsilon_0(e_1(b))-\phi_0(u_{\La})-1\}.
\end{split}
\end{equation}
But $\phi_0(u_\La) \ge \epsilon_0(b) \ge \epsilon_0(e_1(b))-1$
by the fact that $\epsilon_0(b \otimes u_\La)=0$
and Lemma \ref{letter} point \ref{insert} applied after rotation of the
Dynkin diagram.  By \eqref{zero root} the desired result
\eqref{no zero} follows.

Otherwise assume $s=n-1$.  Here $k=n-1$ since $s\le k<n$ with the
inequality holding by the perfectness of $\Bc$. By
\eqref{decreasing position} and the fact that $b$ is a tableau, it must
be the case that $e_1$ acting on $b$
changes a $2$ in the first row of $b$ into a $1$, $e_2$ acting on
$e_1(b)$ changes a $3$ in the second row of $e_1(b)$ into a $2$, etc.
Since $b$ is a tableau with $n-1$ rows with entries between
$1$ and $n$, there are integers $0\le \nu_{n-1} \le \nu_{n-2} \le
\dots \le \nu_1 < \ell'$ such that the $i$-th row of $b$ consists
of $\nu_i$ copies of the letter $i$ and $\ell'-\nu_i$ copies of the
letter $i+1$.  For tableaux $b$ of this very special form, the explicit
formula for $e_0$ in \cite[(3.11)]{S3} yields
$\epsilon_0(b) = \ell' - m_n(b)$ where $m_n(b)$ is the number of occurrences
of the letter $n$ in $b$.  Since $b'=e_{n-1}\dotsm e_1(b)$ also has the
same form (with $\nu_i$ replaced by $\nu_i+1$ for $1\le i\le n-1$)
and $m_n(b')=m_n(b)-1$, it follows that $\epsilon_0(b')=\epsilon_0(b)+1$.
We have
\begin{equation*}
\begin{split}
  \epsilon_0(b' \otimes \haty \otimes u_{\ell\La_0}) &=
  \epsilon_0(b' \otimes u_{\La'}) \\
  &= \max \{ 0, \epsilon_0(b') - \phi_0(u_{\La'}) \} \\
  &= \max \{ 0, \epsilon_0(b)+1 - (\phi_0(u_\La)+1) \} = 0
\end{split}
\end{equation*}
since $b\in\HH(\La,\Bc)$.

Finally, assuming $\ell'=\ell$, it must be shown that $s=k$.
Since the level of $\Bc$ is the same as that of the weights
$\La$ and $\La'$, it follows from the perfectness of $\Bc$ that
both $b$ and $b'$ are uniquely defined by the
property that
$\epsilon(b)=\La$ and $\epsilon(b')=\La'$.
Let $\La=\sum_{i=0}^{n-1} z_i \La_i$.
By the explicit construction of $b$ in Example \ref{min bijection}
\begin{equation*}
  \wt(b) = \sum_{j=1}^k \sum_{i=0}^{n-1} z_i (\La_{i+j}-\La_{i+j-1})
         =\sum_{i=0}^{n-1} z_i (\La_{i+k}-\La_i)
\end{equation*}
with indices taken modulo $n$.
Subtracting the analogous formula for $\wt(b')$,
$\wt(b)-\wt(b') = -\sum_{j=1}^k \alpha_j$.
Using \eqref{crystal string} it follows that $k=s$.
\end{proof}

\begin{proof}[Proof of Theorem \ref{iso}]
First observe that $x\otimes b\in\HH(\La',\Bc\otimes \Bc',\phi(x))$
by \eqref{crystal string}, $b\in \HH(\La',\Bc',\La)$,
and $\epsilon(x)=\La$.  Let $c\in\Bc'$ and $z\in \Bc$ be such that
$x\otimes b \cong c \otimes z$ under the local isomorphism.
Then $c\otimes z\in\HH(\La',\Bc'\otimes\Bc,\phi(x))$
which means that $z$ is $\La'$-restricted. Hence
$z\in \HH(\La',\Bc,\phi(z))$ and $c\in\HH(\phi(z),\Bc',\phi(x))$.
The former together with the perfectness of $\Bc$ implies that $y=z$.
{}From the latter it follows that $\psi^{-k}(c) \in \HH(\La',\Bc',\La)$.
However the set $\HH(\La',\Bc',\La)$ might have multiplicities
so it is not obvious why $b=\psi^{-k}(c)$ or equivalently $c=\psi^k(b)$.

The proof proceeds by an induction that changes
the weight $\La'$ to a weight $\hLap$ that is ``closer to" $\ell\La_0$.
Suppose first that there is a root direction $i\not=0$
such that $\inner{h_i}{\La'}>0$ and 
$\hLap=\La'-\La_i+\La_{i-1}$.  By Lemma \ref{hole}
applied for the weight $\La'$, simple root $\alpha_i$,
and element $x\otimes b\in \HH(\La',\Bc \otimes \Bc')$,
there is an $0\le s<n$ such that
$\hatb=e_{i+s-1} \dotsm e_{i+1} e_i(b) \in \HH(\hLap,\Bc',\hLa)$
where $\hLa=\La-\La_{s+i}+\La_{s+i-1}$ and
$e_{i+s-1} \dotsm e_i(x\otimes b)=x \otimes \hatb$.
Applying Lemma \ref{hole} with $\La$, $\alpha_{s+i}$, and
$x\in\HH(\La,\Bc)$, it follows that
$\hatx=e_{k+s+i-1}\dotsm e_{s+i}(x) \in \HH(\hLa,\Bc)$.

The above computations imply
$e_{k+s+i-1} \dotsm e_i (x \otimes b) = \hatx \otimes \hatb
\in \HH(\hLap,\Bc \otimes \Bc')$.

We have $e_{k+s+i-1} \dotsm e_{i+1} e_i(c \otimes y)\in \HH(\hLap,\Bc'\otimes \Bc)$
since $x\otimes b\mapsto c\otimes y$ under the local isomorphism.
It must be seen which of these raising operators act on the tensor factor
in $\Bc'$ and which act in $\Bc$.
By Lemma \ref{hole} applied with $\La'$, $\alpha_i$,
and $c \otimes y\in\HH(\La',\Bc' \otimes \Bc)$,
it follows that $\haty=e_{k+i-1} \dotsm e_i(y)\in\HH(\hLap,\Bc)$
and that $e_{k+i-1} \dotsm e_i(c \otimes y)=c \otimes \haty$.
Since $\haty\otimes u_{\hLap}$ is an $A_{n-1}^{(1)}$ highest weight vector,
the rest of the raising operators $e_{s+k-1} \dotsm e_{k+i}$ must
act on the first tensor factor.  
Let $\hatc=e_{k+s+i-1} \dotsm e_{k+i}(c)$.
Then $e_{k+s+i-1} \dotsm e_i (c \otimes y) = \hatc \otimes \haty$.
But the local isomorphism is a crystal morphism so
it sends $\hatx \otimes \hatb \mapsto \hatc \otimes \haty$.
By induction $\hatc = \psi^k(\hatb)$.
By \eqref{psi and raising} it follows that $c=\psi^k(b)$.

Otherwise there is no index $i\not=0$ such that $\inner{h_i}{\La'}>0$.
This means $\La'=\ell\La_0$.  
But the sets $\HH(\ell\La_0,\Bc,\La)$ and
$\HH(\ell\La_0,\Bc',\phi(y))$ are singletons
whose lone elements are given by the $A_{n-1}$ highest weight vectors
in $\Bc$ and $\Bc'$ respectively.
Since $\Bc \otimes \Bc'$ is $A_{n-1}$ multiplicity-free
it follows that the sets $\HH(\phi(y),\Bc',\phi(x))$ and
$\HH(\La,\Bc,\phi(x))$ are singletons.
In this case it follows directly that $c=\psi^k(b)$
since both $c$ and $\psi^k(b)$ are elements of the
set $\HH(\phi(y),\Bc',\phi(x))$.
\end{proof}

\subsection{Branching function by restricted generalized
Kostka polynomials}
\label{sec bf K}

The appropriate map from LR tableaux to rigged configurations,
sends the generalized charge of the LR tableau to the
charge of the rigged configuration.  Unfortunately in general
it is not clear what happens when one uses the statistic
coming from the energy function $E(b\otimes y')$ but using
the path $b \otimes y' \otimes y''$.  It is only known that
the statistic $E(b\otimes y' \otimes y'')$ on the path
$b\otimes y' \otimes y''$, is well-behaved.  So to continue the
computation we require that $y''=\emptyset$.  This is
achieved when $\La''=\ell''\La_0$.  So let us assume this.

The other problem is that we do not consider all
paths in $\HH(\ell\La_0,\Bc^{\otimes N} \otimes \Bc',\La)$,
but only those of the form $b \otimes y'$ where
$y'\in\Bc'$ is a fixed path.  Passing to LR tableaux,
this is equivalent to imposing an additional condition that the
subtableaux corresponding to the first several rectangles,
must be in fixed positions.
Conjecture~\ref{conj skew} asserts that the corresponding sets 
of rigged configurations are well-behaved.

The special case that requires no extra work, is when
$\Bc'$ consists of a single perfect crystal.  This is achievable
when $\La'$ has the form $\La'=r\La_s + (\ell'-r)\La_0$;
in this case $\Bc'=\Bc^{s,r}$ and $y'$ is the $\sln$-highest
weight element of $\Bc^{s,r}$.  This is the same as requiring that
the first subtableau of the LR tableau be fixed.  But this is always
the case.  Let $R^{(M)}$ consist of the single rectangle $(r^s)$ followed by
$N=M n$ copies of the rectangle $({\ell'}^k)$ where $\Bc=\Bc^{k,\ell'}$.
Let $\la^{(M)}$ be the partition of the same size as the total size of
$R^{(M)}$, such that $\la^{(M)}$ projects to $\La-\ell\La_0$.
Then the set of paths $\HH(\ell\La_0,\Bc^{\otimes N}\otimes \Bc^{s,r},\La)$
is equal to $\PU^\ell_{\La-\ell\La_0,R^{(M)}}$.  This is summarized by
\begin{equation} \label{Kostka branch}
  b^\La_{\La'\La''}(q) =
  \lim_{M\rightarrow\infty}
  q^{-r s k M - n \ell' \binom{kM}{2}} K^\ell_{\la^{(M)},R^{(M)}}(q),
\end{equation}
where $\La$ is arbitrary, $\La'=r\La_s + (\ell'-r)\La_0$,
and $\La''=\ell''\La_0$.

Inserting expression \eqref{Kmn} for the generalized Kostka polynomial
in \eqref{Kostka branch} and taking the limit yields the following fermionic 
expression for the branching function
\begin{multline}\label{bf ferm}
b_{\La'\La''}^\La(q)=q^{\frac{rs(s-n)}{2n}
 +\frac{1}{2\ell}\sum_{j=1}^n (\la_j-\frac{|\la|}{n})^2}
 \sum_{S\in\SCST(\la')} (-1)^{|S|+1} 
 q^{\frac{1}{2}\vu(S)C^{-1}\otimes C^{-1}\vu(S)}\\
\times
\sum_{\vm} q^{\frac{1}{2}\vm C\otimes C^{-1}\vm-\vm I\otimes C^{-1} \vu(S)}
\Bigl(\prod_{\substack{i=1\\ i\neq \ell'}}^{\ell-1} \prod_{a=1}^{n-1} 
\qbins{m_i^{(a)}+n_i^{(a)}}{m_i^{(a)}}\Bigr)
\Bigl(\prod_{a=1}^{n-1}\frac{1}{(q)_{m_{\ell'}^{(a)}}}\Bigr),
\end{multline}
where $\la$ is any partition which projects to $\La-\ell\La_0$ and
$\vu(S)$ as defined in \eqref{u}. The sum over $\vm$ runs over all
$\vm=\sum_{i=1}^{\ell-1}\sum_{a=1}^{n-1}m_i^{(a)} \ve_i\otimes \ve_a$ 
such that $m_i^{(a)}\in\Z$ and
\begin{equation*}
\ve_{\ell-1}\otimes \ve_a (I\otimes C^{-1}\vm-C^{-1}\otimes C^{-1}\vu(S))
 -\frac{1}{\ell}\sum_{j=1}^a(\la_j-\frac{1}{n}|\la|)\in\Z
\end{equation*}
for all $1\le a\le n-1$. The variables $n_i^{(a)}$ are given by
\begin{equation*}
n_i^{(a)}=\ve_i\otimes \ve_a \bigl\{ -C\otimes C^{-1} \vm+I\otimes C^{-1}
(\vu(s)+\ve_r\otimes \ve_s)\bigr\}
\end{equation*}
for all $1\le a<n$ and $1\le i<\ell$, $i\neq \ell'$.

\section{Proof of Theorem \ref{main}}
\label{sec proof}

To prove Theorem \ref{main} it clearly suffices to
show that there is a bijection
$\psib_R:\RLR^\ell(\la;R)\rightarrow\RC^\ell(\la;R)$
that is charge-preserving, that is,
$\charge_R(T)=\charge(\psib_R(T))$ for all $T\in\RLR^\ell(\la;R)$.
Here we identify $\LRT(\la;R)$
with $\RLR(\la;R)$ via the standardization bijection $\std$.
Also define $\charge'_R:\CLR(\la;R)\rightarrow\N$ by
$\charge'_R =\charge_R \circ \gamma_R$ where
$\charge_R:\RLR(\la;R)\rightarrow\N$.
It will be shown that one of the standard bijections
$\psib_R:\RLR(\la;R)\rightarrow\RC(\la;R)$ is
charge-preserving, and that it restricts to a bijection
$\RLR^\ell(\la;R)\rightarrow\RC^\ell(\la;R)$.

With this in mind let us review the
bijections from LR tableaux to rigged configurations.

\subsection{Bijections from LR tableaux to rigged configurations}

A bijection $\phib_R:\CLR(\la;R)\rightarrow\RC(\la^t;R^t)$ was defined
recursively in \cite[Definition-Proposition 4.1]{KSS}.  It is one of four
natural bijections from LR tableaux to rigged configurations.
\begin{enumerate}
\item Column index quantum: $\phib_R:\CLR(\la;R) \rightarrow \RC(\la^t;R^t)$.
\item Column index coquantum:
$\phit_R:\CLR(\la;R) \rightarrow \RC(\la^t;R^t)$, defined by
$\phit_R=\theta_{R^t} \circ \phib_R$.
\item Row index quantum: $\psib_R:\RLR(\la;R) \rightarrow \RC(\la;R)$,
defined by $\psib_R= \phib_{R^t} \circ\tr$, and
\item Row index coquantum: $\psit_R:\RLR(\la;R) \rightarrow \RC(\la;R)$,
defined by $\psit_R=\theta_R \circ \psib_R$.
\end{enumerate}

Of these four, the one that is compatible with level-restriction is
$\psib$.  First we show that it is charge-preserving.
This fact is a corollary of the difficult result \cite[Theorem 9.1]{KSS}.

\begin{prop} \label{c and psib} $\charge(\psib_R(T))=\charge_R(T)$
for all $T\in\RLR(\la;R)$.
\end{prop}
\begin{proof} Consider the following diagram, which commutes by
the definitions and \cite[Theorem 7.1]{KSS}
\begin{equation*}
\xymatrix{
 & & {\RLR(\la;R)} \ar[dll]_{\gamma_R^{-1}}
	\ar[dl]^{\tr} \ar[ddl]^{\psib_R} \\
 {\CLR(\la;R)} \ar[r]_{\LRtr} \ar[d]_{\phib_R} &
 {\CLR(\la^t;R^t)} \ar[d]^{\phib_{R^t}} & \\
 {\RC(\la^t;R^t)} \ar[r]_{\RCtr} & {\RC(\la;R).} & 
}
\end{equation*}
In particular $\psib_R = \RCtr\circ \phib_R \circ \gamma_R^{-1}$.
Let $T\in\RLR(\la;R)$ and $Q=\gamma_R^{-1}(T)$.  Then, using
$\RCtr\circ\theta_{R^t}=\theta_R\circ\RCtr$,
\begin{equation*}
  \psib_R(T) = \theta_R(\RCtr(\phit_R(Q))).
\end{equation*}
Let $(\nu,J)=\RCtr(\phit_R(Q))$.  Then
\begin{equation*}
\begin{split}
  \charge(\psib_R(T)) &= \charge(\theta_R(\nu,J)) =
  ||R||-\cc(\nu,J) \\
  &= ||R||-\cc(\RCtr(\phit_R(Q))) =
  \cc(\phit_R(Q)) =
  \charge'_R(Q) = \charge_R(T)
\end{split}
\end{equation*}
by Lemma \ref{theta and c}, \eqref{RCtr and cc} and
\cite[Theorem 9.1]{KSS} to pass from $\cc$ to $\charge'_R$.
\end{proof}

In light of Proposition \ref{c and psib}, 
to prove Theorem \ref{main} it suffices to establish the following result.

\begin{theorem} \label{rect res} The bijection
$\psib_R:\RLR(\la;R) \rightarrow \RC(\la;R)$
restricts to a well-defined bijection
$\psib_R:\RLR^\ell(\la;R)\rightarrow\RC^\ell(\la;R)$.
\end{theorem}

Computer data suggests that the bijection $\psib_R$ is not
only well-behaved with respect to level-restriction, but also
with respect to fixing certain subtableaux. It was argued in
Section~\ref{sec bf K} that the branching functions can
be expressed in terms of generating functions of tableaux
with certain fixed subtableaux.

Let $\rho\subset\la$ be partitions,
$R_\rho=((1^{\rho_1^t}),\ldots,(1^{\rho_n^t}))$
and $T_\rho$ the unique tableau in $\RLR(\rho;R_\rho)$.
Define $\RLR^\ell(\la,\rho;R)$ to be the set of tableaux
$T\in\RLR^\ell(\la;R_\rho\cup R)$ such that $T$ restricted
to shape $\rho$ equals $T_\rho$. Recall the set of rigged
configurations $\RC^\ell(\la,\rho;R)$ defined in 
Section~\ref{sec lev rc}.
\begin{conj}\label{conj skew}
The bijection $\psib_R:\RLR(\la;R) \rightarrow \RC(\la;R)$
restricts to a well-defined bijection
$\psib_R:\RLR^\ell(\la,\rho;R)\rightarrow\RC^\ell(\la,\rho;R)$.
\end{conj}

\subsection{Reduction to single rows}

In this section it is shown that to prove Theorem \ref{rect res}
it suffices to consider the case where $R$ consists of single rows.

Recall the nontrivial embedding
$i_R:\LRT(\la;R)\hookrightarrow\LRT(\la;\rows(R))$.
We identify $\LRT(\la;R)$ and $\RLR(\la;R)$ via $\std$,
and therefore have an embedding
$i_R:\RLR(\la;R)\hookrightarrow\RLR(\la;\rows(R))$.

Define a map $j_R:\RC(\la;R)\rightarrow\RC(\la;\rows(R))$
as follows.  Let $(\nu,J)\in\RC(\la;R)$.  For each rectangle of
$R$ having $k$ rows and $m$ columns, add $k-j$ strings
$(m,0)$ of length $m$ and label zero to the rigged partition
$(\nu,J)^{(j)}$ for $1\le j\le k-1$.  The resulting rigged configuration
is $j_R(\nu,J)$.

\begin{prop} \label{psib embed}
The following diagram commutes:
\begin{equation*}
\begin{CD}
  \RLR(\la;R) @>{i_R}>> \RLR(\la;\rows(R)) \\
  @V{\psib_R}VV	@VV{\psib_{\rows(R)}}V \\
  \RC(\la;R) @>>{j_R}> \RC(\la;\rows(R)).
\end{CD}
\end{equation*}
\end{prop}

It must be shown that similar diagrams commute in which
$i_R$ is replaced by either $\il_R$ or $s_p$, the
maps that occur in the definition of $i_R$.

Let $\jl_R:\RC(\la;R) \rightarrow \RC(\la;\Rl)$ be defined by
adding a string $(\mu_1,0)$ to each of the first $\eta_1-1$
rigged partitions in $(\nu,J)\in\RC(\la;R)$.

\begin{lemma} \label{psib less} $\jl_R$ is well-defined and
the following diagram commutes:
\begin{equation*}
\begin{CD}
  \RLR(\la;R) @>{\il_R}>> \RLR(\la;\Rl) \\
  @V{\psib_R}VV	@VV{\psib_{\Rl}}V \\
  \RC(\la;R) @>>{\jl_R}> \RC(\la;\Rl).
\end{CD}
\end{equation*}
\end{lemma}
\begin{proof} Consider the following diagram.
\begin{equation*} 
\xymatrix{
{\RLR(\la;R)} \ar[rrr]^{\il_R} \ar[dd]_{\psib_R} \ar[dr]^{\tr}
& & & {\RLR(\la;\Rl)} \ar[dd]^{\psib_{\Rl}}
\ar[dl]_{\tr} \\
 & {\CLR(\la^t;R^t)} \ar[r]^{\ick} \ar[dl]^{\phib_{R^t}} &
 {\CLR(\la^t;{\Rl}^t)} \ar[dr]_{\phib_{{\Rl}^t}} & \\
{\RC(\la;R)} \ar[rrr]_{\jl_R} & & & {\RC(\la;\Rl)}
}
\end{equation*}
Let us view this diagram as a prism in which the large rectangular
face is the front, and the other faces with four sides are the top and
bottom, and the faces with three sides are the left and right.
We want to show that the front face commutes.  For this it suffices to show
that all other faces commute.  The left and right faces commute by the
definition of $\psib$.  Let us define
$\ick:\CLR(\la^t;R^t)\rightarrow\CLR(\la^t;{\Rl}^t)$ so that
the top face commutes.  It suffices to show the bottom face commutes.
Observe that $\ick$ is the embedding for $\CLR$ that splits off
the first column of the first rectangle in $R^t$.
But then the bottom face commutes by \cite[Lemma 5.4]{KSS} applied
to $R^t$ in place of $R$.
\end{proof}

\begin{lemma} \label{psib perm} The following diagram commutes:
\begin{equation*}
\begin{CD}
  \RLR(\la;R) @>{s_p}>> \RLR(\la;s_p) \\
  @V{\psib_R}VV	@VV{\psib_{s_p R}}V \\
  \RC(\la;R) @= \RC(\la;s_p R).
\end{CD}
\end{equation*}
\end{lemma}
\begin{proof} We use the same kind of diagram as in the previous lemma.
Of course $(s_p R)^t = s_p (R^t)$.
\begin{equation*} 
\xymatrix{
{\RLR(\la;R)} \ar[rrr]^{s_p} \ar[dd]_{\psib_R} \ar[dr]^{\tr}
& & & {\RLR(\la;s_p R)} \ar[dd]^{\psib_{s_p R}}
\ar[dl]_{\tr} \\
 & {\CLR(\la^t;R^t)} \ar[r]^{s_p} \ar[dl]^{\phib_{R^t}} &
 {\CLR(\la^t;(s_p R)^t)} \ar[dr]_{\phib_{(s_p R)^t}} & \\
{\RC(\la;R)} \ar[rrr]_{=} & & & {\RC(\la;s_p R)}
}
\end{equation*}
We argue as in the previous lemma.  The left and right faces commute
by the definition of $\psib$, the top face commutes by
\cite[Proposition 32]{S2},
and the bottom face commutes by \cite[Lemma 8.5]{KSS}.
\end{proof}

\begin{proof}[Proof of Proposition \ref{psib embed}] Consider the diagram
\begin{equation*} 
\xymatrix{
{\RLR(\la;R)} \ar[rrr]^{i_R} \ar[dd]_{\psib_R} \ar[dr]^{\tr}
& & & {\RLR(\la;\rows(R))} \ar[dd]^{\psib_{\rows(R)}}
\ar[dl]_{\tr} \\
 & {\CLR(\la^t;R^t)} \ar[r]^{\Ick} \ar[dl]^{\phib_{R^t}} &
 {\CLR(\la^t;\rows(R)^t)} \ar[dr]_{\phib_{\rows(R)^t}} & \\
{\RC(\la;R)} \ar[rrr]_{j_R} & & & {\RC(\la;\rows(R))}
}
\end{equation*}
The left and right faces commute by the definition of $\psib$.
Let us define $\Ick$ so that the top face
commutes.  It suffices to show the bottom face commutes.
By the previous two lemmas, the bottom face commutes if $j_R$
is given by the composition of maps of the form $\jl_R$ and
the identity map, corresponding to the way that $i_R$ was computed.
But it is easy to see that the effect of this composition of maps
is precisely $j_R$.
\end{proof}

By the definition of $j_R$ and Definition \ref{def levrc} of the
level-restriction for rigged configurations, we have
\begin{equation} \label{lev RC}
  \RC^\ell(\la;R) = \{(\nu,J)\in\RC(\la;R) \mid
  j_R(\nu,J) \in \RC^\ell(\la;\rows(R))\}.
\end{equation}

We now show that Theorem \ref{rect res} follows from the special case
when $R$ consists of single rows.  The proof is a diagram chase using the 
the commutative diagram in Proposition \ref{psib embed}.
Since $\rows(R)$ consists of single rows, it is assumed that
$\psib_{\rows(R)}:\LRT(\la;\rows(R))\rightarrow\RC(\la;\rows(R))$
restricts to a bijection
$\LRT^\ell(\la;\rows(R))\rightarrow \RC^\ell(\la;\rows(R))$.
In particular
$\psib_{\rows(R)}(\LRT^\ell(\la;\rows(R)))=\RC^\ell(\la;\rows(R))$.
Since $\psib_R:\LRT(\la;R)\rightarrow\RC(\la;R)$ is a bijection,
it is enough to show that $\psib_R(\LRT^\ell(\la;R))=\RC^\ell(\la;R)$.
For the inclusion
$\psib_R(\LRT^\ell(\la;R))\subset\RC^\ell(\la;R)$,
suppose that $x\in\LRT^\ell(\la;R)$.  By \eqref{lev LR}
$i_R(x)\in \LRT^\ell(\la;\rows(R))$.  By assumption,
$\psib_{\rows(R)}(i_R(x))\in\RC^\ell(\la;\rows(R))$.
But $\psib_{\rows(R)} \circ i_R = j_R \circ \psib_R$ by
Proposition \ref{psib embed}, so
$j_R(\psib_R(x))\in\RC^\ell(\la;\rows(R))$.
By \eqref{lev RC}, $\psib_R(x)\in\RC^\ell(\la;R)$.
For the other inclusion, suppose $y\in \RC^\ell(\la;R)$.
Let $x\in\LRT(\la;R)$ be the unique
element such that $\psib_R(x)=y$.  Now
$\psib_{\rows(R)}(i_R(x))=j_R(\psib_R(x))=j_R(y)$.
By \eqref{lev RC} $j_R(y)\in\RC^\ell(\la;\rows(R))$.
By assumption $i_R(x)\in\LRT^\ell(\la;\rows(R))$.
By \eqref{lev LR} $x\in \LRT^\ell(\la;R)$, that is,
$y\in\psib_R(\LRT^\ell(\la;R))$.

\subsection{Single row quantum number bijection}

We must prove Theorem \ref{rect res} when $R$ consists of single rows.
For the rest of the paper we shall assume this is the case.
Then $\eta_j=1$ for all $j$, $R_j=(\mu_j)$ for $1\le j\le L$,
$\LRT(\la;R)=\CST(\la;\mu)$, $\LRT^\ell(\la;R)=\CST^\ell(\la;\mu)$,
and $R^t$ consists of single columns.
We also write $\psib_\mu$ for $\psib_R$ in this case.
Again using $\std$ we identify $\LRT(\la;R)$ with $\RLR(\la;R)$,
and $\LRT^\ell(\la;R)$ with its image in $\RLR(\la;R)$ under $\std$.
Now \cite[Section 4.2]{KSS} gives a direct description
of $\phib_{R^t}$ that is particularly simple when $R^t$ consists
of single columns.  This is easily translated to the following algorithm
to compute the bijection $\psib_R:\RLR(\la;R)\rightarrow\RC(\la;R)$.
First, $\nu\in\Conf(\la;R)$ requires that
\begin{equation*}
  |\nu^{(k)}| = \sum_{j>k} \la_j
\end{equation*}
for $k\ge 1$.  The vacancy numbers may be given by
\begin{equation*}
  P^{(k)}_i(\nu) = Q_i(\nu^{(k-1)}) - 2 Q_i(\nu^{(k)}) +
   Q_i(\nu^{(k+1)})
\end{equation*}
where $\nu^{(0)}=\mu$ and (since $\mu$ is not necessarily a partition)
\begin{equation*}
Q_i(\mu) := \sum_j \min\{\mu_j,i\}.
\end{equation*}

Now let us describe the bijection $\psib_R:\RLR(\la;R) \rightarrow
\RC(\la;R)$.  Start with $T\in\RLR(\la;R)$.  Write $T^-$ for the
tableau obtained by removing the largest letter from $T$
(which occurs in row $r$, say) and $\la^-$ for the shape of $T^-$.
Let $R^-=((\mu_1),(\mu_2),\dotsc,(\mu_{L-1}),(\mu_L-1))$.
Since $T^-\in\RLR(\la^-;R^-)$, by induction
$\psib_R(T^-)=(\nub,\Jb)$ is defined.  Let $s^{(r)}= \infty$.
For $k=r-1$ down to $1$, select the longest singular string in
$(\nub,\Jb)^{(k)}$ of length $s^{(k)}$ (possibly of zero length)
such that $s^{(k)} \le s^{(k+1)}$.
With the convention $s^{(0)}=\mu_L-1$,
it can be shown that $s^{(0)} \le s^{(1)}$ as well.
Then $\psib_R(T):=(\nu,J)$ is obtained from $(\nub,\Jb)$ by
lengthening each of the selected strings by one,
and resetting their labels to make them singular with respect
to the vacancy numbers in the definition of $\RC(\la;R)$,
and leaving all other strings unchanged.
Denote the transformation $(\nub,\Jb)\to (\nu,J)$ by $\db^{-1}$.

The inverse of $\db^{-1}$, denoted $\db$, is obtained as follows.
Set $\ell^{(0)}=\mu_L$. Select inductively a singular string of length
$\ell^{(k)}$ in $(\nu,J)^{(k)}$ with $\ell^{(k)}$ smallest such that
$\ell^{(k)}\ge \ell^{(k-1)}$. If no such singular string exists set
$\ell^{(k)}=\infty$. Then $(\nub,\Jb)$ is obtained from $(\nu,J)$ by 
shortening all selected strings by one, making them singular again and 
leaving all other strings unchanged.

\begin{rem} Up to the relabeling bijection $\std$ this
is precisely the description of the bijection
$\CST(\la;\mu)\rightarrow \RC(\la;(\mu_1),\dots,(\mu_L))$
that was given in terms of the map called $\pi_*$ in \cite{KR}.
\end{rem}

\begin{ex} Take $\mu=(2,2,2,2,1)$, $\la=(3,3,2,1)$ and
\begin{equation*}
T=\begin{matrix} 1&2&6\\3&4&8\\5&9&\\7&&\end{matrix}
\quad\text{so that}\quad
T^-=\begin{matrix} 1&2&6\\3&4&8\\5&&\\7&&\end{matrix}
\end{equation*}
and $r=3$. The rigged configuration corresponding to $T^-$ is
\begin{equation*}
\begin{array}{r|c|c|l} \cline{2-3} 0&&&0\\ 
 \cline{2-3} 0&&&0\\ \cline{2-3} 0&*&\multicolumn{2}{l}{0}\\ \cline{2-2}
\end{array}
\qquad
\begin{array}{r|c|l} \cline{2-2} 0&*&0\\ \cline{2-2}0&&0\\ \cline{2-2} 
\end{array}
\qquad
\begin{array}{r|c|l} \cline{2-2} 0&&0\\ \cline{2-2} \end{array}
\end{equation*}
where the labels are written to the right of each part and the vacancy
numbers to the left. The selected strings under $\db^{-1}$ with
$r=3$ are indicated by $*$. Hence the rigged configuration corresponding
to $T$ is
\begin{equation*}
\begin{array}{r|c|c|l} \cline{2-3} 0&&&0\\ 
 \cline{2-3} 0&&&0\\ \cline{2-3} 0&&&0\\ \cline{2-3}
\end{array}
\qquad
\begin{array}{r|c|c|l} \cline{2-3} 1&&&1\\ 
 \cline{2-3} 0&&\multicolumn{2}{l}{0} \\ \cline{2-2} 
\end{array}
\qquad
\begin{array}{r|c|l} \cline{2-2} 0&&0\\ \cline{2-2} \end{array}.
\end{equation*}
\end{ex}

\subsection{Proof of the single row case}

Now we come to the proof of Theorem \ref{rect res} when $R$
is a sequence of single rows. More precisely we will prove the
following theorem.

\begin{theorem}\label{thm bound}
Let $\la=(\la_1,\la_2,\dots,\la_n)$ be a partition of level $\ell$
and $\mu=(\mu_1,\mu_2,\ldots,\mu_L)$ an array
of positive integers not exceeding $\ell$.
Then $(\nu,J)$ is in the image of $\CST^\ell(\la,\mu)$ under 
$\psib_\mu$ if and only if
\begin{enumerate}
\item \label{c1} $\nu_1^{(k)}\le \ell$ for all $1\le k<n$, and 
\item \label{c2} there exists a column-strict tableau $t\in\CST(\la')$
such that
\begin{equation}\label{bound}
\x_i^{(k)}\le P_i^{(k)}(\nu)
 -\sum_{j=1}^{\la_k-\la_n}\chi(i\ge \lt+t_{j,k})
 +\sum_{j=1}^{\la_{k+1}-\la_n}\chi(i\ge \lt+t_{j,k+1})
\end{equation}
for all $1\le k<n$ and $1\le i\le \ell$.
\end{enumerate}
\end{theorem}

\begin{rem}\label{rem col}
The first column of $t\in \CST(\la')$ has length $\la_1-\la_n$. 
Since $t$ is a column-strict tableau over the alphabet 
$\{1,2,\ldots,\la_1-\la_n\}$ this requires that $t_{j,1}=j$.
\end{rem}

\begin{rem}\label{rem pic}
Since $t_{j,k}\le t_{j,k+1}$ the bounds in \eqref{bound} can be rewritten
as
\begin{multline}\label{bound re}
\x_i^{(k)}\le P_i^{(k)}(\nu)\\
 -\sum_{j=1}^{\la_{k+1}-\la_n} \chi(\lt+t_{j,k}\le i<\lt+t_{j,k+1})
 -\sum_{j=\la_{k+1}-\la_n+1}^{\la_k-\la_n} \chi(i\ge \lt+t_{j,k}).
\end{multline}
For the proof of Theorem \ref{thm bound} it will be useful to have
the following graphical description of \eqref{bound re} in mind.
Consider $n-1$ strips of length $\ell$ and height $\la_1-\la_n$
arranged on top of each other. Assign the label $k$ to the $k$-th strip
from the top. Within each strip assign a height label with height $1$ at
the bottom of the strip and height $\la_1-\la_n$ at the top of the strip.
Call the coordinate along the horizontal axis the position.
Then draw a horizontal line from position $\lt+t_{j,k}$ to position
$\lt+t_{j,k+1}$ at height $j$ in the $k$-th strip with a closed dot
at position $\lt+t_{j,k}$ and an open dot at position $\lt+t_{j,k+1}$
to indicate that the first position belongs to the line, whereas the
second one does not. If $t_{j,k}=t_{j,k+1}$ draw an open dot.
If $t_{j,k+1}$ does not exist draw a horizontal line from position
$\lt+t_{j,k}$, indicated by a closed dot, to position $\ell$.
If there is an open dot at position $\lt+t_{j,k+1}$ of height $j$
in strip $k$, then there is also a dot at the same position and height
in strip $k+1$. Connect all such dots by a vertical line.
This way one obtains $\la_1-\la_n$ paths which all end at position
$\ell$. The $j$-th path in strip 1 starts at position $\lt+j$ by 
Remark \ref{rem col}. Furthermore, since $t_{j,k}<t_{j+1,k}$ the paths 
do not intersect. The $k$-th strip contains $\la_k-\la_n$ paths.
The other $\la_1-\la_k$ paths already ended at position $\ell$
in previous strips.

An example for a set of such nonintersecting paths is given in Figure
\ref{fig paths}. It corresponds to $n=5$, $\ell=8$,
$\la=(6,5,4,3,1)$ and
\begin{equation*}
t=\begin{array}{cccc}
 1&1&2&4\\
 2&2&3&5\\
 3&4&4&\\
 4&5&&\\
 5&&&
\end{array}.
\end{equation*}
The dashed lines separate the various strips.

\begin{figure}
\scalebox{1.2}{
\begin{picture}(185,138)(-10,-10)
\Line(0,128)(175,128)
\Line(175,128)(175,0)
\Line(0,0)(175,0)
\Line(0,0)(0,128)
\DashLine(0,32)(175,32){3}
\DashLine(0,64)(175,64){3}
\DashLine(0,96)(175,96){3}
\Text(-4,112)[]{1}
\Text(-4,80)[]{2}
\Text(-4,48)[]{3}
\Text(-4,16)[]{4}
\Line(0,2)(0,-2)
\Line(25,2)(25,-2)
\Line(50,2)(50,-2)
\Line(75,2)(75,-2)
\Line(100,2)(100,-2)
\Line(125,2)(125,-2)
\Line(150,2)(150,-2)
\Line(175,2)(175,-2)
\Text(0,-8)[]{1}
\Text(25,-8)[]{2}
\Text(50,-8)[]{3}
\Text(75,-8)[]{4}
\Text(100,-8)[]{5}
\Text(125,-8)[]{6}
\Text(150,-8)[]{7}
\Text(175,-8)[]{8}
\Line(75,100)(75,68)
\CCirc(75,100){2}{Black}{White}
\Line(75,68)(100,68)
\CCirc(75,68){2}{Black}{Black}
\Line(100,68)(100,36)
\CCirc(100,68){2}{Black}{White}
\Line(100,36)(150,36)
\CCirc(100,36){2}{Black}{Black}
\Line(150,36)(150,4)
\CCirc(150,36){2}{Black}{White}
\Line(150,4)(175,4)
\CCirc(150,4){2}{Black}{Black}
\CCirc(175,4){2}{Black}{Black}
\Line(100,106)(100,74)
\Line(100,74)(125,74)
\Line(124,74)(124,42)
\Line(125,42)(174,42)
\Line(174,42)(174,10)
\CCirc(100,106){2}{Black}{White}
\CCirc(100,74){2}{Black}{Black}
\CCirc(125,74){2}{Black}{White}
\CCirc(125,42){2}{Black}{Black}
\CCirc(175,42){2}{Black}{White}
\CCirc(175,10){2}{Black}{Black}
\Line(125,112)(150,112)
\Line(150,112)(150,48)
\Line(150,48)(175,48)
\CCirc(125,112){2}{Black}{Black}
\CCirc(150,112){2}{Black}{White}
\CCirc(150,80){2}{Black}{White}
\CCirc(150,48){2}{Black}{Black}
\CCirc(175,48){2}{Black}{Black}
\Line(150,118)(175,118)
\Line(174,118)(174,86)
\CCirc(150,118){2}{Black}{Black}
\CCirc(175,118){2}{Black}{White}
\CCirc(175,86){2}{Black}{Black}
\CCirc(175,124){2}{Black}{Black}
\SetColor{Red}
\Line(125,0)(125,128)
\end{picture}
}
\caption{An example for nonintersecting paths illustrating \eqref{bound}
\label{fig paths}}
\end{figure}
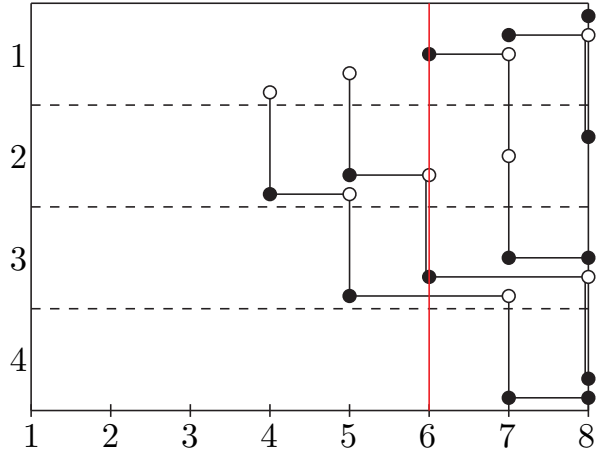

To read off the bound on $\x_i^{(k)}$ from the picture, draw a vertical 
line at position $i$. Suppose that $m$ paths cross this line horizontally 
in strip $k$ (when the vertical line goes through a closed/open dot we 
consider this as crossing/not crossing). Then
$P_i^{(k)}(\nu)-m$ is the maximal possible rigging for strings of
length $i$ in $(\nu,J)^{(k)}$. For example, the vertical line at 
position 6 in Figure \ref{fig paths} crosses one line in strip 1, 
no line in strip 2 and 4, and two lines in strip 3, so that 
$\x_6^{(1)}\le P_6^{(1)}(\nu)-1$,
$\x_6^{(2)}\le P_6^{(2)}(\nu)$, $\x_6^{(3)}\le P_6^{(3)}(\nu)-2$, 
and $\x_6^{(4)}\le P_6^{(4)}(\nu)$.

Recall that the rigging of a singular string of length $i$ in
$(\nu,J)^{(k)}$ equals the vacancy number $P_i^{(k)}(\nu)$.
Hence the above graphical description of the bounds shows
that $(\nu,J)^{(k)}$ cannot contain singular strings of 
length $i$ in the intervals
\begin{equation*}
\begin{aligned}[2]
&\lt+t_{j,k}\le i<\lt+t_{j,k+1} &\qquad&
 \text{for $1\le j\le \la_{k+1}-\la_n$,}\\
\text{and}\quad &\lt+t_{\la_{k+1}-\la_n+1,k}\le i
 &&\text{if $\la_{k+1}<\la_k$,} 
\end{aligned}
\end{equation*}
since in these intervals a vertical line at position $i$ would cross
at least one path. 
Conversely, if $i$ is the length of a singular string in
$(\nu,J)^{(k)}$ then it must be in the complements of these intervals,
that is,
\begin{equation}\label{sing cond}
\begin{aligned}[2]
& 1\le i<\lt+t_{1,k} &&\\
\text{or}\quad &
  \lt+t_{j-1,k+1}\le i<\lt+t_{j,k} &\qquad& 
  \text{for $1<j\le \la_{k+1}-\la_n$,}\\
\text{or}\quad &
  \lt+t_{\la_{k+1}-\la_n,k+1}\le i<\lt+t_{\la_{k+1}-\la_n+1,k} &\qquad& 
  \text{if $\la_{k+1}<\la_k$,}\\
\text{or}\quad &
  \lt+t_{\la_{k+1}-\la_n,k+1}\le i\le \ell &\qquad&
  \text{if $\la_{k+1}=\la_k$.}\\
\end{aligned}
\end{equation}
Since $t_{j,k}\le t_{j,k+1}$ these intervals are pairwise disjoint, but
some of these intervals can of course be empty.
Graphically the conditions in \eqref{sing cond} require that $i$
lies between two paths. More precisely, the first case in \eqref{sing cond} 
states that $i$ lies to the left of the first path, the second condition
requires that $i$ lies between the $(j-1)$-th and $j$-th path,
and the third case applies if there are more than $\la_{k+1}-\la_n$ paths
in the $k$-th strip in which case $i$ lies between paths $\la_{k+1}-\la_n$
and $\la_{k+1}-\la_n+1$. The last condition applies if there are
exactly $\la_{k+1}-\la_n$ paths in strip $k$. None of these ends at $\ell$
in this strip and the condition implies that $i$ lies to the right of
the rightmost path.
\end{rem}

\begin{rem}
We use the following conventions throughout the proof: 
$t_{0,k}=-\lt$ and $t_{j,k}=\la_1-\la_n+1$ for $j>\la_k-\la_n$.
\end{rem}

Without further ado we present the gory details of the proof of Theorem
\ref{thm bound}.

\subsection*{Proof of Theorem \ref{thm bound}}
We prove the theorem by induction on $|\la|$. The theorem is true for
$\la=\emptyset$ since then $T=\emptyset$ and $(\nu,J)=\emptyset$.
In this case $T$ is of level $\ell\ge 0$ and conditions 
\ref{c1} and \ref{c2} are trivially satisfied.

\subsection*{Proof of the forward direction}
Let $T\in\CST^\ell(\la,\mu)$
and $(\nu,J)=\psib_\mu(T)$ its image under the row-wise quantum number 
bijection.
Let $T^-$ be the tableau obtained from $T$ by removing the 
rightmost largest entry. Set $\lab=\shape(T^-)$, $(\nub,\Jb)=\db(\nu,J)$ and
denote by $r$ the row index of the cell $\la/\lab$.

Set $\la^0=\shape(T^{-^{\mu_L}})$.
The tableau $T^-$ is of level $\ell$ since $\lab_1-\la^0_n\le \la_1-\la^0_n
\le \ell$ by the condition that $T$ is of level $\ell$.
By induction the theorem holds for $T^-$ so that $\nub_1^{(k)}\le \ell$ 
and there exists a column-strict tableau $\tb$ of shape 
$(\lab_1-\lab_n,\dots,\lab_{n-1}-\lab_n)^t$ such that by \eqref{bound re}
\begin{multline}\label{bound bar}
\xb_i^{(k)}\le P_i^{(k)}(\nub)\\
 -\sum_{j=1}^{\lab_{k+1}-\lab_n}\chi(\ltb+\tb_{j,k}\le i<\ltb+\tb_{j,k+1})
 -\sum_{j=\lab_{k+1}-\lab_n+1}^{\lab_k-\lab_n}\chi(i\ge \ltb+\tb_{j,k}).
\end{multline}
for all $1\le k<n$ and $1\le i\le \ell$.
Here $\ltb=\ell-\lab_1+\lab_n$, and $\xb_i^{(k)}$ is the largest part of 
$\Jb_i^{(k)}$ and zero if $\Jb_i^{(k)}$ is empty. 
The aim is to show that conditions \ref{c1} and \ref{c2} of the theorem 
hold for $(\nu,J)$.

Denote by $\s^{(k)}$ the length of the selected singular 
string in $(\nub,\Jb)^{(k)}$ under $\db^{-1}$. By definition
$\mu_L-1=\s^{(0)}\le \s^{(1)}\le \cdots\le \s^{(r-1)}$ and 
$\s^{(k)}=\infty$ for $k\ge r$.
We claim that there exist indices $j^{(k)}$ for $0\le k<r$ such that
\begin{align}
\label{jk}
\ltb+\tb_{j^{(k)}-1,k+1} &\le \s^{(k)} <\ltb+\tb_{j^{(k)},k}
 \qquad \text{for $1\le k<r$}\\
\label{j0}
\ltb+\tb_{j^{(0)}-1,1} &\le s^{(0)}<\ltb+\tb_{j^{(0)},1}
\end{align}
and
\begin{equation}\label{order j}
1\le j^{(0)}\le j^{(1)}\le \cdots\le j^{(r-1)}\le \la_r-\la_n+\delta_{r,n},
\end{equation}
where by definition $\tb_{0,k}=-\ltb$.
The proof proceeds by descending induction on $k$ for
$1\le k<r$. We make frequent use of \eqref{sing cond} applied to
$(\nub,\Jb)$ where the first and third line are viewed as
the cases $j=1$ and $j=\lab_{k+1}-\lab_n+1$ of the general interval
appearing in the second line of \eqref{sing cond}.
First assume $k=r-1$. Note that $\lab_r<\lab_{r-1}$ since 
$\lab_r+1=\la_r\le \la_{r-1}=\lab_{r-1}$. Hence the existence of 
$j^{(r-1)}$ follows from \eqref{sing cond} since the last case 
does not apply. In particular we have $j^{(r-1)}\le \lab_r-\lab_n+1
=\la_r-\la_n+\delta_{r,n}$. Now consider $0\le k<r-1$ and
assume that $\ltb+\tb_{j^{(k)}-1,k+1}\le \s^{(k)}$ for some 
$j^{(k)}>j^{(k+1)}$. Then by induction and the column-strictness of $\tb$
\begin{equation*}
\s^{(k+1)}<\ltb+\tb_{j^{(k+1)},k+1}\le \ltb+\tb_{j^{(k)}-1,k+1}
 \le \s^{(k)}
\end{equation*}
which is a contradiction. Hence $j^{(k)}\le j^{(k+1)}$.
Since $\s^{(k)}$ is the length of a singular string and by induction 
$j^{(k)}\le \lab_r-\lab_n+1\le\lab_{k+1}-\lab_n$ for $1\le k<r-1$,
$s^{(k)}$ must be in the first or second set of the intervals in
\eqref{sing cond} with all quantities replaced by their barred counterparts
which proves \eqref{jk}. Equation \eqref{j0} follows since $\tb_{j,1}=j$
by Remark \ref{rem col}.

Let us now prove condition \ref{c1} of the theorem.
By construction $(\nu,J)$ is obtained from $(\nub,\Jb)$ by increasing 
the length of the selected strings in $(\nub,\Jb)^{(k)}$ by one for $1\le k<r$,
making them singular again and leaving all other strings unchanged.
For $r=1$ this means that $(\nu,J)=(\nub,\Jb)$ so that condition \ref{c1} of 
the theorem is satisfied by induction. 
Now assume $r>1$. Since $\tb_{j^{(r-1)},k}\in\{1,2,\ldots,\lab_1-\lab_n\}$
it follows from \eqref{jk} with $k=r-1$ that
$\s^{(r-1)}<\ltb+\tb_{j^{(r-1)},r-1}\le \ltb+\lab_1-\lab_n=\ell$.
Hence $\mu_L-1\le \s^{(1)}\le \cdots \le \s^{(r-1)}<\ell$
which ensures condition \ref{c1} of the theorem for $1<r\le n$.

It remains to prove that the second condition of the theorem holds.
The vacancy numbers of $\nu$ and $\nub$ are related as follows
\begin{equation*}
P_i^{(k)}(\nub)=P_i^{(k)}(\nu)-\chi(\s^{(k-1)}<i\le \s^{(k)})
 +\chi(\s^{(k)}<i\le \s^{(k+1)}).
\end{equation*}
By construction $\x_i^{(k)}\le \xb_i^{(k)}$ for $i\neq \s^{(k)}+1$
and $\x_{\s^{(k)}+1}^{(k)}=P_{\s^{(k)}+1}^{(k)}(\nu)$ for
$1\le k<r$. Hence
\begin{multline}\label{bound inter}
\x_i^{(k)}\le P_i^{(k)}(\nu)-\chi(\s^{(k-1)}<i\le \s^{(k)})
 +\chi(\s^{(k)}<i\le \s^{(k+1)})\\
 -\sum_{j=1}^{\lab_k-\lab_n}\chi(i\ge \ltb+\tb_{j,k})
 +\sum_{j=1}^{\lab_{k+1}-\lab_n}\chi(i\ge \ltb+\tb_{j,k+1})
\end{multline}
for $i\neq \s^{(k)}+1$.
In the remainder of the proof of the forward direction it will
be shown that \eqref{bound} holds for $i=s^{(k)}+1$ and that
\eqref{bound inter} implies \eqref{bound} for $i\neq \s^{(k)}+1$.
We distinguish the cases $r=1$, $1<r<n$, and $r=n$.

\subsection*{Case $r=1$.}
In this case $\lab_1=\la_1-1$, $\lab_k=\la_k$ for
$1<k\le n$, $\ltb=\lt+1$ and $\tb$ is a column-strict tableau over
the alphabet $\{1,2,\ldots,\la_1-\la_n-1\}$.
Furthermore $\s^{(0)}=\mu_L-1$ and $\s^{(k)}=\infty$ for $k\ge 1$. 
By \eqref{bound inter} we have for all $1\le k<n$ and $1\le i\le \ell$
\begin{multline}\label{bound r=1}
\x_i^{(k)}\le P_i^{(k)}(\nu)-\chi(i\ge \mu_L)\delta_{k,1}
 -\sum_{j=1}^{\la_k-\la_n-\delta_{k,1}}\chi(i\ge \lt+1+\tb_{j,k})\\
 +\sum_{j=1}^{\la_{k+1}-\la_n}\chi(i\ge \lt+1+\tb_{j,k+1}).
\end{multline}

Remark \ref{rem col} requires that $t_{j,1}=j$ for
$1\le j\le \la_1-\la_n$. Hence
\begin{multline}\label{ineq inter}
-\chi(i\ge \mu_L)-\sum_{j=1}^{\la_1-\la_n-1}\chi(i\ge \lt+1+\tb_{j,1})\\
\le -\sum_{j=1}^{\la_1-\la_n}\chi(i\ge \lt+t_{j,1})
+\chi(\lt+1\le i<\mu_L).
\end{multline}
If $\mu_L\le\lt+1$ the term $\chi(\lt+1\le i<\mu_L)$ vanishes.
In this case set $t_{j,k}=\tb_{j,k}+1$ for $1<k<n$ and $1\le j\le \la_k-\la_n$ 
which defines a column-strict tableau of shape 
$(\la_1-\la_n,\ldots,\la_{n-1}-\la_n)^t$ over $\{1,2,\ldots,\la_1-\la_n\}$.
Then \eqref{bound r=1} implies \eqref{bound}.

The case $\mu_L>\lt+1$ is considerably harder to establish
due to the extra term $\chi(\lt+1\le i<\mu_L)$.
Our strategy is as follows. The term $\chi(\lt+1\le i<\mu_L)$ can be
absorbed by defining $t_{j,2}$ appropriately except in certain cases. 
In general this introduces extra terms for the bounds at $k=2$. 
These in turn can be absorbed by defining $t_{j,3}$ appropriately (except
in certain cases) and so on. If all $t_{j,k}$ for $1\le k<n$ can be defined 
and all bounds for $1\le k<n$ are written in the form of \eqref{bound} we are
done. In the exceptional cases (when \eqref{bound r=1} does not imply
\eqref{bound}) it can be shown that the corresponding 
tableau $T$ is not of level $\ell$ which contradicts the assumptions.

Let us now plunge into the details.
Define $t_{j,1}=j$ for $1\le j\le \la_1-\la_n$ and set $d=\mu_L-\lt-1$. 
Since $r=1$, we have $\xb_i^{(k)}=\x_i^{(k)}$ and $P_i^{(k)}(\nub,\tb)$ 
equals the right-hand side of \eqref{bound r=1}. Let $a_j^{(1)}$ for
$1\le j\le \la_2-\la_n+1$ be the minimal index
$i\in [\tb_{j-1,2}+1,\tb_{j,2}]\cap [1,d]$
such that $\xb_{\lt+i}^{(1)}=P_{\lt+i}^{(1)}(\nub,\tb)$,
where $\tb_{0,k}=-\ltb$ and $\tb_{j,k}=\lab_1-\lab_n+1$
for $j>\lab_k-\lab_n$.
If no such $i$ exists set $a_j^{(1)}=\tb_{j,2}+1$.
By definition $\xb_{\lt+i}^{(1)}<P_{\lt+i}^{(1)}(\nub,\tb)$ for 
$\tb_{j-1,2}<i<a_j^{(1)}$ and $1\le i\le d$ so that we can sharpen the bounds 
in \eqref{bound r=1} for $k=1$ by adding
$-\sum_{j=1}^{\la_2-\la_n+1}\chi(\lt+\tb_{j-1,2}<i<\lt+a_j^{(1)})
\chi(\lt+1\le i\le \lt+d)$. Note that $\tb_{\la_2-\la_n+1,2}+1=
\lab_1-\lab_n+2=\la_1-\la_n+1=t_{\la_2-\la_n+1,2}$.
The case $a_{\la_2-\la_n+1}^{(1)} < t_{\la_2-\la_n+1,2}$ will be dealt with
later. Suppose that $a_{\la_2-\la_n+1}^{(1)} = t_{\la_2-\la_n+1,2}$.
Then one finds using \eqref{ineq inter}
\begin{multline}\label{bound r1k1}
\x_i^{(1)}\le P_i^{(1)}(\nu)-\sum_{j=1}^{\la_1-\la_n}\chi(i\ge \lt+t_{j,1})
+\sum_{j=1}^{\la_2-\la_n}\chi(i\ge \lt+1+\tb_{j,2})\\
+\sum_{j=1}^{\la_2-\la_n}\chi(\lt+a_j^{(1)}\le i\le \lt+\tb_{j,2}).
\end{multline}
Define $t_{j,2}=\min\{a_j^{(1)},t_{j+1,2}-1\}$ recursively by 
descending $1\le j\le \la_2-\la_n$. 
{}From its definition it is clear that $a_j^{(1)}$ lies in the interval
$[\tb_{j-1,2}+1,\tb_{j,2}+1]$. By descending induction on $j$ it also
follows that $t_{j,2}\in [\tb_{j-1,2}+1,\tb_{j,2}+1]$ and that either
$t_{j,2}=a_j^{(1)}$ or $a_j^{(1)}-1$. The latter case only occurs when 
$a_j^{(1)}=t_{j+1,2}=\tb_{j,2}+1$. In addition there must exist an index 
$j'>j$ such that $a_{j'}^{(1)}=\tb_{j'-1,2}+1$ if $t_{j,2}=a_j^{(1)}-1$.
This is because $t_{j+1,2}=\tb_{j,2}+1$ is at its lower bound in
the interval $[\tb_{j,2}+1,\tb_{j+1,2}+1]$ and this can happen in only two
ways; either $t_{j+1,2}=a_{j+1}^{(1)}$ which proves the assertion with
$j'=j+1$ or $t_{j+1,2}=a_{j+1}^{(1)}-1$ in which case the assertion must
be true by induction since the initial case is
$t_{\la_2-\la_n+1,2}=a_{\la_2-\la_n+1}^{(1)}$.
Note that it also follows by induction that 
$a_{j'}^{(1)}=a_j^{(1)}+j'-j-1$ if $j'$ is minimal.
{}From its definition it follows that $t_{j,2}<t_{j+1,2}$ and furthermore
$t_{j,2}\ge \tb_{j-1,2}+1\ge \tb_{j-1,1}+1=j$
which are the conditions for column-strictness for the first two
columns of $t$. Hence \eqref{bound r1k1} yields \eqref{bound} for $k=1$.

We proceed inductively on $1<k<n$. Assume that by induction
$t_{j,k'}\in [\tb_{j-1,k'}+1,\tb_{j,k'}+1]$ is already defined for 
$1\le k'\le k$. In terms of $t_{j,k}$ \eqref{bound r=1} reads
\begin{multline}\label{hh}
\x_i^{(k)}\le P_i^{(k)}(\nu)-\sum_{j=1}^{\la_k-\la_n}\chi(i\ge \lt+t_{j,k})
+\sum_{j=1}^{\la_{k+1}-\la_n}\chi(i\ge \lt+1+\tb_{j,k+1})\\
+\sum_{j=1}^{\la_k-\la_n}\chi(\lt+t_{j,k}\le i\le \lt+\tb_{j,k}).
\end{multline}

For $1<k<n$ define $a_j^{(k)}$ for $1\le j\le \la_{k+1}-\la_n+1$ to be 
the minimal index $i\in [\tb_{j-1,k+1}+1,\tb_{j,k+1}]\cap 
 \bigcup_{h=1}^{\la_k-\la_n}[t_{h,k},\tb_{h,k}]$
such that $\xb_{\lt+i}^{(k)}=P_{\lt+i}^{(k)}(\nub,\tb)$.
If no such $i$ exists set $a_j^{(k)}=\tb_{j,k+1}+1$.
Note that $\tb_{\la_{k+1}-\la_n+1,k+1}+1=\la_1-\la_n+1=
t_{\la_{k+1}-\la_n+1,k+1}$. The case $a_{\la_{k+1}-\la_n+1}^{(k)}<
t_{\la_{k+1}-\la_n+1,k+1}$ will be dealt with later.
Now assume that $a_{\la_{k+1}-\la_n+1}^{(k)}=t_{\la_{k+1}-\la_n+1,k+1}$, and
define recursively $t_{j,k+1}=\min\{a_j^{(k)},t_{j+1,k+1}-1\}$ on descending 
$1\le j\le \la_{k+1}-\la_n$.
By definition $a_j^{(k)}\in [\tb_{j-1,k+1}+1,\tb_{j,k+1}+1]$. As in the
case $k=1$ it follows by descending induction on $j$ that 
$t_{j,k+1}\in [\tb_{j-1,k+1}+1,\tb_{j,k+1}+1]$ and that either
$t_{j,k+1}=a_j^{(k)}$ or $a_j^{(k)}-1$. 
By definition $t_{j,k+1}<t_{j+1,k+1}$. Let us now show that also
$t_{j,k}\le t_{j,k+1}$ which would prove the column-strictness of $t$.
By definition $t_{h,k}\le a_j^{(k)}\le \tb_{h,k}$ for some $h$. Assume
$h<j$. Then $\tb_{h,k}\ge a_j^{(k)}\ge \tb_{j-1,k+1}+1$ which violates
the column-strictness of $\tb$. Hence $h\ge j$. If $t_{j,k+1}=a_j^{(k)}$
then $t_{j,k+1}\ge t_{h,k}\ge t_{j,k}$ as desired. If $t_{j,k+1}=a_j^{(k)}-1$
then a problem can only occur if $h=j$ and $a_j^{(k)}=t_{j,k}$.
However in this case $a_j^{(k)}=t_{j,k+1}+1=t_{j+1,k+1}=\tb_{j,k+1}+1
>\tb_{j,k}=\tb_{h,k}$ which is a contradiction. 
This proves $t_{j,k}\le t_{j,k+1}$.
By the same arguments as in the case $k=1$ there must exist
an index $j'>j$ such that $a_{j'}^{(k)}=\tb_{j'-1,k+1}+1$ if
$t_{j,k+1}=a_j^{(k)}-1$. For minimal $j'$ it follows again by
induction that $a_{j'}^{(k)}=a_j^{(k)}+j'-j-1$.

Since by definition $\x_{\lt+i}^{(k)}=\xb_{\lt+i}^{(k)}
<P_{\lt+i}^{(k)}(\nub,\tb)$ for $\tb_{j-1,k+1}<i<a_j^{(k)}$
and $t_{h,k}\le i\le \tb_{h,k}$ for some $1\le h\le \la_k-\la_n$,
one can add 
\begin{equation*}
-\sum_{j=1}^{\la_{k+1}-\la_n+1}
\chi(\lt+\tb_{j-1,k+1}<i<\lt+a_j^{(k)})
\sum_{h=1}^{\la_k-\la_n}\chi(\lt+t_{h,k}\le i\le \lt+\tb_{h,k})
\end{equation*}
to \eqref{hh}. If $a_{\la_{k+1}-\la_n+1}^{(k)}=t_{\la_{k+1}-\la_n+1,k+1}$
then the sum of this term and 
$\sum_{j=1}^{\la_k-\la_n}\chi(\lt+t_{j,k}\le i\le \lt+\tb_{j,k})$
does not exceed $\sum_{j=1}^{\la_{k+1}-\la_n}\chi(\lt+a_j^{(k)}\le i\le
\lt+\tb_{j,k+1})$ and \eqref{bound} is proven for $1<k<n$.

It remains to treat the case when there exists a $1\le k<n$
such that $a_{\la_{k+1}-\la_n+1}^{(k)}<\la_1-\la_n+1$. Let $\kappa$
be minimal with this property.
We will show that in this case $T$ is not of level $\ell$ which contradicts
the assumptions.

We claim that there exist indices $h_k$ and $j_k$ for
$1\le k\le \kappa$ such that
\begin{align}
\label{j k}
 \tb_{j_k-1,k+1}+1&\le a_{j_k}^{(k)}\le \tb_{j_k,k+1},\\
\label{h k}
 t_{h_k,k}&\le a_{j_k}^{(k)}\le \tb_{h_k,k},
\end{align}
$h_1\ge h_2\ge \cdots\ge h_\kappa$ and $h_k\ge j_k\ge h_{k+1}$.
The inequalities \eqref{j k} and \eqref{h k} hold for $k=\kappa$ with 
$j_{\kappa}=\la_{\kappa+1}-\la_n+1$ and some $h_{\kappa}$ by the definition 
of $\kappa$. Now suppose that $k<\kappa$ and that $h_{k'}$ and $j_{k'}$ for 
$k<k'\le\kappa$ satisfying \eqref{j k} and \eqref{h k} have been
defined by induction.
Recall that either $t_{j,k+1}=a_j^{(k)}$ or $a_j^{(k)}-1$.
First assume that $t_{h_{k+1},k+1}=a_{h_{k+1}}^{(k)}$. This
implies in particular that $a_{h_{k+1}}^{(k)}= t_{h_{k+1},k+1}\le 
\tb_{h_{k+1},k+1}$ by \eqref{h k} and hence
$\tb_{h_{k+1}-1,k+1}+1\le a_{h_{k+1}}^{(k)}
\le \tb_{h_{k+1},k+1}$ by the definition of $a_j^{(k)}$. Set $j_k=h_{k+1}$ 
and choose $h_k$ such that \eqref{h k} holds which must be possible by the 
definition of $a_j^{(k)}$. Also $h_k\ge j_k=h_{k+1}$ since otherwise 
$\tb_{h_k,k}\le\tb_{j_k-1,k+1}$ by the column-strictness of $\tb$ which 
yields a contradiction since then \eqref{j k} and \eqref{h k} cannot
hold simultaneously.
Next assume $t_{h_{k+1},k+1}=a_{h_{k+1}}^{(k)}-1$. 
Let $j_k>h_{k+1}$ be minimal such that 
$a_{j_k}^{(k)}=\tb_{j_k-1,k+1}+1\le \tb_{j_k,k+1}$; the existence of
$j_k$ was proved before. In addition it was shown that
$a_{j_k}^{(k)}=a_{h_{k+1}}^{(k)}+j_k-h_{k+1}-1$.
The existence of $h_k$ follows again from the definition of $a_j^{(k)}$. 
As before $h_k\ge j_k\ge h_{k+1}$.

By definition $\x_{\lt+a_{j_k}^{(k)}}^{(k)}=\xb_{\lt+a_{j_k}^{(k)}}^{(k)}
=P^{(k)}_{\lt+a_{j_k}^{(k)}}(\nub,\tb)$.
Since $P_i^{(k)}(\nub,\tb)$
is given by the right-hand side of \eqref{hh}, it follows from
\eqref{j k}, \eqref{h k} and the fact that 
$t_{j,k}\in [\tb_{j-1,k}+1,\tb_{j,k}+1]$ that
\begin{equation}\label{JP}
\x_{\lt+a_{j_k}^{(k)}}^{(k)}=P_{\lt+a_{j_k}^{(k)}}^{(k)}(\nu)-h_k+j_k
\qquad \text{for $1\le k\le \kappa.$}
\end{equation}

Define $T^b=T^{-^{\mu_L-b}}$ for $0\le b\le \mu_L$ with corresponding 
rigged configurations $(\nu^b,J^b)=\db^{\mu_L-b}(\nu,J)$. 
Let $r_b$ be the row index of the cell $T^b/T^{b-1}$.
Denote the 
length of the selected string in $(\nu^b,J^b)^{(k)}$ under $\db$ 
by $\ell_b^{(k)}$. 
We claim that \eqref{JP} implies
\begin{equation}
\label{ell bound}
\ell^{(k)}_{\lt+j_k}\le \lt+a_{j_k}^{(k)} \quad \text{for $1\le k\le 
 \kappa$.}
\end{equation}
This is shown by induction on $k$. By construction
$b=\ell_b^{(0)}\le \ell_b^{(1)}\le \ell_b^{(2)}\le \cdots \le \ell_b^{(n-1)}$
and $\ell_1^{(k)}<\ell_2^{(k)}<\ldots<\ell_{\mu_L}^{(k)}$. In addition
\begin{equation}\label{vac change}
P_i^{(k)}(\nu^{b-1})=P_i^{(k)}(\nu^b)
-\chi(\ell_b^{(k-1)}\le i<\ell_b^{(k)})+\chi(\ell_b^{(k)}\le i<\ell_b^{(k+1)}).
\end{equation}

If $\ell_{\lt+a_{j_1}^{(1)}}^{(1)}\le \lt+a_{j_1}^{(1)}$ then
\eqref{ell bound} follows immediately since $j_1\le a_{j_1}^{(1)}$
and $\ell_{b-1}^{(1)}<\ell_b^{(1)}$. Hence assume
that $\ell_{\lt+a_{j_1}^{(1)}}^{(1)}>\lt+a_{j_1}^{(1)}$. 
Since $\ell_b^{(0)}=b$, the vacancy number at $i=\lt+a_{j_1}^{(1)}$
is decreased by one with each application of $\db$ until 
$\ell_b^{(1)}\le \lt+a_{j_1}^{(1)}$ because of \eqref{vac change} at $k=1$. 
By \eqref{JP} it takes $h_1-j_1$ applications
of $\db$ until there is a singular string of length $\lt+a_{j_1}^{(1)}$
in the first rigged partition. Since $a_{j_1}^{(1)}=h_1$ by
\eqref{h k} and Remark \ref{rem col} this means that the singular string
occurs at $b=\lt+a_{j_1}^{(1)}-h_1+j_1=\lt+j_1$ which 
proves \eqref{ell bound} at $k=1$. 

Now consider the cases $1<k\le \kappa$ and assume that \eqref{ell bound} holds 
for $k'<k$. First assume that $t_{h_k,k}=a_{h_k}^{(k-1)}$. In this case 
$j_{k-1}=h_k$ and by \eqref{h k} 
$a_{j_{k-1}}^{(k-1)}=t_{h_k,k}\le a_{j_k}^{(k)}$ so that
$\ell_{\lt+j_{k-1}}^{(k-1)}\le\lt+a_{j_k}^{(k)}$ by \eqref{ell bound} 
at $k-1$. 
If $\ell^{(k)}_{\lt+j_{k-1}}\le \lt+a_{j_k}^{(k)}$ there is nothing to show
since $\ell_{b-1}^{(k)}<\ell_b^{(k)}$ and $j_{k-1}\ge j_k$.
Hence assume that $\ell^{(k)}_{\lt+j_{k-1}}> \lt+a_{j_k}^{(k)}$.
Again by \eqref{JP} and \eqref{vac change} it takes $h_k-j_k$
applications of $\db$ until there is a singular string of length
$\lt+a_{j_k}^{(k)}$ in the $k$-th rigged partition. Since $j_{k-1}=h_k$
the singular string occurs at $b=\lt+j_{k-1}-(h_k-j_k)=\lt+j_k$ which
proves \eqref{ell bound}.

Next assume that $t_{h_k,k}=a_{h_k}^{(k-1)}-1$. Then
$a_{j_{k-1}}^{(k-1)}=a_{h_k}^{(k-1)}+j_{k-1}-h_k-1=t_{h_k,k}+j_{k-1}-h_k$ 
so that by \eqref{h k}
$a_{j_{k-1}}^{(k-1)}\le a_{j_k}^{(k)}+j_{k-1}-h_k$.
Since $\ell_{\lt+j_{k-1}}^{(k-1)}\le \lt+a_{j_{k-1}}^{(k-1)}$
it takes at most $j_{k-1}-h_k$ applications of $\db$ before there is
a singular string in the $(k-1)$-th rigged partition of length not
exceeding $\lt+a_{j_k}^{(k)}$. After that, by \eqref{JP} and 
\eqref{vac change}, it takes $h_k-j_k$ applications of $\db$
until there is a singular string of length $\lt+a_{j_k}^{(k)}$ in the 
$k$-th rigged partition. Hence altogether the existence of a
singular string of length $\lt+a_{j_k}^{(k)}$ is assured at
$b=\lt+j_{k-1}-(j_{k-1}-h_k)-(h_k-j_k)=\lt+j_k$ which 
concludes the proof of \eqref{ell bound}.

Recall that $j_\kappa=\la_{\kappa+1}-\la_n+1$. Therefore \eqref{ell bound} 
implies that $\ell^{(\kappa)}_b$ is finite for $1\le b\le \lt+\la_{\kappa+1}
-\la_n+1$. If $\ell_b^{(k)}$ is finite this means that $r_b>k$ so that
$r_b>\kappa$ for $1\le b\le \lt+\la_{\kappa+1}-\la_n+1$.
Since at most $\la_{\kappa+1}-\la_n$ boxes can be removed from 
$T^{\lt+\la_{\kappa+1}-\la_n+1}$ in rows with index $\kappa<r_b<n$
it follows that $r_{\lt+1}=n$. This implies
$\la^0_n\le \la_n-\lt-1$ where $\la^0=\shape(T^0)$. Hence 
$\la_1-\la_n^0\ge \ell+1$ which contradicts the assumption that $T$ is of 
level $\ell$. This concludes the proof of the case $r=1$.

\subsection*{Case $1<r<n$.}
In this case $\lab_r=\la_r-1$, $\lab_k=\la_k$
for $k\neq r$ and $\ltb=\lt$. It is convenient to introduce
$p^{(1)}=j^{(1)}-1$,
\begin{equation}\label{p}
p^{(k)}=\begin{cases}
\max\{j^{(k)}-1,p^{(k-1)}\} & 
 \text{for $\s^{(k)}<\lt +\tb_{j^{(k)},k-1}-1$,}\\
j^{(k)} & \text{for $\s^{(k)}\ge \lt+\tb_{j^{(k)},k-1}-1$,}
\end{cases}
\end{equation}
for $1< k<r$ and $p^{(r)}=\la_r-\la_n$. Note that for $1\le k<r$ either
$p^{(k)}=j^{(k)}$ or $p^{(k)}=j^{(k)}-1$ and that $p^{(k)}\le p^{(k+1)}$.
Define $t_{j,1}=j$ for $1\le j\le \la_1-\la_n$,
\begin{equation}
\label{t def}
t_{j,k}=\begin{cases} \tb_{j,k} & 
 \text{for $1\le j<j^{(k-1)}$ and $p^{(k)}< j\le \la_k-\la_n$,}\\
\s^{(k-1)}-\lt+1 & \text{for $j=j^{(k-1)}=p^{(k-1)}$,}\\
\max\{\tb_{j-1,k},\tb_{j,k-1}\} & \text{for $p^{(k-1)}<j\le p^{(k)}$,}
\end{cases}
\end{equation}
for $1<k\le r$ and $t_{j,k}=\tb_{j,k}$ for $r<k<n$ and $1\le j\le \la_k-\la_n$.

By \eqref{jk} we have $\tb_{j^{(k-1)}-1,k}<\tb_{j^{(k-1)},k-1}$
for $1<k\le r$ so that 
\begin{equation}\label{t jk}
t_{j^{(k-1)},k}=\tb_{j^{(k-1)},k-1} \qquad
 \text{for $p^{(k-1)}=j^{(k-1)}-1<p^{(k)}$.}
\end{equation}

It needs to be shown that $t$ indeed defines a column-strict tableau
over the alphabet $\{1,2,\ldots,\la_1-\la_n\}$. Since
$\tb_{j,k}\in\{1,2,\ldots,\la_1-\la_n\}$ and $\s^{(k)}<\ell$ for all
$1\le k<r$ the condition $t_{j,k+1}\in\{1,2,\ldots,\la_1-\la_n\}$ might
only be violated if $p^{(k)}=j^{(k)}$ and $\s^{(k)}<\lt$
for $1\le k< r$. By \eqref{jk} the latter condition requires
$j^{(k)}=1$ so that $1=j^{(1)}=\dots=j^{(k)}$ by \eqref{order j}.
Since $0\le \s^{(1)}\le \cdots\le \s^{(k)}<\lt$ the first condition
in \eqref{p} applies for $p^{(h)}$ for $2\le h\le k$. However, since
$p^{(1)}=j^{(1)}-1=0$ this implies that $p^{(k)}=0$ which contradicts
the requirement $j^{(k)}=p^{(k)}$. This shows that 
$t_{j,k+1}\in\{1,2,\ldots,\la_1-\la_n\}$.

Next we check that $t$ is column-strict.
The condition $t_{j,k}<t_{j+1,k}$ only needs to be checked for
$1<k\le r$ and $j^{(k-1)}-1\le j\le p^{(k)}$ since in all other
cases it automatically follows from the column-strictness of $\tb$.
First assume $p^{(k-1)}=j^{(k-1)}$. Then 
$t_{j^{(k-1)}-1,k}=\tb_{j^{(k-1)}-1,k}<\s^{(k-1)}-\lt+1=t_{j^{(k-1)},k}$
by \eqref{jk}. Furthermore $t_{j^{(k-1)},k}=\s^{(k-1)}-\lt+1
\le \tb_{j^{(k-1)},k-1}<t_{j^{(k-1)}+1,k}$ by \eqref{jk} and \eqref{t def}.
Next assume $p^{(k-1)}=j^{(k-1)}-1$. Then for $p^{(k-1)}<p^{(k)}$,
$t_{j^{(k-1)}-1,k}=\tb_{j^{(k-1)}-1,k}<\tb_{j^{(k-1)},k-1}
=t_{j^{(k-1)},k}$ by \eqref{jk} and \eqref{t jk}. For $p^{(k-1)}=p^{(k)}$ 
the column-strictness is trivial.
Furthermore $\tb_{j-1,k}<\tb_{j,k}\le \max\{\tb_{j,k},\tb_{j+1,k-1}\}$
and $\tb_{j,k-1}< \tb_{j+1,k-1}\le \max\{\tb_{j,k},\tb_{j+1,k-1}\}$
so that $t_{j,k}< t_{j+1,k}$ for $p^{(k-1)}<j<p^{(k)}$.
And finally $t_{p^{(k)},k}=\max\{\tb_{p^{(k)}-1,k},\tb_{p^{(k)},k-1}\}
\le \tb_{p^{(k)},k}<\tb_{p^{(k)}+1,k}=t_{p^{(k)}+1,k}$. 

The conditions $t_{j,k}\le t_{j,k+1}$ only need to be verified
for $j^{(1)}\le j\le p^{(2)}$ and $k=1$, 
for $j^{(k-1)}\le j\le p^{(k+1)}$ and $1< k<r$, and for
$j^{(r-1)}\le j$ and $k=r$. First assume $k=1$. Then
$t_{j,1}=j\le \max\{\tb_{j-1,2},\tb_{j,1}\}=t_{j,2}$ for 
$j^{(1)}\le j\le p^{(2)}$. Now assume $1<k<r$. 
For $p^{(k-1)}=j^{(k-1)}<j^{(k)}$ we have 
$t_{j^{(k-1)},k}=\s^{(k-1)}-\lt+1\le \tb_{j^{(k-1)},k-1}
\le \tb_{j^{(k-1)},k+1}=t_{j^{(k-1)},k+1}$ by \eqref{jk}.
For $p^{(k-1)}=j^{(k-1)}=j^{(k)}$ we have
$t_{j^{(k-1)},k}=\s^{(k-1)}-\lt+1\le \s^{(k)}-\lt+1=t_{j^{(k)},k+1}$.
For $p^{(k-1)}<j<j^{(k)}$ one obtains 
$t_{j,k}=\max\{\tb_{j-1,k},\tb_{j,k-1}\}\le \tb_{j,k+1}=t_{j,k+1}$.
Next assume $p^{(k-1)}<p^{(k)}=j^{(k)}$. Then the second case
of \eqref{p} applies so that $\s^{(k)}-\lt+1\ge \tb_{j^{(k)},k-1}$. 
By \eqref{jk} also $\s^{(k)}-\lt+1\ge \tb_{j^{(k)}-1,k+1}\ge 
\tb_{j^{(k)}-1,k}$. This implies
$t_{j^{(k)},k}=\max\{\tb_{j^{(k)}-1,k},\tb_{j^{(k)},k-1}\}
\le \s^{(k)}-\lt+1=t_{j^{(k)},k+1}$.
And finally for $p^{(k)}<j\le p^{(k+1)}$ we have
$t_{j,k}=\tb_{j,k}\le \max\{\tb_{j-1,k+1},\tb_{j,k}\}=t_{j,k+1}$.
In a similar fashion one shows that $t_{j,r}\le t_{j,r+1}$.

Hence $t$ forms a column-strict tableau of shape 
$(\la_1-\la_n,\ldots,\la_{n-1}-\la_n)^t$.

We will now show that \eqref{bound} holds with $t$ as 
defined in \eqref{t def}. First assume that $i=\s^{(k)}+1$. 
In this case $\x_i^{(k)}=P_i^{(k)}(\nu)$ if $1\le k<r$.
Hence it needs to be shown that in this case 
$P_i^{(k)}(\nu,t)=P_i^{(k)}(\nu)$.
To this end it suffices
to show that there exists an index $j$ such that
\begin{equation}\label{j}
\lt+t_{j-1,k+1}\le \s^{(k)}+1<\lt+t_{j,k}.
\end{equation}
Assume that $p^{(k)}=j^{(k)}$, so that $\lt+t_{j^{(k)},k+1}=\s^{(k)}+1$.
By \eqref{jk} $\s^{(k)}+1\le \lt+\tb_{j^{(k)},k}
<\lt+\tb_{j^{(k)}+1,k}=\lt+t_{j^{(k)}+1,k}$ so that \eqref{j} holds with
$j=j^{(k)}+1$. Next assume that $p^{(k)}=j^{(k)}-1$. Then
by \eqref{p}, 
$\s^{(k)}+1<\lt+\tb_{j^{(k)},k-1}\le \lt+\tb_{j^{(k)},k}=\lt+t_{j^{(k)},k}$.
Furthermore by \eqref{jk}, 
$\lt+t_{j^{(k)}-1,k+1}=\lt+\tb_{j^{(k)}-1,k+1}\le \s^{(k)}+1$
which implies \eqref{j} with $j=j^{(k)}$.
In summary \eqref{j} holds for $j=p^{(k)}+1$.

\begin{figure}
\begin{picture}(300,60)(0,0)
\DashLine(0,10)(300,10){3}
\DashLine(0,40)(300,40){3}
\Line(10,50)(10,15)
\Line(10,15)(30,15)
\Line(30,15)(30,0)
\CCirc(10,15){2}{Black}{Black}
\CCirc(30,15){2}{Black}{White}
\Line(70,50)(70,20)
\Line(70,20)(120,20)
\Line(120,20)(120,0)
\CCirc(70,20){2}{Black}{Black}
\CCirc(120,20){2}{Black}{White}
\Line(100,50)(100,25)
\Line(100,25)(140,25)
\Line(140,25)(140,0)
\CCirc(100,25){2}{Black}{Black}
\CCirc(140,25){2}{Black}{White}
\Line(180,50)(180,30)
\Line(180,30)(210,30)
\Line(210,30)(210,0)
\CCirc(180,30){2}{Black}{Black}
\CCirc(210,30){2}{Black}{White}
\Line(250,50)(250,35)
\Line(250,35)(290,35)
\Line(290,35)(290,0)
\CCirc(250,35){2}{Black}{Black}
\CCirc(290,35){2}{Black}{White}
\SetColor{Blue}
\ZigZag(50,17.5)(70,17.5){2}{5}
\CCirc(50,17.5){2}{Blue}{White}
\CCirc(70,17.5){2}{Blue}{White}
\ZigZag(140,27.5)(180,27.5){2}{10}
\CCirc(140,27.5){2}{Blue}{Blue}
\CCirc(180,27.5){2}{Blue}{White}
\ZigZag(210,32.5)(230,32.5){2}{5}
\CCirc(210,32.5){2}{Blue}{Blue}
\CCirc(230,32.5){2}{Blue}{Blue}
\SetColor{Red}
\Line(50,0)(50,50)
\Line(230,0)(230,50)
\Text(50,-5)[]{$s^{(k)}$}
\Text(230,-5)[]{$s^{(k+1)}$}
\end{picture}
\caption{Illustration of the extra terms in \eqref{add1} and \eqref{add2}
\label{fig add}}
\end{figure}
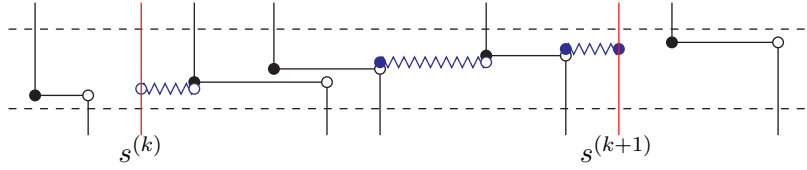

It remains to show that for $i\neq \s^{(k)}+1$ the bounds 
\eqref{bound inter} imply \eqref{bound} with $t$ as in \eqref{t def}.
First assume $1\le k<r$.
For $i$ such that $\lt+\tb_{j-1,k+1}\le i<\lt+\tb_{j,k}$
with $1\le j\le \la_r-\la_n$ \eqref{bound bar} simply reads
$\xb_i^{(k)}\le P_i^{(k)}(\nub)$. By construction there are no singular
strings of length $\s^{(k)}<i\le \s^{(k+1)}$ in $(\nub,\Jb)^{(k)}$.
Hence, for $1\le k\le r-2$ we can sharpen the bounds in \eqref{bound bar} 
and therefore also those in \eqref{bound inter} by adding the terms
\begin{equation}\label{add1}
-\chi(\s^{(k)}<i\le \min\{\s^{(k+1)},\lt+\tb_{j^{(k)},k}-1\})
\end{equation}
if $j^{(k)}=j^{(k+1)}$ and
\begin{multline}\label{add2}
-\chi(\s^{(k)}<i<\lt+\tb_{j^{(k)},k})
-\sum_{j=j^{(k)}}^{j^{(k+1)}-2}\chi(\lt+\tb_{j,k+1}\le i<\lt+\tb_{j+1,k})\\
-\chi(\lt+\tb_{j^{(k+1)}-1,k+1}\le i\le 
 \min\{\s^{(k+1)},\lt+\tb_{j^{(k+1)},k}-1\})
\end{multline}
if $j^{(k)}<j^{(k+1)}$.
In terms of the paths, this corresponds to adding a horizontal line segment
(which is equivalent to extra minus signs) in the $k$-th strip in the 
interval $\s^{(k)}<i\le \s^{(k+1)}$ whenever there is a horizontal gap 
between two neighboring paths. 
An example is given in Figure \ref{fig add}. It depicts the $k$-th
strip and the zigzag lines correspond to the added line segments.

The sum of extra terms \eqref{add1} or \eqref{add2} and 
$\chi(\s^{(k)}<i\le \s^{(k+1)})$ does not exceed
\begin{multline}\label{plus terms}
\sum_{j=j^{(k)}}^{p^{(k+1)}}
 \chi(\lt+\max\{\tb_{j-1,k+1},\tb_{j,k}\}\le i<\lt+\tb_{j,k+1})\\
=\sum_{j=j^{(k)}}^{p^{(k+1)}}\chi(\lt+t_{j,k+1}\le i<\lt +\tb_{j,k+1})
 -\chi(\s^{(k)}<i<\lt+\tb_{p^{(k)},k}).
\end{multline}
To obtain the first line of \eqref{plus terms} we have used
$\tb_{j^{(k)}-1,k+1}<\tb_{j^{(k)},k}$ by \eqref{jk} for the term
$j=j^{(k)}$, the definition \eqref{p} of $p^{(k+1)}$ and
$\s^{(k+1)}<\lt+\tb_{j^{(k+1)},k+1}$ when $p^{(k+1)}=j^{(k+1)}$
which follows from \eqref{jk}. When $p^{(k)}=j^{(k)}$ the second
line follows directly using \eqref{t def}.
When $p^{(k)}=j^{(k)}-1$ use that $\lt+\tb_{j^{(k)}-1,k}\le
\lt+\tb_{j^{(k)}-1,k+1}\le \s^{(k)}$ by \eqref{jk} so that the
last term vanishes.

Similarly for $k=r-1$ the bounds in \eqref{bound inter} can be sharpened 
by adding
\begin{equation*}
-\chi(\s^{(r-1)}<i<\lt+\tb_{j^{(r-1)},r-1})
-\sum_{j=j^{(r-1)}}^{\la_r-\la_n-1}\chi(\lt+\tb_{j,r}\le i<\lt+\tb_{j+1,r-1}).
\end{equation*}
Together with $\chi(\s^{(r-1)}<i\le \s^{(r)})=\chi(\s^{(r-1)}<i)$ 
this yields by similar reasons as before
\begin{multline}\label{plus terms r-1}
\sum_{j=j^{(r-1)}}^{\la_r-\la_n-1}
 \chi(\lt+t_{j,r}\le i<\lt +\tb_{j,r})
 +\chi(i\ge \lt+t_{\la_r-\la_n,r})\\
-\chi(\s^{(r-1)}<i<\lt+\tb_{p^{(r-1)},r-1}).
\end{multline}

Note that $\tb_{j-1,k}\le t_{j,k}\le \tb_{j,k}$ for 
$j^{(k-1)}\le j\le p^{(k)}$ and $1\le k<r$.
In addition $\s^{(k-1)}<\lt +t_{j^{(k-1)},k}$ for $1\le k\le r$.
For $k=1$ this follows from \eqref{j0} and Remark \ref{rem col}, and
for $1<k\le r$ and $p^{(k-1)}=j^{(k-1)}$ this follows from \eqref{t def} and
for $p^{(k-1)}=j^{(k-1)}-1$ one exploits the first condition of \eqref{p}
and \eqref{t jk}.
Also $\s^{(k)}\ge \lt+\tb_{j^{(k)}-1,k}$ 
for $1\le k<r$ thanks to \eqref{jk}. This implies for $1\le k<r$
\begin{multline}\label{minus terms}
-\chi(\s^{(k-1)}<i\le \s^{(k)})
\le -\sum_{j=j^{(k-1)}}^{p^{(k)}}\chi(\lt+t_{j,k}\le i<\lt+\tb_{j,k})\\
+\chi(\s^{(k)}<i<\lt+\tb_{p^{(k)},k}).
\end{multline}
{}From \eqref{plus terms} (or \eqref{plus terms r-1} for $k=r-1$) and
\eqref{minus terms} it is straightforward to see that \eqref{bound inter}
implies \eqref{bound} for $1\le k<r$.

Now consider $k=r$. Then $\s^{(r)}=\infty$,
$\s^{(r-1)}<\lt+t_{j^{(r-1)},r}$ as shown above \eqref{minus terms}
and $\tb_{j,r}\le t_{j+1,r}$ for $j^{(r-1)}\le j<\la_r-\la_n$ so that
\begin{equation*}
\begin{split}
&-\chi(\s^{(r-1)}<i\le \s^{(r)})
 -\sum_{j=1}^{\lab_r-\lab_n}\chi(i\ge \lt+\tb_{j,r})\\
\le & -\chi(i\ge \lt+t_{j^{(r-1)},r})
 -\sum_{j=1}^{j^{(r-1)}-1} \chi(i\ge \lt +t_{j,r})
 -\sum_{j=j^{(r-1)}}^{\la_r-\la_n-1} \chi(i\ge \lt+t_{j+1,r})\\
=&-\sum_{j=1}^{\la_r-\la_n} \chi(i\ge \lt+t_{j,r}).
\end{split}
\end{equation*}
For $r<k<n$ we have $\s^{(k-1)}=\s^{(k)}=\infty$ and 
$\tb_{j,k}=t_{j,k}$ so that 
$-\chi(\s^{(k-1)}<i\le \s^{(k)})=0$ and
$-\sum_{j=1}^{\lab_k-\lab_n}\chi(i\ge \ltb+\tb_{j,k})
=-\sum_{j=1}^{\la_k-\la_n}\chi(i\ge \lt+t_{j,k})$.

For $r\le k<n$ we have $\s^{(k)}=\s^{(k+1)}=\infty$
and $\tb_{j,k+1}=t_{j,k+1}$ so that $\chi(\s^{(k)}<i\le \s^{(k+1)})=0$
and $\sum_{j=1}^{\lab_{k+1}-\lab_n}\chi(i\ge \ltb+\tb_{j,k+1})
=\sum_{j=1}^{\la_{k+1}-\la_n}\chi(i\ge \lt+t_{j,k+1})$.

This concludes the proof that \eqref{bound inter} implies \eqref{bound}
for $1<r<n$.

\subsection*{Case $r=n$.}
In this case $\lab_k=\la_k$ for $1\le k<n$,
$\lab_n=\la_n-1$, $\ltb=\lt-1$ and $\tb$ is a tableau over the 
alphabet $\{1,2,\ldots,\la_1-\la_n+1\}$. It follows from \eqref{order j}
that $j^{(0)}=\cdots=j^{(n-1)}=1$. Note that this requires in
particular that $\s^{(0)}=\mu_L-1<\ltb+1=\lt$. Define
$t_{j,k}=\tb_{j+1,k}-1$ for $1\le k<n$ and $1\le j\le \la_k-\la_n$.
Then by the column-strictness of $\tb$ we have
$t_{j,k}<t_{j+1,k}$ and $t_{j,k}\le t_{j,k+1}$. Note in particular
that $t_{j,1}=\tb_{j+1,1}-1=j$ so that $t$ is a column-strict
tableau over the alphabet $\{1,2,\ldots,\la_1-\la_n\}$.
In addition it follows from \eqref{jk} that
$0\le \s^{(k)}+1\le \ltb+\tb_{1,k}<\ltb+\tb_{2,k}=\lt+t_{1,k}$
so that \eqref{j} holds for $j=1$. This ensures \eqref{bound}
for $i=\s^{(k)}+1$.
Using the fact that there are no singular strings of length 
$\s^{(k)}<i<\ltb+\tb_{1,k}$ in $(\nub,\Jb)^{(k)}$
and that $\s^{(k+1)}<\ltb+\tb_{1,k+1}$ by \eqref{jk} the term
$\chi(\s^{(k)}<i\le \s^{(k+1)})$ in \eqref{bound inter} can 
be safely replaced by $\chi(\ltb+\tb_{1,k}\le i<\ltb+\tb_{1,k+1})$.
Furthermore dropping the term $-\chi(\s^{(k-1)}<i\le \s^{(k)})$ 
equation \eqref{bound inter} becomes
\begin{equation*}
\x_i^{(k)}\le P_i^{(k)}(\nu)
-\sum_{j=2}^{\la_k-\la_n+1}\chi(i\ge \ltb+\tb_{j,k})
+\sum_{j=2}^{\la_{k+1}-\la_n+1}\chi(i\ge \ltb+\tb_{j,k+1}).
\end{equation*}
Using $\ltb=\lt-1$ and the definition of $t$ this is exactly \eqref{bound}.

This concludes the proof of the forward direction of the theorem.

\subsection*{Proof of the reverse direction}

Let us now prove the reverse direction. To this end consider a
rigged configuration $(\nu,J)$ corresponding to a column-strict tableau $T$
of shape $\la$ and content $\mu$ under $\psib_\mu$
which satisfies $\nu_1^{(k)}\le\ell$ for all $1\le k<n$ and \eqref{bound}. 
We need to show that then $T$ is of level $\ell$. This is equivalent to
showing that $T^-$ is of level $\ell$ and that $\la_1-\la_n^0\le \ell$
if $r=1$, where $r$ is the row index of the cell $\la/\lab$,
$\lab=\shape(T^-)$ and $\la^0=\shape(T^{-^{\mu_L}})$.
By induction the statement that $T^-$ is of level $\ell$ 
is equivalent to the statement that $\ltb=\ell-\lab_1+\lab_n\ge 0$ and
$(\nub,\Jb)=\db(\nu,J)$ satisfies $\nub_1^{(k)}\le \ell$ 
for all $1\le k<n$ and 
\begin{equation}\label{bound bar rev}
\xb_i^{(k)}\le P_i^{(k)}(\nub)
 -\sum_{j=1}^{\lab_k-\lab_n}\chi(i\ge \ltb+\tb_{j,k})
 +\sum_{j=1}^{\lab_{k+1}-\lab_n}\chi(i\ge \ltb+\tb_{j,k+1})
\end{equation}
for some column-strict tableau $\tb$ of shape $(\lab_1-\lab_n,\ldots,
\lab_{n-1}-\lab_n)^t$. To prove $\ltb\ge 0$ it suffices to show that 
$r=n$ cannot occur when $\lt=0$.

Let $\ell^{(k)}$ ($1\le k<r$) be the length of the selected singular
string in $(\nu,J)^{(k)}$ under $\db$. By definition
$\mu_L=\ell^{(0)}\le \ell^{(1)}\le \ell^{(2)}\le \cdots 
\le \ell^{(r-1)}\le \ell$,
and the rigged configuration $(\nub,\Jb)$ is obtained from $(\nu,J)$
by shortening the selected strings by one, making them singular again
and leaving all other strings unchanged.
Since $\nu_1^{(k)}\le \ell$ this immediately implies $\nub_1^{(k)}\le \ell$
for all $1\le k<n$.

The vacancy numbers are related by
\begin{equation}\label{vac ch}
P_i^{(k)}(\nub)=P_i^{(k)}(\nu)-\chi(\ell^{(k-1)}\le i<\ell^{(k)})
 +\chi(\ell^{(k)}\le i<\ell^{(k+1)}).
\end{equation}
Furthermore $\xb_i^{(k)}\le \x_i^{(k)}$ for $i\neq \ell^{(k)}-1$ and
$\xb_{\ell^{(k)}-1}^{(k)}=P_{\ell^{(k)}-1}^{(k)}(\nub)$ for $1\le k<r$
so that \eqref{bound} implies
\begin{multline}\label{bound rev}
\xb_i^{(k)}\le P_i^{(k)}(\nub)+\chi(\ell^{(k-1)}\le i<\ell^{(k)})
 -\chi(\ell^{(k)}\le i<\ell^{(k+1)})\\
 -\sum_{j=1}^{\la_k-\la_n}\chi(i\ge \lt+t_{j,k})
 +\sum_{j=1}^{\la_{k+1}-\la_n}\chi(i\ge \lt+t_{j,k+1})
\end{multline}
for $i\neq \ell^{(k)}-1$. 

Since $\ell^{(k)}$ is the length of a singular string in $(\nu,J)^{(k)}$
it must be in one of the intervals in \eqref{sing cond}. Let $j^{(k)}$
for $1\le k<r$ be the index such that
\begin{equation}\label{jk rev}
\lt+t_{j^{(k)}-1,k+1}\le \ell^{(k)} <\lt+t_{j^{(k)},k}
\end{equation}
where recall that $t_{0,k+1}=-\lt$ and $t_{j,k}=\la_1-\la_n+1$
for all $j>\la_k-\la_n$.
By similar arguments as in the derivation of \eqref{order j} one finds
\begin{equation}\label{jk order rev}
1\le j^{(1)}\le \cdots\le j^{(r-1)}\le \la_r-\la_n+1.
\end{equation}

\subsection*{Case $r=1$.}
In this case $\lab_k=\la_k$ for $1<k\le n$,
$\lab_1=\la_1-1$, $\lt=\ltb-1$
and $t_{j,k}\in\{1,2,\ldots,\la_1-\la_n\}=\{1,2,\ldots,\lab_1-\lab_n+1\}$.
Let $a^{(k)}$ ($1\le k<n$) be maximal such that $t_{j,k}=j$ for all 
$1\le j\le a^{(k)}$. It follows from Remark \ref{rem col} and 
$t_{j,k}\le t_{j,k+1}$ that $0\le a^{(n-1)}\le \cdots\le a^{(2)}\le a^{(1)}
=\la_1-\la_n$. Set $\tb_{j,1}=j$ for $1\le j\le \lab_1-\lab_n$ and 
\begin{equation*}
\tb_{j,k}=\begin{cases}
 j & \text{for $1\le j\le a^{(k)}$,}\\
 t_{j,k}-1 & \text{for $a^{(k)}<j\le \lab_k-\lab_n$,}
\end{cases}
\end{equation*}
for $1<k<n$.
The definition of $a^{(k)}$ and column-strictness of $t$ ensure 
the column-strictness of $\tb$. Note that the terms $j=1,2,\ldots,a^{(k+1)}$
in the two sums in \eqref{bound rev} cancel each other.
Recall that $\ell^{(0)}=\mu_L$ and $\ell^{(k)}=\infty$ for $k\ge 1$.
Assume $k=1$. The term $\chi(\ell^{(0)}\le i<\ell^{(1)})=\chi(i\ge \mu_L)$ 
in \eqref{bound rev} can be replaced by $\chi(i\ge \lt+t_{a^{(2)}+1,1})
=\chi(i\ge \lt+a^{(2)}+1)$. If $\mu_L\le \lt+a^{(2)}+1$ this follows from
the fact that by construction there are no singular strings 
of length $\ge \mu_L$ in $(\nu,J)^{(1)}$. For $\mu_L> \lt+a^{(2)}+1$ 
we have $\chi(i\ge \mu_L)\le \chi(i\ge \lt+a^{(2)}+1)$. 
Hence using $\lt=\ltb-1$
\begin{equation*}
\begin{split}
\xb_i^{(1)}\le & P_i^{(1)}(\nub)+\chi(i\ge \ltb+a^{(2)})
-\sum_{j=a^{(2)}+1}^{\lab_1-\lab_n+1}\chi(i\ge \lt+j)\\
&+\sum_{j=a^{(2)}+1}^{\lab_{2}-\lab_n}\chi(i\ge \lt+t_{j,2})\\
=& P_i^{(1)}(\nub)-\sum_{j=a^{(2)}+1}^{\lab_1-\lab_n}\chi(i\ge \ltb+j)
+\sum_{j=a^{(2)}+1}^{\lab_{2}-\lab_n}\chi(i\ge \ltb+\tb_{j,2})
\end{split}
\end{equation*}
which is \eqref{bound bar rev} for $k=1$. Now assume $1<k<n$.
Since $\ell^{(k)}=\infty$ for $1\le k\le n$ the terms involving $\ell^{(k)}$
in \eqref{bound rev} vanish and
\begin{equation*}
\begin{split}
\xb_i^{(k)}\le& P_i^{(k)}(\nub)
-\sum_{j=a^{(k+1)}+1}^{\lab_k-\lab_n}\chi(i\ge \lt+t_{j,k})
+\sum_{j=a^{(k+1)}+1}^{\lab_{k+1}-\lab_n}\chi(i\ge \lt+t_{j,k+1})\\
\le & P_i^{(k)}(\nub)
-\sum_{j=a^{(k+1)}+1}^{\lab_k-\lab_n}\chi(i\ge \ltb+\tb_{j,k})
+\sum_{j=a^{(k+1)}+1}^{\lab_{k+1}-\lab_n}\chi(i\ge \ltb+\tb_{j,k+1})
\end{split}
\end{equation*}
which is \eqref{bound bar rev} for $1<k<n$. This concludes the proof
that \eqref{bound} implies \eqref{bound bar rev} for $r=1$.

\subsection*{Case $1<r<n$.}
Here $\lab_k=\la_k$ for $k\neq r$, $\lab_r=\la_r-1$
and $\ltb=\lt$.
Set $p^{(r)}=j^{(r)}=\la_r-\la_n$ and
\begin{equation}\label{p rev}
p^{(k)}=\begin{cases}
 j^{(k)}-1 & \text{for $\ell^{(k)}\le \lt+t_{j^{(k)}-1,k+2}$,}\\
 \min\{j^{(k)},p^{(k+1)}\} & \text{for $\ell^{(k)}> \lt+t_{j^{(k)}-1,k+2}$,}
\end{cases}
\end{equation}
for $1\le k<r$ where recall that 
$t_{j,r+1}=\la_1-\la_n+1$ for $j>\la_{r+1}-\la_n$. 
Note that $p^{(k)}=j^{(k)}$ or $j^{(k)}-1$ and $p^{(k)}\le p^{(k+1)}$
due to \eqref{jk order rev}.
Define $\tb_{j,1}=j$ for $1\le j\le \lab_1-\lab_n$,
\begin{equation}\label{tbar}
\tb_{j,k}=\begin{cases}
 t_{j,k} & \text{for $1\le j<p^{(k-1)}$ and 
  $j^{(k)}\le j\le \lab_k-\lab_n$,}\\
 \min\{t_{j,k+1},t_{j+1,k}\} & \text{for $p^{(k-1)}\le j<p^{(k)}$,}\\
 \ell^{(k)}-\lt & \text{for $j=p^{(k)}=j^{(k)}-1$,}
\end{cases}
\end{equation}
for $1<k\le r$ and $\tb_{j,k}=t_{j,k}$ for $r<k<n$ and 
$1\le j\le \lab_k-\lab_n$. 
Recall that $t_{j,k}=\la_1-\la_n+1$ for $j>\la_k-\la_n$.

It needs to be shown that $\tb$ is a tableau over the
alphabet $\{1,2,\ldots,\lab_1-\lab_n\}$. Since
$t_{j,k}\in\{1,2,\ldots,\la_1-\la_n\}=\{1,2,\ldots,\lab_1-\lab_n\}$
the only problematic case is the third case in \eqref{tbar}.
Condition \ref{c1} of the theorem implies that $\ell^{(k)}\le \ell$ 
for $1\le k<r$ so that $\ell^{(k)}-\lt\le \lab_1-\lab_n$. By \eqref{jk rev}
the condition $1\le \ell^{(k)}-\lt$ can only be violated if $j^{(k)}=1$.
Assume that $j^{(k)}=1$ for some $1\le k<r$.
Let $h$ be maximal such that $j^{(k)}=1$ for all $1\le k\le h$.
Then $\ell^{(k)}>\lt+t_{j^{(k)}-1,k+2}=0$ for all $1\le k\le h$ so that
the second case in \eqref{p rev} applies. 
If $h<r-1$ we have $p^{(h+1)}\ge j^{(h+1)}-1\ge 1$ by the maximality of $h$.
If $h=r-1$, $p^{(r)}=j^{(r)}=\la_r-\la_n\ge 1$. In both cases it follows that
$p^{(k)}=j^{(k)}=1$ for all $1\le k\le h$. Hence by \eqref{tbar} the case
$\tb_{p^{(k)},k}=\ell^{(k)}-\lt<1$ does not occur.
This proves that $\tb_{j,k}\in\{1,2,\ldots,\lab_1-\lab_n\}$.

It remains to show that $\tb$ is column-strict. The condition
$\tb_{j,k}<\tb_{j+1,k}$ only needs to be considered for
$p^{(k-1)}-1\le j<j^{(k)}$ and $1<k<r$ and for 
$p^{(r-1)}-1\le j\le \la_{r+1}-\la_n$ and $k=r$ by the column-strictness 
of $t$. In these cases $\tb_{j,k}<\tb_{j+1,k}$ can be deduced from the 
following inequalities:
\begin{align*}
(a)\quad & t_{j,k}<\min\{t_{j+1,k+1},t_{j+2,k}\},\\
(b)\quad & \min\{t_{j,k+1},t_{j+1,k}\}\le t_{j+1,k}<t_{j+2,k},\\
   \quad & \min\{t_{j,k+1},t_{j+1,k}\}\le t_{j,k+1}<t_{j+1,k+1},\\
(c)\quad & \min\{t_{j^{(k)}-2,k+1},t_{j^{(k)}-1,k}\}\le t_{j^{(k)}-2,k+1}
       <t_{j^{(k)}-1,k+1}\le \ell^{(k)}-\lt,\\
   \quad & \ell^{(k)}-\lt<t_{j^{(k)},k} \qquad \text{for $1\le k<r$},\\
(d)\quad & \min\{t_{j^{(k)}-1,k+1},t_{j^{(k)},k}\}=
       t_{j^{(k)}-1,k+1}<t_{j^{(k)},k} \qquad \text{for $1\le k<r$},
\end{align*}
where \eqref{jk rev} was employed extensively. 
The condition $\tb_{j,k}\le \tb_{j,k+1}$ needs to be verified for
$k=1$ and $p^{(1)}\le j<j^{(2)}$, for $1<k<r$ and $p^{(k-1)}\le j<j^{(k+1)}$
and for $k=r$ and $p^{(r-1)}\le j\le \la_{r+1}-\la_n$. 
In these cases $\tb_{j,k}\le \tb_{j,k+1}$ can be deduced from the 
following inequalities:
\begin{align*}
(a) \quad & \min\{t_{j,k+1},t_{j+1,k}\}\le t_{j,k+1},\\
(b) \quad & t_{j,k}\le \min\{t_{j,k+2},t_{j+1,k+1}\},\\
(c) \quad & t_{j^{(k+1)}-1,k}\le t_{j^{(k+1)}-1,k+2}\le \ell^{(k+1)}-\lt,
\end{align*}
where again \eqref{jk rev} was employed. In addition for $1<k<r$ we have 
$\ell^{(k)}-\lt\le \min\{t_{j^{(k)}-1,k+2},t_{j^{(k)},k+1}\}$
if $\ell^{(k)}\le \lt+t_{j^{(k)}-1,k+2}$. 
If $\ell^{(k)}>\lt+t_{j^{(k)}-1,k+2}$ then $p^{(k)}=j^{(k)}-1$
is only possible if $p^{(k)}=p^{(k+1)}=j^{(k)}-1$ which implies that 
$j^{(k)}=j^{(k+1)}$. However in this case
$\tb_{j^{(k)}-1,k}=\ell^{(k)}-\lt\le \ell^{(k+1)}-\lt=\tb_{j^{(k)}-1,k+1}$.
This proves the column-strictness of $\tb$.

By definition $\xb_{\ell^{(k)}-1}^{(k)}=P_{\ell^{(k)}-1}^{(k)}(\nub)$
for $1\le k<r$. Hence we need to check that 
$P_i^{(k)}(\nub,\tb)=P_i^{(k)}(\nub)$ for $i=\ell^{(k)}-1$.
It suffices to show that there exists an index $j$ such that
\begin{equation}\label{ccond}
\lt+\tb_{j-1,k+1}\le \ell^{(k)}-1<\lt+\tb_{j,k}.
\end{equation}
Assume that $p^{(k)}=j^{(k)}-1$. Then $\lt+\tb_{j^{(k)}-1,k}=\ell^{(k)}$
and by \eqref{jk rev}
$\ell^{(k)}\ge \lt+t_{j^{(k)}-1,k+1}>\lt+t_{j^{(k)}-2,k+1}
=\lt+\tb_{j^{(k)}-2,k+1}$ so that \eqref{ccond} holds for $j=j^{(k)}-1$.
Now assume $p^{(k)}=j^{(k)}$. Then by \eqref{p rev},
$\ell^{(k)}>\lt+t_{j^{(k)}-1,k+2}\ge \lt+t_{j^{(k)}-1,k+1}=
\lt+\tb_{j^{(k)}-1,k+1}$.
Furthermore by \eqref{jk rev} $\ell^{(k)}<\lt+t_{j^{(k)},k}
=\lt+\tb_{j^{(k)},k}$ so that \eqref{ccond} holds for $j=j^{(k)}$.

For $i\neq \ell^{(k)}-1$ we need to show that \eqref{bound rev} implies
\eqref{bound bar rev}. Note that for $1\le k<r$ we have
$\tb_{j,k+1}\le t_{j+1,k+1}$ for $p^{(k)}\le j<j^{(k+1)}$
since $\min\{t_{j,k+2},t_{j+1,k+1}\}\le t_{j+1,k+1}$ and
$\ell^{(k+1)}-\lt<t_{j^{(k+1)},k+1}$ for $1\le k<r-1$ by \eqref{jk rev}.
In addition $\ell^{(k+1)}\ge \lt+\tb_{j^{(k+1)}-1,k+1}$. For
$p^{(k+1)}=j^{(k+1)}-1$ this follows directly from \eqref{tbar},
and for $p^{(k+1)}=j^{(k+1)}$ we have $\ell^{(k+1)}\ge 
\lt+t_{j^{(k+1)}-1,k+2}\ge \lt+\tb_{j^{(k+1)}-1,k+1}$ by \eqref{jk rev}.
Since furthermore $\ell^{(k)}<\lt+t_{j^{(k)},k+1}$ by 
\eqref{jk rev} we have for $1\le k<r$
\begin{multline}\label{b min}
-\chi(\ell^{(k)}\le i<\ell^{(k+1)})
\le \chi(\lt+t_{p^{(k)},k+1}\le i<\ell^{(k)})\\
-\sum_{j=p^{(k)}}^{j^{(k+1)}-1}\chi(\lt+t_{j,k+1}\le i<\lt+\tb_{j,k+1})
-\delta_{k+1,r}\chi(i\ge \lt+t_{\la_r-\la_n,r}),
\end{multline}
where the last term occurs since $\ell^{(r)}=\infty$.
Observe that $t_{j,k+1}\le \tb_{j,k+1}$ for $p^{(k)}\le j<j^{(k+1)}$ and
$1\le k<r$ since $\min\{t_{j,k+2},t_{j+1,k+1}\}\ge t_{j,k+1}$ and
$\ell^{(k+1)}-\lt\ge t_{j^{(k+1)}-1,k+2}\ge t_{j^{(k+1)}-1,k+1}$ 
for $1\le k<r-1$ by \eqref{jk rev}. Hence using \eqref{b min} and
\eqref{tbar} we have for $1\le k<r$
\begin{equation}\label{bb min}
\begin{split}
&-\chi(\ell^{(k)}\le i<\ell^{(k+1)})
+\sum_{j=1}^{\la_{k+1}-\la_n}\chi(i\ge \lt+t_{j,k+1})\\
\le &\sum_{j=1}^{\lab_{k+1}-\lab_n}\chi(i\ge \lt+\tb_{j,k+1})
 +\chi(\lt+t_{p^{(k)},k+1}\le i<\ell^{(k)}).
\end{split}
\end{equation}

Since $t_{j,1}=j$ by Remark \ref{rem col} and 
$t_{j,1}\le t_{j,2}$ equation \eqref{jk rev} implies that
either $\ell^{(0)}\le \ell^{(1)}<\lt+1$ for $j^{(1)}=1$ or 
$\ell^{(1)}=\lt+j^{(1)}-1$ and $t_{j,2}=j$ for 
$1\le j<j^{(1)}\le \la_r-\la_n+1$. Note that in both cases
\eqref{bound} reads $\x_i^{(1)}\le P_i^{(1)}(\nu)$ for
$1\le i<\ell^{(1)}$. Since by construction there are no singular strings of 
length $\ell^{(0)}\le i<\ell^{(1)}$ in $(\nu,J)^{(1)}$ we can add the 
term $-\chi(\ell^{(0)}\le i<\ell^{(1)})$ to \eqref{bound} for $k=1$. 
This has the effect that the term 
$\chi(\ell^{(0)}\le i<\ell^{(1)})$ in \eqref{bound rev} for $k=1$ can be 
dropped. Note that for $j^{(1)}=1$ we have $p^{(1)}=1$ so that the term
$\chi(\lt+t_{p^{(1)},2}\le i<\ell^{(1)})$ in \eqref{bb min} is zero.
For $j^{(1)}>1$ this term is also zero since 
$t_{p^{(1)},2}\ge t_{j^{(1)}-1,2}=j^{(1)}-1$ and $\ell^{(1)}=\lt+j^{(1)}-1$.
Since in addition $-\sum_{j=1}^{\la_1-\la_n}\chi(i\ge \lt+t_{j,1})
=-\sum_{j=1}^{\lab_1-\lab_n}\chi(i\ge \lt+\tb_{j,1})$,
this proves that \eqref{bound rev} implies \eqref{bound bar rev} for $k=1$.

Now assume that $1<k<r$. By construction there are no
singular strings of length $\ell^{(k-1)}\le i<\ell^{(k)}$ in
$(\nu,J)^{(k)}$. Therefore the bounds \eqref{bound} and hence also
the bounds in \eqref{bound rev} can be sharpened by adding
\begin{equation*}
-\chi(\max\{\ell^{(k-1)},\lt+t_{j^{(k)}-1,k+1}\}\le i<\ell^{(k)})
\end{equation*}
for $j^{(k-1)}=j^{(k)}$ and
\begin{multline*}
-\chi(\max\{\ell^{(k-1)},\lt+t_{j^{(k-1)}-1,k+1}\}\le i<\lt+t_{j^{(k-1)},k})\\
-\sum_{j=j^{(k-1)}+1}^{j^{(k)}-1}\chi(\lt+t_{j-1,k+1}\le i<\lt+t_{j,k})
-\chi(\lt+t_{j^{(k)}-1,k+1}\le i<\ell^{(k)})
\end{multline*}
for $j^{(k-1)}<j^{(k)}$. Adding these to 
$\chi(\ell^{(k-1)}\le i<\ell^{(k)})$ does not exceed
\begin{multline*}
\sum_{j=p^{(k-1)}}^{j^{(k)}-1}
 \chi(\lt+t_{j,k}\le i<\lt+\min\{t_{j,k+1},t_{j+1,k}\})\\
=\sum_{j=p^{(k-1)}}^{j^{(k)}-1}\chi(\lt+t_{j,k}\le i<\lt+\tb_{j,k})
 -\chi(\lt+t_{p^{(k)},k+1}\le i<\ell^{(k)}).
\end{multline*}
Using again that $t_{j,k}\le \tb_{j,k}$ for $p^{(k-1)}\le j<j^{(k)}$
this can be combined with the term 
$ -\sum_{j=1}^{\la_k-\la_n}\chi(i\ge \lt+t_{j,k})$ of \eqref{bound rev}
to yield
$-\sum_{j=1}^{\lab_k-\lab_n}\chi(i\ge \lt+\tb_{j,k})
-\chi(\lt+t_{p^{(k)},k+1}\le i<\ell^{(k)})$.
Together with \eqref{bb min} this proves that \eqref{bound rev} implies
\eqref{bound bar rev} for $1<k<r$.

Consider $k=r$. Recall that $\la_{r+1}<\la_r$ so that \eqref{bound} 
implies $\x_i^{(r)}\le P_i^{(r)}(\nu)-1$ for 
$i\ge \lt+t_{\la_{r+1}-\la_n+1,r}$.
By construction there are no singular strings of length $i\ge \ell^{(r-1)}$ 
in $(\nu,J)^{(r)}$. Hence for $j^{(r-1)}\le \la_{r+1}-\la_n+1$
the bounds in \eqref{bound} and \eqref{bound rev} can be sharpened by adding
\begin{multline*}
-\chi(\max\{\ell^{(r-1)},\lt+t_{j^{(r-1)}-1,r+1}\}\le i<\lt+t_{j^{(r-1)},r})\\
-\sum_{j=j^{(r-1)}+1}^{\la_{r+1}-\la_n+1}\chi(\lt+t_{j-1,r+1}\le i<\lt+t_{j,r})
\end{multline*}
which added to $\chi(i\ge \ell^{(r-1)})$ does not exceed
\begin{equation}\label{com r}
\begin{split}
&\sum_{j=p^{(r-1)}}^{\la_{r+1}-\la_n}
 \chi(\lt+t_{j,r}\le i<\lt+\min\{t_{j,r+1},t_{j+1,r}\})
 +\chi(i\ge \lt+t_{\la_{r+1}-\la_n+1,r})\\
=&\sum_{j=p^{(r-1)}}^{\la_r-\la_n-1}
 \chi(\lt+t_{j,r}\le i<\lt+\min\{t_{j,r+1},t_{j+1,r}\})
 +\chi(i\ge \lt+t_{\la_r-\la_n,r})
\end{split}
\end{equation}
where in the last line we used that $\min\{t_{j,r+1},t_{j+1,r}\}=
t_{j+1,r}$ for $j>\la_{r+1}-\la_n$ 
since by definition $t_{j,r+1}=\la_1-\la_n+1$ in this case. 
The last line of \eqref{com r}
also makes sense for $j^{(r-1)}>\la_{r+1}-\la_n+1$ since then
$p^{(r-1)}=j^{(r-1)}-1$ and $\lt+t_{j^{(r-1)}-1,r}\le \ell^{(r-1)}$
by \eqref{jk rev}. The last line of \eqref{com r} combined with
$-\sum_{j=1}^{\la_r-\la_n}\chi(i\ge \lt+t_{j,r})$ yields
$-\sum_{j=1}^{\lab_r-\lab_n}\chi(i\ge \lt+\tb_{j,r})$ using \eqref{tbar}.
For $r<k<n$ the term $\chi(\ell^{(k-1)}\le i<\ell^{(k)})$
vanishes and $-\sum_{j=1}^{\la_k-\la_n}\chi(i\ge \lt+t_{j,k})
=-\sum_{j=1}^{\lab_k-\lab_n}\chi(i\ge \lt+\tb_{j,k})$. 
Similarly $-\chi(\ell^{(k)}\le i<\ell^{(k+1)})$ is zero for $r\le k<n$
and $\sum_{j=1}^{\la_{k+1}-\la_n}\chi(i\ge \lt+t_{j,k+1})
=\sum_{j=1}^{\lab_{k+1}-\lab_n}\chi(i\ge \lt+\tb_{j,k+1})$.
Together these results prove \eqref{bound bar rev} for $r\le k<n$.

This concludes the proof of the reverse direction of the theorem
for $1<r<n$.

\subsection*{Case $r=n$.}
In this case $\lab_k=\la_k$ for $1\le k<n$, $\lab_n=\la_n-1$
and $\ltb=\lt-1$. Then by \eqref{jk order rev} it follows that 
$j^{(1)}=\cdots=j^{(n-1)}=1$. In particular from \eqref{jk rev}, 
$\ell^{(1)}<\lt+t_{1,1}=\lt+1$ which yields a contradiction
when $\lt=0$ since by assumption $\ell^{(1)}\ge \mu_L\ge 1$. Hence the case
$r=n$ cannot occur when $\lt=0$. Define $\tb_{1,k}=1$
and $\tb_{j,k}=t_{j-1,k}+1$ for $1<j\le \lab_k-\lab_n$ for all
$1\le k<n$. Since $t_{j,k}\in\{1,2,\ldots,\lab_1-\lab_n-1\}$
it follows that $\tb_{j,k}\in\{1,2,\ldots,\lab_1-\lab_n\}$.
The column-strictness of $t$ immediately implies the column-strictness
of $\tb$.

Since $\mu_L\le \ell^{(k)}<\lt+t_{1,k}$ and there are no singular strings
of length $\ell^{(k-1)}\le i<\ell^{(k)}$ in $(\nu,J)^{(k)}$ we may drop
the term $\chi(\ell^{(k-1)}\le i<\ell^{(k)})$ in \eqref{bound rev}.
In addition dropping the term $-\chi(\ell^{(k)}\le i<\ell^{(k+1)})$
\eqref{bound rev} implies for $i\neq \ell^{(k)}-1$
\begin{equation*}
\begin{split}
\xb_i^{(k)}\le & P_i^{(k)}(\nub)
 -\sum_{j=1}^{\lab_k-\lab_n-1}\chi(i\ge \ltb+1+t_{j,k})
 +\sum_{j=1}^{\lab_{k+1}-\lab_n-1}\chi(i\ge \ltb+1+t_{j,k+1})\\
=& P_i^{(k)}(\nub)
 -\sum_{j=2}^{\lab_k-\lab_n}\chi(i\ge \ltb+\tb_{j,k})
 +\sum_{j=2}^{\lab_{k+1}-\lab_n}\chi(i\ge \ltb+\tb_{j,k+1}).
\end{split}
\end{equation*}
The terms $j=1$ can be added to both sums since they just cancel each
other so that we have \eqref{bound bar rev} for $i\neq \ell^{(k)}-1$.

Finally consider the case $i=\ell^{(k)}-1$.
We have $\ell^{(k)}<\lt+t_{1,k}$ for $1\le k<n$ so that
$\ell^{(k)}-1<\ltb+\tb_{2,k}$. Since the terms $j=1$ in the two
sums cancel, \eqref{bound bar rev} for $i=\ell^{(k)}-1$ reduces to
$\xb_{\ell^{(k)}-1}^{(k)}\le P_{\ell^{(k)}-1}^{(k)}(\nub)$,
or equivalently $P_{\ell^{(k)}-1}^{(k)}(\nub,\tb)=
P_{\ell^{(k)}-1}^{(k)}(\nub)$ as desired.

This concludes the proof that $T^-$ is of level $\ell$. 

\subsection*{Zu guter Letzt}
It remains to show that $\la_1-\la_n^0\le \ell$.

Define $T^b=T^{-^{\mu_L-b}}$ for $0\le b\le \mu_L$ with corresponding 
rigged configurations $(\nu^b,J^b)=\db^{\mu_L-b}(\nu,J)$, and
$\la^b=\shape(T^b)$. Let $(x^b)^{(k)}_i$ be the largest rigging
occurring for the strings of length $i$ in $(\nu^b,J^b)$ and let
$r_b$ be the row index of the cell
$\la^b/\la^{b-1}$ for $1\le b\le \mu_L$. 
Then $n\ge r_1\ge r_2\ge \cdots\ge r_{\mu_L}\ge 1$.
Denote the length of the selected string in $(\nu^b,J^b)^{(k)}$ under $\db$ 
by $\ell_b^{(k)}$. Let $1\le \beta\le \mu_L$ be maximal such that $r_\beta=n$.
If no such $\beta$ exists set $\beta=0$. Then
$\la_n^0=\la_n-\beta$. Hence proving $\la_1-\la_n^0\le \ell$ is
equivalent to showing that $\beta\le \lt$.

If $\mu_L\le \lt$ then also $\beta\le \mu_L\le \lt$. Hence assume
that $\mu_L=\lt+d$ with $d\ge 1$. For $b>\lt$ set $\bb=b-\lt$.
We will show by descending induction on $\lt<b\le \lt+d$ that
\begin{multline}\label{b ind}
(\x^b)_i^{(k)}\le P_i^{(k)}(\nu^b)\\
 -\sum_{j=1}^{\min\{\bb,\la_k-\la_n\}}\chi(i\ge \lt+t_{j,k})
 +\sum_{j=1}^{\min\{\bb,\la_{k+1}-\la_n\}}\chi(i\ge \lt+t_{j,k+1})
\end{multline}
for all $1\le k<n$ and $1\le i< \lt+t_{\bb,k+1}$ 
where recall that $t_{j,k}=\la_1-\la_n+1$ if $j>\la_k-\la_n$.

Let us first show that \eqref{b ind} holds for $b=\lt+d$. This follows
directly from \eqref{bound} using that
$-\sum_{j=1}^{\la_k-\la_n}\chi(i\ge \lt+t_{j,k})\le
-\sum_{j=1}^{\min\{d,\la_k-\la_n\}}\chi(i\ge \lt+t_{j,k})$
and $\sum_{j=1}^{\la_{k+1}-\la_n}\chi(i\ge \lt+t_{j,k+1})=
\sum_{j=1}^{\min\{d,\la_{k+1}-\la_n\}}\chi(i\ge \lt+t_{j,k+1})$
thanks to the fact that by assumption $1\le i<\lt+t_{d,k+1}$ and 
$t_{j,k+1}<t_{j+1,k+1}$.

Now assume \eqref{b ind} to be true for some $\lt<b\le \mu_L$.
We will prove that then $r_b<n$ and that \eqref{b ind} holds for $b-1$
if $b>\lt+1$. We claim that
\begin{equation}\label{lb cond}
\ell_b^{(k)}\ge \lt+t_{\bb,k+1}
\end{equation}
for all $1\le k<n$ where by definition $t_{j,n}=\la_1-\la_n+1$. Assume 
the opposite, namely let $1\le\kappa<n$ be the smallest index such that
$\ell_b^{(\kappa)}<\lt+t_{\bb,\kappa+1}$. By \eqref{b ind} there are
no singular strings of lengths $\lt+t_{\bb,\kappa}\le i<\lt+t_{\bb,\kappa+1}$
in $(\nu^b,J^b)^{(\kappa)}$ so that $\ell_b^{(\kappa)}<\lt+t_{\bb,\kappa}$.
By the minimality of $\kappa$ and the fact that $\ell_b^{(0)}=b=\lt+t_{\bb,1}$
this implies that $\ell_b^{(\kappa)}<\ell_b^{(\kappa-1)}$ which is a 
contradiction. This proves \eqref{lb cond}.
Note that similar to Remark \ref{rem pic} equation \eqref{b ind} can be 
interpreted in terms of $\bb$ non-intersecting paths which all end at 
position $\ell$. In this language the condition \eqref{lb cond}
states that $\ell_b^{(k)}$ is to the right of the $\bb$-th path.
Since all paths end at $\ell$ and there are no parts of length
greater than $\ell$ in $\nu_b^{(k)}$ this implies that $r_b<n$.
More precisely, $r_b\le k$ if $\bb>\la_{k+1}-\la_n$ for all $1\le k<n$.

Let us now prove \eqref{b ind} at $b-1$. It follows from \eqref{vac ch} that
\begin{equation*}
P_i^{(k)}(\nu^b)=P_i^{(k)}(\nu^{b-1})+\chi(\ell_b^{(k-1)}\le i<\ell_b^{(k)})
-\chi(\ell_b^{(k)}\le i<\ell_b^{(k+1)}).
\end{equation*}
By \eqref{lb cond}, $\ell_b^{(k-1)}\ge \lt+t_{\bb,k}$. 
Hence $\chi(\ell_b^{(k-1)}\le i<\ell_b^{(k)})\le \chi(\lt+t_{\bb,k}\le i)$ 
so that for $\bb\le \la_k-\la_n$
\begin{multline*}
\chi(\ell_b^{(k-1)}\le i<\ell_b^{(k)})
 -\sum_{j=1}^{\min\{\bb,\la_k-\la_n\}}\chi(i\ge \lt+t_{j,k})\\
\le -\sum_{j=1}^{\min\{\bb-1,\la_k-\la_n\}}\chi(i\ge \lt+t_{j,k})
\end{multline*}
as desired. When $\bb>\la_k-\la_n$ then $\ell^{(k-1)}_b=\infty$
so that the above inequality still holds. Since we only consider 
$1\le i<\lt+t_{\bb-1,k+1}$ the term $-\chi(\ell_b^{(k)}\le i<\ell_b^{(k+1)})$
does not contribute by \eqref{lb cond} and in addition $\bb$ can be replaced by
$\bb-1$ in $\sum_{j=1}^{\min\{\bb,\la_{k+1}-\la_n\}}\chi(i\ge \lt+t_{j,k+1})$.
This proves \eqref{b ind} for $b-1$. Since we have shown that \eqref{b ind}
implies that $r_b<n$ for $\lt<b\le \mu_L$ it follows that $\beta\le \lt$.

This concludes the proof of Theorem \ref{thm bound}.

\end{document}